\documentclass[12pt,a4paper]{article} 

\usepackage[T1]{fontenc}
\usepackage[latin1]{inputenc}
\usepackage{amssymb}
\usepackage{latexsym}
\usepackage{amsmath}
\usepackage{mathtools}
\usepackage{amscd}
\usepackage{graphicx}

\usepackage{proof}

\usepackage[hmargin=3cm,vmargin=3.5cm]{geometry}

\newtheorem{theorem}{Theorem}[section]
\newtheorem{corollary}[theorem]{Corollary}
\newtheorem{proposition}[theorem]{Proposition}
\newtheorem{construction}[theorem]{Construction}

\newtheorem{definition}[theorem]{Definition}

\newtheorem{prop}[theorem]{Proposition}
\newtheorem{lemma}[theorem]{Lemma}

\newtheorem{example}[theorem]{Example}

\newtheorem{remark}[theorem]{Remark}

\newtheorem{innercustomthm}{Theorem}
\newenvironment{customthm}[1]
  {\renewcommand\theinnercustomthm{#1}\innercustomthm}
  {\endinnercustomthm}

\newtheorem{innercustomlemma}{Lemma}
\newenvironment{customlemma}[1]
  {\renewcommand\theinnercustomlemma{#1}\innercustomlemma}
  {\endinnercustomlemma}

\newcommand{\qed}[0]{\Box}

\newcommand{\dontshow}[1]{}

\input diagxy.tex
\xyoption{curve}

\begin{document}

\newcommand{\sbrac}[1]{[\![#1]\!]}
\newcommand{\Hom}{{\rm Hom}}
\newcommand{\True}{\top}
\newcommand{\False}{\bot}
\newcommand{\drule}[2]{\frac{#1}{#2}}
\newcommand{\pgate}{\vdash}
\newcommand{\sequent}{\Longrightarrow}
\newcommand{\sequentunder}[1]{\xRightarrow{\;\;#1\;\;}}

\newcommand{\nat}{{\mathbb N}}
\newcommand{\nattype}[0]{{\rm N}}
\newcommand{\context}{{\rm context}}
\newcommand{\type}{{\rm type}}
\newcommand{\form}{{\rm formula}}

\newcommand{\Ty}{{\rm Ty}}
\newcommand{\Tm}{{\rm Tm}}

\newcommand{\vf}[0]{{\varphi}}
\newcommand{\pto}[0]{\rightharpoondown}
\newcommand{\piev}[0]{{\rm ev}}
\newcommand{\memof}[0]{\, \epsilon \,}
\newcommand{\subseteqof}[0]{\, \dot{\subseteq} \,}

\newcommand{\mono}[0]{\to/ >->/}

\newcommand{\emptycontext}{\langle\rangle}

\title{Categories with families and first-order logic with dependent sorts\footnote{A longer preprint version of this paper appeared under the title {\em Categories with families, FOLDS, and logic enriched type theory} (Palmgren 2016).}}
\author{Erik Palmgren \\ Department of Mathematics, Stockholm University}
\date{Final version June 29, 2019}

\maketitle

\begin{abstract} 
First-order logic with dependent sorts, such as  Makkai's first-order logic with 
dependent sorts (FOLDS), or  Aczel's and Belo's dependently typed (intuitionistic) first-order logic (DFOL), may be regarded as logic enriched dependent type theories.
Categories with families (cwfs) is an established semantical structure for dependent
type theories, such as Martin-L\"of type theory.  We introduce in this article a notion of {\em hyperdoctrine over a cwf}, and show how FOLDS  and DFOL fit in this semantical framework.  A soundness and completeness theorem is proved for DFOL. The semantics is functorial in the sense of Lawvere, and uses a 
dependent version of the Lindenbaum-Tarski algebra for a DFOL theory. Agreement with standard first-order semantics is established. Applications of DFOL to constructive mathematics and categorical foundations are given. A key feature is a {\em local propositions-as-types principle.}

{\em AMS MSC 2010 classification:} 
                                                 03B15,
                                                 03F50, 
                                                 03G30, 
                                                 18C10, 
                                                68Q55 
                                                
{\em Keywords:} Intuitionistic first-order logic, dependent types, categorical logic, models of type theory
\end{abstract}

\tableofcontents

\section{Introduction}

Dependent families of sets often occur in mathematical notation, for instance as the fibers of a map $f:Y \rightarrow X$,
\begin{equation} \label{mapfiber}
f^{-1}(x) \quad (x\in X)
\end{equation}
$${\rm inj}_x(u) \in Y \quad (x\in X, u \in f^{-1}(x))$$
or as  the hom-sets of a category
$${\rm Hom}(A,B) \quad (A,B\in{\rm Ob})$$
$$1_A \in {\rm Hom}(A,A)  \quad (A \in {\rm Ob}).$$
Such expressions are however not directly supported in first-order logic. One may formulate them in a roundabout way using set-valued functions as in set theory, or directly using dependent types (or sorts) as in Martin-L\"of type theory. Extensions of first-order logic with  
dependent types have been considered by Makkai (1995) and Aczel (2004).
The dependently sorted system of the first-order logic FOLDS was introduced and studied 
in a series of papers (Makkai 1995, 1998, 2013). Its purpose is to provide a  foundation for (higher) category theory that can also serve as a theory of abstract sets. Unlike set theory there is no global notion of equality in the system. A stronger system DFOL was introduced by Belo (2008) following the idea of  Aczel (2004).  It is based on the generalized algebraic theories of Cartmell 1986. In contrast to (Makkai 1995) it allows function symbols.  The system is an extension of many-sorted intuitionistic logic. Moreover an important property is that the system DFOL is straightforwardly interpretable in Martin-L\"of type theory, and thus possible to use in connection with proof assistants such as Coq or Agda.

A {\em category with families} (cwf) is a semantical structure for interpreting dependent type theories in categories (Dybjer 1996, Hofmann 1997, Jacobs 1999). This structure goes back to the {\em contextual categories} of Cartmell (1986). Cwfs have a rich structure which can encode also ordinary categories. Clairambault and Dybjer (2014) show that cwfs equipped with the appropriate type constructions can capture cartesian (or finitely complete) categories, and locally cartesian closed categories, up to bi-equivalence between 2-categories.
A purpose of the present paper is to show how DFOL type systems gives rise to cwfs, and to establish basic properties of these systems.  We
also introduce {\em hyperdoctrines over cwfs} as a general semantics
of first-order logic with dependent sorts. A soundness and completeness theorem is proved for  DFOL. Hyperdoctrines over cwfs should also be a suitable semantic structure for
logic enriched type theories, such as that of Aczel--Gambino.
A complicating aspect of first-order dependent signatures is that the type and function symbol declarations in general interact with each other. One cannot separate them as in standard many-sorted logic. However, every dependent first-order type-and-function signature $\Sigma$ generates a free cwf ${\cal F}_\Sigma$. A {\em structure } $({\cal C},M)$ for the $\Sigma$-signature is defined as a cwf morphism 
$M: {\cal F}_\Sigma\, \to\, {\cal C}$.  Every first-order dependent theory $T$ over a signature $\Sigma$ with 
predicates $\Pi$ generates a hyperdoctrine ${\cal H}_{\Sigma,\Pi,T}$  which is essentially a 
dependently typed Lindenbaum-Tarski algebra.
A  dependent first-order model of $T$ consists of a $\Sigma$-structure
$M: {\cal F}_{\Sigma} \to {\cal C}$ and a hyperdoctrine ${\cal D}$ with an $M$-based hyperdoctrine morphism $G: {\cal H}_{\Sigma,\Pi,T} \to {\cal D}$.
Also these morphisms can be constructed incrementally by induction on the signature.
This shows that this functorial notion of model extends the usual non-dependent version
of model of first-order signature. 

\medskip
The paper is organized as follows. Before the formal development of type systems and DFOL, we illustrate in Section \ref{UDFOL} their use with examples.
This includes treatment of partial functions, the local propositions-as-types property, and   a constructive version of Lawvere's Elementary Theory of the Category of Sets (ETCS). After technical preliminaries in Section \ref{prelsec}, the type systems are set up in Section  \ref{ddtsec}.  Section \ref{catcontsec} shows that the natural syntactic category of contexts coming from the type systems forms a cwf. The extension problem for signatures is considered in Section \ref{extprob}.
A cwf morphism ${\cal F}_\Sigma\, \to\, {\cal C}$ is uniquely determined by its values on  the signature  $\Sigma$ (Lemma \ref{cwfuniq}). In Theorems \ref{extbytype} and \ref{extbyfun} we show that such morphisms can be constructed incrementally by induction on the signature.   Hyperdoctrines over cwfs are formulated in Section \ref{semcomplsec}. The corresponding
   logic DFOL is studied in Section \ref{DFOL}, and a completeness theorem (existence of universal model) and soundness theorem is proved.

\subsubsection*{Acknowledgements} This research was in part financially supported by the Swedish Research Council (VR).
 I am indebted to members of the Stockholm Logic Seminar for comments: Peter LeFanu Lumsdaine, Per Martin-L\"of, and especially H{\aa}kon Robbestad Gylterud for
 discussions that lead to a better understanding of FOLDS signatures in terms of 
 type theories (Palmgren 2016). Thanks goes to the anonymous referee for suggestions that helped improve the presentation of the paper.

\section{Using dependently typed first-order logic} \label{UDFOL}

Before entering the formal development we introduce some examples of the use of DFOL, which should hopefully be intuitively understandable. The basic judgement forms of the type systems are,
$$A \; {\rm type} \; (x_1 : A_1, \ldots, x_n : A_n) $$
which means that {\em $A$ is a type in the context $x_1 : A_1, \ldots, x_n : A_n$,}
and 
$$a : A \; \; (x_1 : A_1, \ldots, x_n : A_n) $$
which means that {\em $a$ is an element of $A$ in the context $x_1 : A_1, \ldots, x_n : A_n$.} For two formulas $\varphi$ and $\psi$ in a context $\Gamma = x_1 : A_1, \ldots, x_n : A_n$, we have the judgement that 
{\em $\varphi$ entails $\psi$} written
$$\varphi \sequentunder{\Gamma} \psi$$
or 
$$ \varphi\Longrightarrow \psi \; \; (\Gamma).$$
The logical rules of DFOL are intuitionistic; see Section \ref{DFOL}. For flexibility we use a categorical adjunction based presentation of logical rules as in Johnstone (2002). 

\medskip
\begin{example} \label{SigmaCat}
{\em
The signature $\Sigma_{\rm Cat}$ for a category is given by 
 the sequence of declarations  below
\begin{equation}
\begin{array}{ll}
 {\sf Ob} \; {\rm type}\quad &() \\
 X \rightarrow Y\; {\rm type}\quad &(X:{\sf Ob}, Y: {\sf Ob})\\
 1(X): X \rightarrow X \quad &(X:{\sf Ob})\\
\circ(g,f): X \rightarrow Z \quad &(X:{\sf Ob}, Y: {\sf Ob},Z:{\sf Ob}, 
g: Y \rightarrow Z, f:X \rightarrow Y)
\end{array}
\end{equation}
}
\end{example}

\medskip
As there is no standard equality on the sorts, we have to axiomatize the equality relations explicitly. Following e.g. (Makkai 2013) we impose
equality on the arrows, but not on the objects
$$f= g \; {\rm formula} \; (X,Y:{\sf Ob}, f,g: X \rightarrow Y).$$
Axioms stating that it is an equivalence relation are added
\begin{equation}
\begin{array}{ll}
 \top \sequent f= f \quad &(X,Y:{\sf Ob}, f: X \rightarrow Y)\\
 f= g  \sequent   g=  f \quad &(X,Y:{\sf Ob}, f,g: X \rightarrow Y) \\
  f= g  \land   g= h \sequent   f= h \quad &(X,Y,Z:{\sf Ob}, f,g,h: X \rightarrow Y)\\
\end{array}
\end{equation}
and that composition respects the equality, where we write $\circ(g,f)$ as 
$g \circ f$,
\begin{equation}
f= h \land   g= k \sequent   g \circ f= k \circ h 
\quad  (X,Y,Z:{\sf Ob}, f,h: X \rightarrow Y,  g,k: Y \rightarrow Z) 
\end{equation}
\begin{equation}
\begin{array}{ll}
  \top \sequent 1(Y) \circ  f= f \quad &(X,Y:{\sf Ob}, f: X \rightarrow Y) \\
   \top \sequent  f \circ  1(X)= f \quad &(X,Y:{\sf Ob}, f: X \rightarrow Y) \\
   \top \sequent  h \circ  (g \circ f) = (h \circ  g) \circ f 
\quad &(X,Y,Z,W:{\sf Ob}, f: X \rightarrow Y, g: Y \rightarrow Z, h: Z \rightarrow W)\;  \\
\end{array}
\end{equation}

This is a principal example of a DFOL theory. It is continued in Section \ref{CETCSsec}
showing how to axiomatize certain categorical constructions.

\medskip
Already the type system of DFOL has some power to encode a fragment of logic, using type inhabitation as truth, as in the following example.

\begin{example} \label{semigp}
{\em
The following
is a signature for a binary operation $m$ with a congruence relation $E$.
Inhabitedness of the type $E(x,y)$ is understood as that the relation $E(x,y)$ holds.
$$
\begin{array}{ll}
A \; {\rm type} \;  \; &() \\
m(x,y): A \;  \; &(x,y:A) \\
E(x,y) \; {\rm type}\;  \; &(x,y:A) \\
\rho(x): E(x,x)\;  \; &(x:A)\\
\sigma(p): E(y,x)\;  \; &(x,y:A, p:E(x,y)) \\
\tau(p,q): E(x,z)\;  \; &(x, y, z:A, p:E(x,y), q:E(y,z)) \\
\gamma(p,q): E(m(x,y),m(u,v))\;  \; &(x, y, u, v:A, p:E(x,u), q:E(y,v)) \\
\alpha(x,y,z) : E(m(x,m(y,z)), m(m(x,y),z)) \; \; &(x, y, z:A)
\end{array}
$$
Now we can, for example, add generators $a, b, c: A$ and a relation  $ab = bc$ to specify a particular semigroup. This is done by introducing the axiom
$${\rm ax}_1: E(m(a,b), m(b,c)) \; ()$$
Using this method we may encode the semigroup word problem which entails that type inhabitation is not decidable; see Remark \ref{decision}.
}
\end{example}

\subsection{Local propositions-as-types} \label{usesDFOL}

A feature of the dependently typed logic is that it has a {\em local propositions-as-types} principle. 
This means the possibility to extend the signature with a type for the proof objects of a particular formula. This is useful in the definition of functions whose domains
are restricted by a formula.

Consider a fixed signature. Suppose $\Gamma=x_1:A_1,\ldots,x_n:A_n$ is a context and that $\phi$ is a formula in that context.
Add a new dependent type to the signature
$$ F(x_1,\ldots,x_n) \; {\rm type} \quad (\Gamma)$$
Then add two axioms over the extended signature
$$ \top \sequent \phi \quad (\Gamma, p:F(x_1,\ldots,x_n))$$
$$ \phi \sequentunder{\Gamma} (\exists p:F(x_1,\ldots,x_n))\top.$$
The truth of $\phi(x_1,\ldots,x_n)$ has thus been encoded as inhabitedness of $F(x_1,\ldots,x_n)$.
We call $F$ a {\em proposition-type} for $\phi$.

\subsection{Treatment of partial functions}

In classical (non-dependent, many-sorted) logic a partial functional relation can be extended to a function.
Suppose
\begin{equation}\label{phiext}
x=_A u \land y=_B z \land \phi(x,y) \sequent \phi(u,z)
\end{equation}
and
\begin{equation}\label{phifun}
\phi(x,y) \land \phi(x,z) \sequent y=_B z
\end{equation}
Assume that $b$ is some constant in $B$. We can introduce a total function symbol $f:A \rightarrow B$ with
defining axiom 
$$f(x)=_B y \leftrightarrow \phi(x,y) \lor y=_B b \land \lnot (\exists y:B)\, \phi(x,y).$$
In intuitionistic logic this is not possible, unless the domain of definition of the relation $\phi$ is decidable.
However in DFOL we can do this with the help of dependent types.
Suppose again (\ref{phiext}) and (\ref{phifun}) holds.
Introduce a proposition-type $D$ for the domain of definition of $\phi$
$$D(x) \; {\rm type} \quad (x:A)$$
and an axiom
$$\phi(x,y) \sequent (\exists p: D(x))\top\quad (x:A, y:B)$$
and a dependent function symbol $f$
$$f(x,p):B \quad (x:A, p:D(x))$$
with the axiom
$$\top \sequent \phi(x,f(x,p)) \quad (x:A, p:D(x))$$

This works constructively as desired, as we can interpret $D(x)$ in e.g.\ Martin-L\"of type theory as $(\Sigma y:B)\phi(x,y)$ and $f(x,p)$ as $\pi_1(p)$, the first projection of $p$.

\subsection{Equality on dependent types} \label{depeqsec}

If we wish to introduce equality on dependent types there is the established method
of  using setoids (type with equivalence relations) and proof-irrelevant families of setoids; see for instance (Palmgren 2012b).

Suppose that $A$ is a type with an equivalence relation $=_A$. For a dependent type $B$ over $A$ we introduce an equivalence relation $=_{B,x}$, as follows
\begin{equation} \label{eq20}
\begin{array}{ll}
u =_{B, x} v\; {\rm formula} \quad &(x:A, u,v: B(x))\\
 \top \sequent  u =_{B, x} u  \quad &(x:A, u: B(x)) \\
  u =_{B, x} v \sequent v =_{B, x} u \quad &(x:A, u,v: B(x)) \\
u =_{B, x} v \land  v =_{B, x} w \sequent u =_{B, x} w \quad &(x:A, u,v,w: B(x))
\end{array}
\end{equation}
How do elements of $B(x)$ and $B(y)$ compare if $x=_A y$ is true? We assume  that $=_A$ is associated with a proposition-type $ E_A$ as follows
\begin{equation} \label{eq21}
\begin{array}{ll}
E_A(x,y) \; {\rm type} \quad &(x,y:A) \\
 \top \sequent  x=_Ay  \quad &(x,y:A, p:E_A(x,y)) \\
x=_Ay  \sequent (\exists p:E_A(x,y)) \top \quad &(x, y:A) \\
\end{array}
\end{equation}
We now introduce transport functions ${\rm tr}_{p}$ to be able to relate elements of $B(x)$ and $B(y)$ when $x=_Ay$ 
\begin{equation} \label{eq22}
\begin{array}{ll}
{\rm tr}_{p}(u): B(y) \quad &(x,y:A,p: E_A(x,y), u:B(x) )\\
 u=_{B, x} v \sequent  {\rm tr}_{p}(u) =_{B, y} {\rm tr}_{p}(v) \quad &(x,y:A,p: E_A(x,y), u,v:B(x))\\
  \top \sequent {\rm tr}_{p}(u) =_{B, y} {\rm tr}_{q}(u) \quad &(x,y:A,p,q: E_A(x,y), u :B(x)) \\
  \top \sequent {\rm tr}_{p}(u) =_{B, x} u \quad &(x:A,p: E_A(x,x), u:B(x))\\
  \top \sequent {\rm tr}_{q}({\rm tr}_{p}(u)) =_{B, z} {\rm tr}_{r}(u) \quad &(x,y,z:A,q: E_A(y,z), p: E_A(x,y), \\
 & \qquad r: E_A(x,z),u:B(x))  \\
\end{array}
\end{equation}
Note that second sequent says that the transport function is independent of the particular value of the proof object in $E_A(x,y)$. This makes the family
$B$ over $A$ {\em proof-irrelevant.}

These definitions can be extended to arbitrary contexts of setoids, see (Maietti 2009).

We may use the above to axiomatize the fibers of a function $f:C \to A$ as in (\ref{mapfiber}). Let $=_C$ be denote the
equality relation on $C$. We think now of $B(x)$ as $f^{-1}(x)$. Introduce  injective functions
$$i_x(u):C \quad (x:A, u:B(x))$$ 
satisfying
\begin{equation}
\begin{array}{ll}
 \top \sequent f(i_x(u))=_A x \quad & (x:A, u:B(x)) \\
 \top \sequent  i_y({\rm tr}_{p}(u)) =_C i_x(u) \quad &(x,y:A, p: E_A(x,y), u:B(x)) \\
u=_{B, x} v \sequent i_x(u) =_C i_x(v) \quad &(x:A, u,v:B(x)) \\
i_x(u) =_C i_x(v) \sequent u=_{B, x} v \quad &(x:A, u,v:B(x)) \\
f(w) =_A x \sequent (\exists u:B(x))\; i_x(u) =_C w \quad &(w:C, x:A) \\
\end{array}
\end{equation}

\subsection{Constructive elementary theory of the category of sets} \label{CETCSsec}

Lawvere's Elementary Theory of the Category of Sets  ETCS (Lawvere 1964) is a first-order theory that axiomatizes a sufficiently rich category of sets for mathematics. In (Palmgren 2012a) a constructive version of this theory was proposed. The theory, CETCS, used equality of objects, which is not entirely natural for straightforward interpretation as setoids in Martin-L\"of type theory. Here  we provide a version  which avoids equality of objects.


The signature for a CETCS category is a follows. In addition to $\Sigma_{\rm Cat}$ we add declarations for terminal object and pullbacks 
\begin{enumerate}
\item $1:{\sf Ob} \quad ()$
\item ${\sf P}(f,g): {\sf Ob} \quad (X,Y,Z:{\sf Ob}, f:X \rightarrow Z, g: Y \rightarrow Z)$
\item $    {\sf p}_1(f,g):
{\sf P}(f,g) \rightarrow X  \quad (- " -)$
\item $ {\sf p}_2(f,g):
{\sf P}(f,g) \rightarrow Y  \quad (- " -)$
\end{enumerate}
and similarly for initial object and pushouts. Further we add the $\Pi$-construction
\begin{itemize}
\item[5.] $ \Pi(f,g):{\sf Ob} \quad (X,Y,Z:{\sf Ob}, f:X \rightarrow Y, g: Z \rightarrow X) $
\item[6.] $ \pi(f,g): \Pi(f,g) \rightarrow Y \quad ( - " -)$ 
\item[7.] $  {\rm ev}(f,g):{\sf P}(f,\pi(f,g)) \rightarrow Z \quad (- " -)$
\end{itemize}
Moreover constants for the presentation axiom
\begin{itemize}
\item[8.] ${\rm Pre}(X): {\sf Ob} \quad (X: {\sf Ob})$
\item[9.] ${\sf pre}(X): {\rm Pre}(X) \rightarrow X \quad (X: {\sf Ob})$
\end{itemize}

The theory CETCS of (Palmgren 2012a) was presented in ordinary many-sorted intuitionistic logic with equality on all sorts. We here
explain how to reformulate the theory using dependent sorts, and thereby avoiding equality on objects. Below only a few crucial and
illustrative axioms are considered. The remaining axioms follow the same pattern.

\medskip
That the object $1$ is terminal is formalized as
\begin{equation}
\begin{array}{l}
 \top \sequent f = g \quad (X:{\sf Ob},f,g: X \rightarrow 1)\\
\\
\top \sequent (\exists f:X \rightarrow 1) \top  \quad(X:{\sf Ob}) \\
\end{array}
\end{equation}
(We could as well expanded the signature with an explicit function for selecting a map to the terminal $!_x:X \rightarrow 1$.)
The existence of pullbacks is given by these  axioms for ${\sf P}$, ${\sf p}_1$ and ${\sf p}_2$:
\begin{equation}
\begin{array}{l}
 \top \sequent  f \circ {\sf p}_1(f,g)  
=  g \circ {\sf p}_2(f,g) 
\qquad (X,Y,Z:{\sf Ob}, f:X \rightarrow Z, g: Y \rightarrow Z) \\ \\
 f \circ h = g \circ k \sequent  (\exists t:W \rightarrow {\sf P}(f,g))\, {\sf p}_1(f,g) \circ t = h \land 
{\sf p}_2(f,g) \circ t = k \\
\qquad \qquad (X,Y,Z,W:{\sf Ob}, f:X \rightarrow Z, g: Y \rightarrow Z, h:W \rightarrow X, k: W \rightarrow Y)  \\ \\
 {\sf p}_1(f,g) \circ t = {\sf p}_1(f,g) \circ s \land {\sf p}_2(f,g) \circ t = {\sf p}_2(f,g) \circ s
\sequent  t= s\\
\qquad \qquad  (X,Y,Z,W:{\sf Ob}, f:X \rightarrow Z, g: Y \rightarrow Z,  t,s: W \rightarrow{\sf P}(f,g)) \\
\end{array}
\end{equation}
We have similar axioms for other limits and colimits.

Recall that a morphism $x:1 \to X$ is called an {\em element} of $X$. We write this as $x \in X$. 
A morphism $f: X \to Y$ is 
called {\em surjective} if for every $y\in Y$, there is $x\in X$ such that $f \circ x = y$. The following axiom states that $1$ is a {\em strong generator}
of the category:
\begin{itemize}
\item[(G)] Every surjective mono $X \to Y$ is an isomorphism.
\end{itemize}
This statement can straightforwardly be written as a first-order formula.

We consider as in (Palmgren 2012a) a sequence of maps $\alpha_1:P \to X_1,\ldots,\alpha_1:P \to X_n$ that are jointly monic as a relation between the objects $X_1,\ldots,X_n$
and for elements $x_1\in X_1,\ldots,x_n\in X_n$ we write $$(x_1,\ldots,x_n) \memof (\alpha_1,\ldots,\alpha_n)$$ if there is $p \in P$ with
$\alpha_1 \circ p = x_1,\ldots,\alpha_n \circ p = x_n$. A relation $\alpha_1:P \to X_1, \alpha_2:P \to X_2$ is called a {\em partial function
 from $X_1$ to $X_2$} if $\alpha_1$ is mono.

\medskip
Existence of dependent products is formulated as follows (Palmgren 2012a): For any 
mappings  $Y \to^g X \to^f I $  we have a commutative diagram
\begin{equation}  \label{defpid}
\bfig                                                                          
\qtriangle/<-`->`/<700,700>[Y``;\piev(f,g)`g`]         
\square(760,0)<1000,700>[\qquad\qquad{\sf P}(f,\pi(f,g))`\Pi(f,g)`X`I;{\sf p}_1`{\sf p}_2`\pi(f,g)`f]
\efig
\end{equation}
where 
for any element $i \in I$  and any partial function
$\psi=_{\rm def}(\xi: R \to X,\upsilon:R \to Y)$
such that
\begin{itemize}
\item[(a)] for all $x\in X$, $y \in Y$,  
$(x,y) \memof \psi$ implies $g \circ y = x$ and $f \circ x = i,$
\item[(b)] for all $x\in X$, with $f \circ x = i$, there is 
$y \in Y$ with $(x,y) \memof \psi$,
\end{itemize}
then there is a unique $s \in \Pi(f,g)$ so that $\pi(f,g) \circ s = i$ and 
for all  $x\in X$, $y \in Y$,  
\begin{equation} \label{equival1}
(s,x,y) \memof \alpha   \text{ iff }
(x,y) \memof \psi.
\end{equation}
Here $\alpha=_{\rm def}({\sf p}_1:{\sf P} \to \Pi(f,g),{\sf p}_2:{\sf P} \to X ,\piev(f,g):{\sf P} \to Y)$, where ${\sf P}= {\sf P}(f,\pi(f,g))$.

\section{Preliminaries}  \label{prelsec}

The mathematical background theory used in development of the results of this article is constructive, in the sense that it uses only intuitionistic logic, and set constructions which are generalized predicative, and thus, in particular, no power sets are used. We leave it unspecified whether it is constructive set theory CZF (Aczel and Rathjen 2010) or Martin-L\"of type theory. We employ terminology from constructive mathematics of the Bishop school (Bridges and Richman 1987). So for instance a set is {\em discrete} if for each pair of its elements $x$ and $y$, $x=y$ or $\lnot x = y$.

To treat the issue of fresh variables relative to a context of variable declarations
we introduce the following notions.
 
\begin{definition} {\em
Let $V$ be an infinite discrete set. A  {\em fresh variable provider (fvp) for $V$} is a pair of functions $({\rm \varphi}, {\sf fr})$ which to each finite subset $X\subseteq V$, assigns an inhabited subset ${\rm \varphi}(X) \subseteq V \setminus X$, and an element  ${\sf fr}(X) \in {\rm \varphi}(X)$. }
\end{definition}
The first function
provides the variables available after a context has used the variables $X$, the second function
selects one of these variables. We call the tuple ${\cal V}=(V, ({\rm \varphi}, {\sf fr}))$ a {\em variable system.}  A prime example of a fvp is defined by 
\begin{equation} \label{unrestr}
{\rm \varphi}_{\infty}(\{x_1,\ldots,x_n\}) = 
\{ y \in V :  \mbox{$\lnot y = x_i$ for all $i=1,\ldots,n$} \},
\end{equation}
and ${\sf fr}_\infty(\{x_1,\ldots,x_n\})$ selects some element in the set.
Another example is when $V= \nat_+$,  
${\sf fr}_1(\{x_1,\ldots,x_n\}) = {\rm max} \{1, x_1+1,\ldots,x_n+1\}$,
$${\rm \varphi}_1(\{x_1,\ldots,x_n\}) = \{ {\sf fr}_1(\{x_1,\ldots,x_n\})\},$$
Here ${\rm \varphi}_1(\emptyset) = \{1\}$. A variable system $ {\cal V} =(V, ({\rm \varphi}, {\sf fr}))$ is said to have  the {\em de Bruijn property} if ${\rm \varphi}(\{x_1,\ldots,x_n\}) = \{ {\sf fr}(\{x_1,\ldots,x_n\})\}$. This will have the effect that the $n$th declared variable of any two contexts will be the same.  The system is called {\em unrestricted} if $\varphi$ satisfies (\ref{unrestr}). If ${\cal V} =(V, ({\rm \varphi}, {\sf fr}))$ is an arbitrary variable system, we write
${\cal V}^{(1)}$ for $(V, ({\rm \varphi}_{\sf fr}, {\sf fr}))$ where ${\rm \varphi}_{\sf fr}(\{x_1,\ldots,x_n\}) = \{ {\sf fr}(\{x_1,\ldots,x_n\})\}$. This system has the de Bruijn property.
We write ${\cal V}^{(\infty)}$ for $(V, ({\rm \varphi}_{\infty}, {\sf fr}))$, which is an unrestricted variable system.
 
We define  untyped terms or preterms inductively. Let ${\rm Ter}(M,V)$ denote the set of {\em preterms} which is the smallest set  $T$ such that $V \subseteq T$, and for $t_1,\ldots,t_n \in T$, $n \ge 0$, and any $\sigma \in M$, $\sigma(t_1,\ldots, t_n) \in T$.  We write just $\sigma$ for $\sigma()$. Syntactic equality of terms is written $s \equiv t$. We assume that $M$ is a discrete set and that $V$ is a
 discrete set equipped with a fvp. We use standard notation for simultaneous substitution $s[t_1,\ldots,t_n/x_1,\ldots,x_n]$ and the set  ${\rm V}(E)$ of variables in an expression $E$. We will use the following property as well. 
 
 \begin{prop} For every $s \in  {\rm Ter}(M,V)$, with ${\rm V}(s)= \{x_1,\ldots,x_n\}$,
 if $t_1,\ldots, t_n,$ $t'_1,\ldots, t'_n \in {\rm Ter}(M,V)$,
 $$s[t_1,\ldots,t_n/x_1,\ldots,x_n] \equiv s[t'_1,\ldots,t'_n/x_1,\ldots,x_n]
 \Longrightarrow t_1 \equiv t'_1, \ldots, t_n \equiv t'_n.
 $$
 \end{prop}
 {\flushleft \bf Proof.} Induction on $s$. $\qed$

\section{Dependently typed terms} \label{ddtsec}

The following class of type systems is essentially those of Belo (2008), and is close to those of Cartmell (1986) except that they not use equality of types. A difference to (Belo 2008) is that we allow hidden arguments to function and type symbols, similarly to certain proof assistants (Coq and Agda), a device which can shorten expressions considerably. Our basic dependent
type systems can code a limited form of Horn logic and equational logic, 
using the propositions-as-types principle.  Later we introduce full first-order logic on top of the type systems.

\subsection{Presyntax}
Let $V$ be an infinite discrete  set with a fvp $({\rm \varphi},{\sf fr})$, so that ${\cal V} = (V,({\rm \varphi},{\sf fr}))$ forms a variable system.
Let $F$ and $T$ be two discrete disjoint sets and let $M= F \cup T$. We call 
${\cal S}=({\cal V},F,T)$ a {\em symbol system}, where $F$ is intended to be the set of function symbols and $T$ to be the set of type symbols. For any syntactic expression $E$ we denote by ${\rm V}(E)$ the set of variables in $E$.

A term $t \in {\rm Ter}(M,V)$ is called a {\em preelement} over ${\cal S}$ if it does not contain symbols from $T$; it is a {\em pretype} over ${\cal S}$
if it is of the form $S(t_1,\ldots,t_n)$ where all $t_1,\ldots,t_n$ are preelements and $S \in T$.

\begin{definition} \label{precon}
{\em 
A {\em precontext} over ${\cal S}$ is a sequence $\Gamma=\langle x_1:A_1,\ldots,x_n:A_n\rangle$, $n\ge 0$, where $A_1,\ldots,A_n$ are pretypes over ${\cal S}$,  $x_1,\ldots,x_n$ are all in $V$ with $x_k \in {\rm \varphi}(\{x_1,\ldots,x_{k-1}\})$ for all $k=1,\ldots,n$, 
and  ${\rm V}(A_k) \subseteq \{x_1,\ldots,x_{k-1}\}$ for all $k=1,\ldots,n$.  Each $x_i:A_i$ is called a {\em variable declaration} and declares the variable $x_i$ to be of type $A_i$. We define the set of variables of a precontext as ${\rm V}(\Gamma) = \{x_1,\ldots,x_n\}$. The ordered sequence of variables of the precontext  $\Gamma$ is denoted 
$${\rm OV}(\Gamma) = x_1,\ldots,x_n \text{ or  as } \check{\Gamma} = x_1,\ldots,x_n.$$ 
The {\em length $|\Gamma|$} of the context $\Gamma$ is $n$.
}
\end{definition}
Below we will often use the following notation for substitution  
$$E[a_1,\ldots,a_n/\Gamma] =_{\rm def}E[a_1,\ldots,a_n/{\rm OV}(\Gamma)],$$
as well as $E[\bar{a}/\Gamma]$ where $\bar{a} = a_1,\ldots,a_n$.
Define the following operations on contexts using the fvp. For a precontext
$\Gamma$ let
$${\rm Fresh}(\Gamma) = {\varphi}({\rm V}(\Gamma)) \qquad
{\rm fresh}(\Gamma) = {\rm fr}({\rm V}(\Gamma)).$$
The former is the set of variables that are fresh for $\Gamma$ and the latter
is a choice of one variable that is fresh for $\Gamma$.
The set of {\em top-most variables} ${\rm TV}(\Gamma)$ of a precontext $\Gamma$ is defined by induction on the length of $\Gamma$:
\begin{eqnarray*}
{\rm TV}(\langle \rangle) &= &  \emptyset\\
{\rm TV}(\langle \Gamma,x:A\rangle) &= &  ({\rm TV}(\Gamma) \setminus {\rm V}(A)) \cup \{x\}\\
\end{eqnarray*}
Thus for example
$${\rm TV}(\langle x:S,y:T(x), z:R(x,y),u:U(x)\rangle) = \{z,u\}.$$
The values of all other variables will be shown to be derivable from the top-most variables in a precise sense (Theorem \ref{UniqType}).
Precontexts (and contexts) $\langle x_1:A_1,\ldots,x_n:A_n\rangle$ will usually be written 
$$x_1:A_1,\ldots,x_n:A_n$$
without brackets. We allow abbreviations when the same type occurs several times in juxtaposition, so that for instance
$$x:A,y,z:B,u:C,v,w:D$$
abbreviates
$$x:A,y:B,z:B,u:C,v:D,w:D.$$

\begin{definition} \label{typedecl}
{\em 
A {\em type predeclaration}  over ${\cal S}$  is a pair $(\Gamma, S,\bar{i})$ where $\Gamma$ is a precontext over ${\cal S}$ and $S \in T$, and moreover $\bar{i}=i_1,\dots,i_k$ is a strictly increasing sequence of indexes in $\{1,\ldots,n\}$ such that 
\begin{equation} \label{hiding}
{\rm TV}(\Gamma) \subseteq \{x_{i_1},\ldots,x_{i_k}\} 
\end{equation}
where ${\rm OV}(\Gamma)= x_1,\ldots,x_n$.
Such $\bar{i}$  is called a {\em determining sequence}, and  we let ${\rm DS}(\Gamma)$ denote the set of these. 
 If $\Gamma = \langle x_1:A_1,\ldots,x_n:A_n\rangle$, we write $(\Gamma, S, i_1,\ldots i_k)$   as 
 \begin{equation} \label{Sdecl}
S(x_{i_1},\ldots,x_{i_k}) \; {\rm type} \qquad (\Gamma).
\end{equation}
The constant {\em declared} is $S$. The relaxed condition (\ref{hiding}) allows for hiding variables whose values are derivable. In case $\bar{i}=1,2,\ldots,n$, we say that the declaration is  {\em regular}. This is to say, all variables are used and
come in the same order as in the context. 
}
\end{definition}
 
 \begin{definition} \label{typedecl}
{\em 
 A {\em function predeclaration} over ${\cal S}$ is a triple $(\Gamma, f, \bar{i}, U)$ where $\Gamma$ is a precontext over ${\cal S}$ and $f \in F$  where $U$ is a pretype over ${\cal S}$  with 
${\rm V}(U) \subseteq {\rm V}(\Gamma)$, and  $\bar{i} \in {\rm DS}(\Gamma)$. If $\Gamma = \langle x_1:A_1,\ldots,x_n:A_n\rangle$, we write $(\Gamma, f,  i_1,\ldots,i_k, U)$   as 
 \begin{equation} \label{fdecl}
f( x_{i_1},\ldots,x_{i_k}):U \quad (\Gamma).
\end{equation}
The constant {\em declared} is $f$. Again,  in case $\bar{i}=1,2,\ldots,n$, we say that the declaration is
 {\em regular}. 
}
\end{definition}

For regular declarations we often omit
the determining sequence $\bar{i}$. If $D$ is any predeclaration, then we denote the 
constant declared by ${\rm decl}(D)$.

\medskip
For instance, the predeclarations tuples for $A$, $E$ and $\sigma$ of 
Example \ref{semigp} would read 
$(\langle\rangle,\varepsilon, A)$, $(\langle x:A, y:A \rangle,(1,2), E)$ and
$(\langle x:A, y:A, z:E(x,y) \rangle, (3), \sigma)$, respectively.
\medskip

\begin{definition} \label{presig}
{\em
A (unordered)  {\em presignature} over ${\cal S}$ is a set $\Sigma$ of type predeclarations and function predeclarations satisfying the {\em consistency condition}: if $D_1, D_2 \in \Sigma$, ${\rm decl}(D_1) = {\rm decl}(D_2)$, then $D_1=D_2$.  A  presignature
 $\Sigma$ is  {\em regular} if each of its predeclarations is regular.}
\end{definition}

Below we shall use the following notation:  for $\bar{a}= a_1,\ldots,a_n$ and 
$\bar{i} = i_1,\ldots,i_k$, with $1 \le i_1,\ldots,i_k \le n$, write
$$\bar{a}_{\,\bar{i}}= a_{i_1},\ldots,a_{i_k}.$$

\begin{definition} \label{judgexpr}
{\em
A {\em judgement expression} is any of the following three expressions
\begin{itemize} 
\item[] $\Gamma \; {\rm context}$ 
\item[] $A \; {\rm type}\;  (\Gamma)$
\item[] $a : A\; (\Gamma)$
\end{itemize}
where $\Gamma$ is a precontext, $A$ is a pretype and, $a$ is a preelement.
}
\end{definition}

\subsection{Basic dependent type systems}

Let a symbol system ${\cal S}=({\cal V}, F ,T)$ be fixed. 
For any presignature $\Sigma$ over this symbol system, 
let ${\cal J}(\Sigma)$ be the smallest set of judgement 
expressions closed under following derivation rules (R1) -- (R5) below.

$$\infer[({\rm R1})]{\langle \rangle \; {\rm context}}{ }$$

$$\infer[({\rm R2}) \quad  x \in {\rm Fresh}(\Gamma)]{\Gamma,x:A\; {\rm context}}{ \Gamma\; {\rm context} &
    A\; {\rm type}\; (\Gamma) }$$

$$\infer[({\rm R3})]{x_i:A_i  \;(x_1:A_1,\ldots, x_n:A_n)}{  x_1:A_1,\ldots, x_n:A_n \; {\rm context} }$$

To formulate the next two rules we introduce the notion of context map.
Let $\Delta$ and $\Gamma= x_1:A_1,\ldots,x_n:A_n$ be two precontexts. 
 We write 
 \begin{equation} \label{cmap}
 (a_1,\ldots,a_n): \Delta \to \Gamma 
 \end{equation}
 for the conjunction of these $n+2$ judgement expressions:
\begin{equation} \label{cmpjudgement}
\begin{array}{lll} 
 &\Delta  \; {\rm context}  \\
 &\Gamma  \; {\rm context}  \\
 &  a_1: A_1 \; (\Delta) \\ 
&  a_2: A_2[a_1/x_1] \;(\Delta) \\
&\quad \vdots \\
& a_n: A_n[a_1,\ldots,a_{n-1}/x_1,\ldots,x_{n-1}] \;(\Delta)
 \end{array}
\end{equation}
Note that (\ref{cmap}) is to be thought of saying that $\bar{a}=(a_1,\ldots,a_n)$ is a context map  or substitution from $\Delta$ into $\Gamma$.  For $n=0$, the conjunction is
equivalent to the judgement $\Delta  \; {\rm context}$. The sequence $\bar{a}$
 can be inserted into correct function and type declarations according to the next rules.

$$\infer[    ({\rm R4}) \quad \begin{array}{l} (\Gamma,S, \bar{i}) \mbox{ in  } \Sigma\\
  \end{array}   ]{S(\bar{a}_{\,\bar{i}}) \; {\rm type} \;(\Delta) }{\bar{a}: \Delta \to \Gamma 
}$$

$$\infer[({\rm R5}) \quad \begin{array}{l} (\Gamma,f,\bar{i},U) \mbox{ in  } \Sigma \\
 \end{array}]{ f(\bar{a}_{\,\bar{i}}) : U[\bar{a}/\Gamma] \;(\Delta) }{ 
 \bar{a}: \Delta \to\Gamma  & U[\bar{a}/\Gamma] \; {\rm type} \;(\Delta)   \\
}$$

Note that rules R4 and R5  have $n+2$ and $n+3$ premisses, respectively.

\medskip
We introduce the following notation for later convenience: for  $\Gamma = x_1:A_1,\ldots,x_n:A_n$ and $k \le n$ write
$$\Gamma^k = x_1:A_1,\ldots,x_k:A_k$$
and 
$$\bar{a}^k = (a_1,\ldots,a_k).$$
Thus, for instance, the general term in (\ref{cmpjudgement}) above is
$$a_k: A_k[\bar{a}^{k-1}/\Gamma^{k-1}] \; (\Delta).$$

\begin{remark} {\em Belo (2007) uses the notation $(\Gamma) \; A$ for our $A\; {\rm type} \;(\Gamma)$, and  $(\Gamma) \; a: A$ for our $a: A\; (\Gamma)$.}
\end{remark}

\medskip
We define following (Belo 2007):
\begin{definition} \label{signatur}
{\em  A pre-signature $\Sigma$ to be a {\em signature}
if the following {\em correctness conditions} hold:
\begin{itemize}
\item[(i)] If $(\Gamma, S, \bar{i}) \in \Sigma$,  
then $(\Gamma  \; {\rm context})\in {\cal J}(\Sigma)$,
\item[(ii)]  If $(\Gamma, f, \bar{i}, U) \in \Sigma$, then  
$(U \; {\rm type}\;  (\Gamma)) \in {\cal J}(\Sigma)$.
\end{itemize}
}
\end{definition}

\medskip
Returning briefly to Example \ref{semigp} for illustration: to prove that the list of predeclarations forms a correct signature we build it up inductively. We note that 
$\Sigma_1=\{(\langle\rangle, A,\varepsilon)\}$ is
a signature since $(\langle\rangle \; {\rm context}) \in {\cal J}(\emptyset)$. Next to
prove $\Sigma_2= \Sigma_1 \cup \{(\langle x:A, y:A \rangle, E, (1,2))\}$ is a signature, 
we need to show that $(\langle x:A, y:A \rangle\; {\rm context}) \in {\cal J}(\Sigma_1)$.
This is done by first deriving $(\langle x:A\rangle\; {\rm context}) \in {\cal J}(\Sigma_1)$ with
rule (R2). Then $((): \langle x:A\rangle \to \langle\rangle) \in {\cal J}(\Sigma_1)$ is a context map, so (R4) may be applied to obtain $(A \; {\rm type}\; (x:A)) \in {\cal J}(\Sigma_1)$.
Thus by (R2), we get the desired $(\langle x:A, y:A \rangle\; {\rm context}) \in {\cal J}(\Sigma_1)$. Derivations like these are much simplified by the weaking, strengthening and substitution theorems established in Section \ref{strthms} below.

\begin{remark} \label{ctxstandard}
{\em Note that ${\cal J}(\Sigma)$ depends only on $V$ and the function $\varphi$ of the variable system (not ${\sf fr}$).  }
\end{remark}

\medskip
The {\em height of a derivation} is the height of a derivation tree, take the height of the derivation (R1) to be 0.   If ${\cal U}$ is a judgement we write
$\vdash_n {\cal U}$ for ${\cal U}$ is derivable with a derivation tree of height at most $n$.  
We will also be interested in what the minimum number of applications of particular rule $R$.
Write $\vdash^R_n {\cal U}$ for ${\cal U}$ is derivable with a derivation tree where 
at most $n$ applications of rule $R$ occur along any branch.
For a signature $\Sigma$ we denote by 
\begin{equation} \label{gradedJ}
{\cal J}_n(\Sigma) = \{ {\cal U} \in {\cal J}(\Sigma) : \quad \vdash_n {\cal U} \}.
\end{equation}
For the compound context map judgement,
we write $$\vdash_{m_1,m_2,m_3,\ldots,m_{n+2}} (a_1,\ldots,a_n): \Delta \to \Gamma$$
if the height of the $n+2$ derivation trees  of the judgements in  (\ref{cmpjudgement}) are 
at most $m_1,m_2,m_3,\ldots,m_{n+2}$, respectively.

\medskip
The following are evident.

\begin{lemma} \label{Jmono}
For pre-signatures $\Sigma \subseteq \Sigma'$, 
${\cal J}(\Sigma) \subseteq {\cal J}(\Sigma')$. $\qed$
\end{lemma}

\begin{lemma} 
\begin{itemize}
\item[(i)] The empty set is a signature.
\item[(ii)] If $\Sigma_i$ $(i \in I)$ is a family of signatures, then  $\cup_{i \in I} \Sigma_i$ is a signature whenever it satisfies the consistency condition. 
\end{itemize}
$\qed$
\end{lemma}

 We have moreover two
possibilities to extend a signature.
\begin{lemma}  \label{onext}
Let $\Sigma$ be a signature. Suppose that
$(\Gamma\; {\rm context}) \in {\cal J}(\Sigma)$.
\begin{itemize}
\item[(i)] If $S \in T$ not declared in $\Sigma$  and $\bar{i} \in  {\rm DS}(\Gamma)$, 
then $\Sigma \cup \{ (\Gamma,S,\bar{i})\}$ is a signature.
\item[(ii)] Suppose that $(U \; {\rm type}\; (\Gamma)) \in {\cal J}(\Sigma)$.
If $f \in F$ not declared in $\Sigma$  and $\bar{i} \in  {\rm DS}(\Gamma)$, 
then $\Sigma \cup \{ (\Gamma,f, \bar{i},U)\}$ is a signature.

\end{itemize}
$\qed$
\end{lemma}

A signature is called {\em inductive} if it is built from the empty signature using Lemma
\ref{onext}.

\begin{remark}{\em Ordinary many-sorted signatures are those where the
contexts of type declarations are always empty, i.e. the types do not depend
on other types.
}
\end{remark}

\begin{remark} {\em If we restrict signatures to  consist of only type declarations, we get a type system akin to  that of 
FOLDS  (Makkai 1995); see Section 4 of  (Palmgren 2016) for a proof of equivalence with Makkai's category-theoretic formulation.}
\end{remark}

\begin{remark} {\em The dependency between type declaration and function declarations may in general alternate in a signature. A rich class of declarations
is obtained by declaring a single universal dependent type:
\begin{itemize}
\item[] $U\; {\rm type}\; ()$,
\item[] $T(x)\; {\rm type}\; (x:U)$.
\end{itemize}
Further dependent types can then be simulated by declaring new terms in the "universe" $U$, for instance:
\begin{itemize}
\item[] $a: U\; ()$,
\item[] $b(y): U \; (y:T(a))$.
\end{itemize}
Now we can derive
$$T(b(y)) \; {\rm type} \; (y:T(a)),$$
which in effect works as a new dependent type: $B(y) \; {\rm type} \; (y:A)$.
}
\end{remark}

\begin{lemma}  Let $\Sigma$ be a signature. Suppose that $D$ is one of the declarations
obtained in Lemma \ref{onext}, and that it is regular. Then if ${\cal U} \in {\cal J}(\Sigma \cup \{D\})$, and ${\cal U}$
does not contain the symbol ${\rm decl}(D)$, then already 
${\cal U} \in {\cal J}(\Sigma)$. 
\end{lemma}
{\bf \flushleft Proof.} Induction on derivations. The only rule applications 
in a derivation where the symbol ${\rm decl}(D)$ can disappear from a 
judgement expression is R3 and R4. Now since the declarations are 
regular, nothing is lost going from $\bar{a}$ to $\bar{a}_{\bar{i}}$.
Also by assumption $(\Gamma\; {\rm context}) \in {\cal J}(\Sigma)$. Thus
$\bar{a}: \Delta \to \Gamma$ does not contain ${\rm decl}(D)$, and is therefore
in ${\cal J}(\Sigma)$ by inductive hypothesis. $\qed$

\dontshow{
\begin{lemma}  Let $\Sigma$ be a signature. Suppose that
$(\Gamma\; {\rm context}) \in {\cal J}(\Sigma)$, 
$S \in T$ not declared in $\Sigma$  and $\bar{i} \in  {\rm DS}(\Gamma)$.
For the signature $\Sigma'=\Sigma \cup \{ (\Gamma, S, \bar{i})\}$, 
\begin{itemize}
\item[(i)] If $(b: S(\bar{a}) \; (\Delta)) \in {\cal J}(\Sigma')$, then the judgement
is the result of the assumption rule and $b$ is a variable $x$ and 
$$\Delta = \Delta_1, x:S(\bar{a}), \Delta_2,$$ 
\item[(ii)] If $(\Delta_1, x:S(\bar{a}), \Delta_2\; {\rm context}) \in {\cal J}(\Sigma')$,
then there is $\bar{b}$ with $\bar{a} =\bar{b}_{\bar{i}}$ and
$$(\bar{b}: \Delta_1 \to \Gamma) \in {\cal J}(\Sigma').$$
\end{itemize}
$\qed$
\end{lemma}}

\begin{proposition} \label{uniqder} For regular signatures $\Sigma$ , each judgement  in
${\cal J}(\Sigma)$ has a unique derivation. 
\end{proposition}
{\bf \flushleft Proof.} If the signature is regular, then the whole
sequence $(a_1,\ldots,a_n)$ can be read off from the conclusion of the rules
(R4) and (R5). $\qed$

\medskip
To obtain a decidability result for correctness of judgement expressions we require that a signature $\Sigma$ can be {\em indexed} in the sense that there is a function 
$$I: F \cup T \to \Sigma \cup \{ *\}$$ so that
$$D \in \Sigma \Longleftrightarrow (\exists s \in F \cup T) D=I(s),$$
and ${\rm decl}(I(s)) = s$ if $I(s)\ne *$.

\medskip
\begin{proposition} \label{typechk}
Suppose that $\Sigma$ is an indexed regular signature.
Then there is an algorithm which decides whether a judgement expression belongs to
${\cal J}(\Sigma)$.
\end{proposition}
{\bf \flushleft Proof.} Proposition \ref{uniqder} indicates a  deterministic procedure.
By the indexing we can find whether a symbol has a declaration, and if so, what
are its components. $\qed$

\medskip
We believe that Propositions \ref{uniqder} and \ref{typechk} are true for {\em any} indexed signatures, but shall not pursue this question here.

\medskip
For arbitrary contexts $\Gamma$ and $\Delta$ define $\Gamma \le \Delta$ to hold if every variable declaration
$x:A$ in $\Gamma$ is also in $\Delta$.  Now the following is evident.

\begin{lemma} Consider derivations relative to a fixed signature.
Suppose $\vdash_k \Gamma$ and $\vdash_d \Delta$ are contexts with $\Gamma \le \Delta$. If $\Gamma = x_1:A_1,\ldots,x_n:A_n$,
then $$\vdash_{k,d, d+1,\ldots,d+1} (x_1,\ldots,x_n): \Delta \to \Gamma$$
 is a context map.
\end{lemma}
{\flushleft \bf Proof.} For each $i=1,\ldots,n$ we have by (R3) and $\Gamma \le \Delta$
$$\vdash_{d+1} x_i:A_i \; (\Delta).$$
Now since $A_i  \equiv A_i[x_1,\ldots,x_{i-1}/x_1,\ldots,x_{i-1}]$ we get
 $$\vdash_{k,d, d+1,\ldots,d+1} (x_1,\ldots,x_n): \Delta \to \Gamma$$
as desired. $\qed$

\medskip
We show that elements have unique types up to syntactic equality. 
Moreover it is shown that the positions indicated by the top-most variables determine
the values of the other positions in context maps.

\begin{theorem} \label{UniqType}
 (Unique typing lemma)  Consider derivations relative to a fixed signature $\Sigma$.
  Let $\Delta = y_1:B_1,\ldots, y_n: B_n$ be a context.
\begin{itemize}
\item[(a)] If $a:A \; (\Delta)$ and $a:A' \; (\Delta)$, then $A \equiv A'$
\item[(b)] If $(a_1,\ldots,a_n), (b_1,\ldots,b_n): \Gamma \rightarrow \Delta$ are contexts maps, with $a_i \equiv b_i$
for each $y_i \in {\rm TV}(\Delta)$, then $(a_1,\ldots,a_n) \equiv (b_1,\ldots,b_n)$.
\end{itemize}
\end{theorem}
{\flushleft \bf Proof.} See Appendix A.  $\qed$

\begin{remark} {\em
By this result  the hidden arguments of a function or type declaration, as in (\ref{Sdecl}) and (\ref{fdecl}), are unique if  they exists.
Indeed if 
$$S(a_{i_1},\ldots,a_{i_k}) \; {\rm type} \; (\Delta)$$
has been derived by R4 in two ways using substitutions $(a_1,\ldots,a_n)$ and $(b_1,\ldots,b_n)$ with $a_{i_j} =b_{i_j}$ for
all $j=1,\ldots,k$. Then $(a_1,\ldots,a_n) \equiv (b_1,\ldots,b_n)$.
Similarly suppose that 
$$f(a_{i_1},\ldots,a_{i_k}): U[a_1,\ldots,a_n/x_1,\ldots,x_n] \;(\Delta)$$
 and 
$$f(b_{i_1},\ldots,b_{i_k}): U[b_1,\ldots,b_n/x_1,\ldots,x_n] \;(\Delta)$$
had been derived using R5 and
that  $f(a_{i_1},\ldots,a_{i_k})= f(b_{i_1},\ldots,b_{i_k})$. Then again $(a_1,\ldots,a_n) \equiv (b_1,\ldots,b_n)$.
In either case, the remaining arguments can be obtained uniquely from $a_{i_1},\ldots,a_{i_k}$.}
\end{remark}

\begin{remark} \label{decision}
{\em {\em Decision problems.} The following basic decision problems are of interest. Suppose $\Sigma$
is an indexed signature and that $(\Gamma \; {\rm context}) \in {\cal J}(\Sigma)$. 
\begin{itemize}
\item[] (Correctness of types) For a pretype $A$, decide whether 
$(A \; {\rm type}\; (\Gamma))\in {\cal J}(\Sigma).$
\item[] (Type checking) For a pretype $A$ and a preelement $a$, decide whether 
\begin{equation} \label{ainA}
(a: A \; (\Gamma))\in 
{\cal J}(\Sigma)
\end{equation}
holds.
\item[] (Type inference)  For a preelement $a$, is there a pretype $A$ such that (\ref{ainA}) holds?
\item[] (Type inhabitation) For a pretype $A$, is there a preelement $a$ such that (\ref{ainA})?
\end{itemize}
The semigroup word problem can be encoded as a type inhabitation problem (cf.\ Example \ref{semigp}). Therefore the last problem cannot be expected to be decidable.
We sketch the decision procedures for the remaining problems:

Correctness of types: The pretype $A$ has the form $S(\bar{b})$ for some $\bar{b}$ and we must have some $(\Delta, \bar{i},S) \in \Sigma$ with $|\bar{i}| = |\bar{b}|$, otherwise 
the type is not correct. Furthermore we must find $(\bar{a}: \Gamma \to \Delta) \in {\cal J}(\Sigma)$ so that $\bar{a}_{\bar{i}} = \bar{b}$. We know 
$(\Gamma \; {\rm context}) \in {\cal J}(\Sigma)$ and 
$(\Delta \; {\rm context}) \in {\cal J}(\Sigma)$, so writing $\Delta=y_1:B_1,\ldots,y_m:B_m$ need to find $\bar{a} =a_1,\ldots,a_m$ with 
\begin{equation} 
\begin{array}{lll} 
 &  a_1: B_1 \; (\Gamma) \\ 
&  a_2: B_2[a_1/y_1] \;(\Gamma) \\
&\quad \vdots \\
& a_m: B_m[a_1,\ldots,a_{m-1}/y_1,\ldots,y_{m-1}] \;(\Gamma)
 \end{array}
\end{equation}
The sequence $\bar{i}$ contains the positions of the topmost variables in $\Delta$, and
in particular $m$, so $a_m$ is determined. By induction, we can then find the type $a_m$
say $C$ such that $a_m:C \; (\Gamma)$. Suppose that $j_1,\ldots,j_k$ is the increasing sequence $k\ge 0$ such that $y_{j_1},\ldots,y_{j_k}$ are exactly the variables that occur in
$B_m$. Thus we must have for some  $a_{j_1},\ldots,a_{j_k}$ that
$$C=B_m[a_{j_1},\ldots,a_{j_k}/y_{j_1},\ldots,y_{j_k}].$$
If $k=0$, and $m>0$, $m-1$ must be in $\bar{i}$ and we can continue the procedure with 
$a_{m-1}$. If $k>0$, then $a_{j_1},\ldots,a_{j_k}$ have been determined, and the procedure continues with these. This is repeated until all of $\bar{a}$ has been determined.

Type inference: Let $a$ be given. If $a$ is a variable, then the type of $a$ is $A$ only if $a:A$ occurs in $\Gamma$.  Suppose that $a = f(\bar{b})$. The we must have
some $(\Delta,f,\bar{i},U) \in \Sigma$ such that $|\bar{i}| = |\bar{b}|$ otherwise $a$ cannot be typed. As in the correctness of types problem above we check if there is
$(\bar{a}:\Gamma \to \Delta) \in {\cal J}(\Sigma)$ so that $\bar{b} = \bar{a}_{\bar{i}}$.
Then type of $a$ must be $U[\bar{a}/\Delta]$.

Type checking: Let $a$ and $A$ be given. We check using type inference whether $a$ has some type $B$, and then check if $A=B$.
}
\end{remark}

\subsection{Structural theorems} \label{strthms}

The following structural results about the type system are rather standard and occur in one form or another in (Cartmell 1978) and (Belo 2008). We postpone the detailed proofs
for most of them to the Appendix.

We establish a standard substitution lemma as well as some structural rules (weakening, strengthening and interchange) for the type system.

\begin{theorem}(Substitution lemma) \label{substitution} Consider derivations in ${\cal J}(\Sigma)$ for a fixed signature $\Sigma$. Suppose that $\vdash_{\bar{\ell}} \bar{s}:\Theta \to \Gamma$.
\begin{itemize}
\item[(a)] If $\vdash_k B \; {\rm type}\; (\Gamma)$, then 
$\vdash_{k+{\rm max}(\bar{\ell})} B[\bar{s}/\Gamma] \; {\rm type}\; (\Theta)$.
\item[(b)] If $\vdash_k b: B\; (\Gamma)$, then 
$\vdash_{k+{\rm max}(\bar{\ell})} b[\bar{s}/\Gamma] :B[\bar{s}/\Gamma] \;  (\Theta)$.
\end{itemize}
\end{theorem}
{\flushleft \bf Proof.}  See Appendix A. $\qed$

\begin{theorem}(Weakening lemma) \label{weakening} Consider derivations in ${\cal J}(\Sigma)$ for a fixed signature $\Sigma$ with unrestricted variable system. Suppose that 
$\vdash B \; {\rm type}\;(\Gamma)$. Let $y$ be a variable not in ${\rm V}(\Gamma, \Theta)$. 
\begin{itemize}
\item[(a)] If $\vdash \Gamma, \Theta \; {\rm context}$, then $\vdash \Gamma, y:B, \Theta \; {\rm context}$ 
\item[(b)] If $\vdash A \; {\rm type}\; (\Gamma,\Theta)$, then $\vdash A \; {\rm type}\; (\Gamma,y:B, \Theta)$
\item[(c)] If $\vdash a:A\; (\Gamma,\Theta)$, then $\vdash a:A \; (\Gamma,y:B, \Theta)$
\end{itemize}
\end{theorem}
{\flushleft \bf Proof.}  See Appendix A. $\qed$

\medskip
Assuming that signatures are regular we can reverse this process:

\begin{theorem}(Strengthening lemma) \label{strengthening}
Consider derivations in ${\cal J}(\Sigma)$ for a fixed signature $\Sigma$ with unrestricted variable system.
Let $y$ be a variable not in ${\rm V}(\Theta, A, a)$. 
\begin{itemize}
\item[(a)] If $\vdash \Gamma, y:B, \Theta \; {\rm context}$, then $\vdash \Gamma, \Theta \; {\rm context}$. 
\item[(b)] If $\vdash A \; {\rm type}\; (\Gamma,y:B, \Theta)$, then  $\vdash A \; {\rm type}\; (\Gamma,\Theta) $.
\item[(c)] If $\vdash a:A \; (\Gamma,y:B, \Theta)$, then $\vdash a:A\; (\Gamma,\Theta)$.
\end{itemize}
\end{theorem}
{\flushleft \bf Proof.}  See Appendix A. $\qed$

\medskip
For regular signatures,  variable declarations may be interchanged provided they satisfy a simple
syntactic condition.

\begin{theorem}(Interchange lemma) \label{interchange}
Consider derivations in ${\cal J}(\Sigma)$ for a fixed signature $\Sigma$ with unrestricted variable system. Suppose that $x \notin {\rm V}(C)$. 
\begin{itemize}
\item[(a)] If $\vdash \Gamma,x;B, y:C, \Theta \; {\rm context}$, then $\vdash \Gamma,y:C, x:B, \Theta \; {\rm context} $. 
\item[(b)] If $\vdash A \; {\rm type}\; (\Gamma,x:B, y:C, \Theta) $, then  $\vdash A \; {\rm type}\; (\Gamma,y:C, x:B,\Theta)$.
\item[(c)] If $\vdash a:A \; (\Gamma, x:B,y:C, \Theta) $, then $\vdash a:A\; (\Gamma,y:C, x:B,\Theta)$.
\end{itemize}
\end{theorem}
{\flushleft \bf Proof.}  See Appendix A. $\qed$

\subsection{Presuppositions}

Two of the judgement forms have {\em presuppositions.} Define for judgement expressions:
\begin{itemize}
\item $A \; {\rm type}\; (\Gamma)$ presupposes that $\Gamma \; {\rm context}$,
\item $a: A\; (\Gamma)$ presupposes that $A\; {\rm type}\; (\Gamma)$, and 
$\Gamma \; {\rm context}$.
\end{itemize}

\begin{theorem} \label{presup} Consider derivations in ${\cal J}(\Sigma)$ for a fixed signature $\Sigma$. 
\begin{itemize}
\item[(a)] If $\vdash_k A \; {\rm type}\; (\Gamma)$, then $\vdash_{k-1} \Gamma\; {\rm context}$.
\item[(b)] If $\vdash_k \Gamma, x: A\; {\rm context}$, then $\vdash_{k-1} A \; {\rm type}\; (\Gamma)$. 
\item[(c)] If $\vdash_k \Gamma, x: A\; {\rm context}$, then 
$\vdash_{k-1} \Gamma\; {\rm context}$.
\item[(d)] If $\vdash_k x:A \; (\Gamma)$ and $x$ is a variable, then $\vdash_{k-2} A \; {\rm type}\; (\Gamma)$.
\item[(e)] If $\vdash_k f(\bar{b}):A \; (\Gamma)$, then $\vdash_{k-1} A \; {\rm type}\; (\Gamma)$.
\item[(f)] If $\vdash_k a:A \; (\Gamma)$, then $\vdash_{k-1} A \; {\rm type}\; (\Gamma)$.
\end{itemize} 
\end{theorem}
{\flushleft \bf Proof.} Cases (a) -- (c) and (e) are immediate as the presuppositions are found
in the premisses of the last applied rule.  

Case (d): only rule (R3) can achieve
$\vdash_k x:A \; (\Gamma)$, so $\vdash_{k-1} \Gamma\; {\rm context}$. Using (a) and (b)
at least once we get  $\vdash_{k-2} A \; {\rm type}\; (\Gamma)$. 

Case (f): follows from (d) and (e). $\qed$

\dontshow{
\medskip
What happens for more general signatures? Consider applications of (R4) and (R5) for regular declarations in ${\cal J}(\Sigma)$. If 
$$\vdash_k S(\bar{a}) \; {\rm type} \; (\Gamma)$$
then $k=1+\max(\bar{m})$, where using the assumptions of (\ref{cmpjudgement})
$$\vdash_{\bar{m}} \bar{a}: \Gamma \to \Delta.$$
Thus
$$\vdash_{m_1} \Gamma \; {\rm context}$$
and
$$\vdash_{m_2} \Delta \; {\rm context}$$
and furthermore  for $k=1,\ldots,n$,
$$\vdash_{m_{2+k}} a_k: A_k[\bar{a}^{k-1}/\Gamma^{k-1}] \; (\Delta).$$
}

\subsection{Variations}

In the type system of (Belo 2008) the (R5) rule is instead
$$\infer[({\rm R5^*}) \quad \begin{array}{l} (\Gamma,f,\bar{i},U) \mbox{ in  } \Sigma \\
 \end{array}]{ f(\bar{a}_{\,\bar{i}}) : U[\bar{a}/\Gamma] \;(\Delta) }{ 
 \bar{a}: \Delta \to\Gamma  }$$
 (However, that system uses only regular declarations.) 
 Define ${\cal J}^*(\Sigma)$ just as ${\cal J}(\Sigma)$ but using the (R5*) rule instead of
 (R5). We have evidently ${\cal J}(\Sigma) \subseteq {\cal J}^*(\Sigma)$, since (R5*) has fewer premisses than (R5). Towards proving the reverse inclusion we have the following results.
 
 \begin{lemma}(Substitution lemma) \label{substitution_star} Consider derivations in ${\cal J}^*(\Sigma)$ for a fixed signature $\Sigma$. Suppose that $\vdash_{\bar{m}} \bar{s}:\Theta \to \Gamma$.
\begin{itemize}
\item[(a)] If $\vdash_k B \; {\rm type}\; (\Gamma)$, then 
$\vdash_{k+{\rm max}(\bar{m})} B[\bar{s}/\Gamma] \; {\rm type}\; (\Theta)$.
\item[(b)] If $\vdash_k b: B\; (\Gamma)$, then 
$\vdash_{k+{\rm max}(\bar{m})} b[\bar{s}/\Gamma] :B[\bar{s}/\Gamma] \;  (\Theta)$.
\end{itemize}
\end{lemma}
{\flushleft \bf Proof.}  Similar to the proof of Theorem \ref{substitution}. $\qed$

\begin{lemma}(Substitution lemma) \label{substitution_star} Consider derivations in ${\cal J}^*(\Sigma)$ for a fixed signature $\Sigma$. Let $R$ be an application of rule R5* with respect to a fixed function constant. Suppose that $\vdash^R_{\bar{m}} \bar{s}:\Theta \to \Gamma$.
\begin{itemize}
\item[(a)] If $\vdash^R_k B \; {\rm type}\; (\Gamma)$, then 
$\vdash^R_{k+{\rm max}(\bar{m})} B[\bar{s}/\Gamma] \; {\rm type}\; (\Theta)$.
\item[(b)] If $\vdash^R_k b: B\; (\Gamma)$, then 
$\vdash^R_{k+{\rm max}(\bar{m})} b[\bar{s}/\Gamma] :B[\bar{s}/\Gamma] \;  (\Theta)$.
\end{itemize}
\end{lemma}
{\flushleft \bf Proof.}  As above, but we get an overestimation of the height. $\qed$

\begin{lemma} \label{presup_star} Consider derivations in ${\cal J}^*(\Sigma)$ for a fixed signature $\Sigma$. 
\begin{itemize}
\item[(a)] If $\vdash_k A \; {\rm type}\; (\Gamma)$, then $\vdash_{k-1} \Gamma\; {\rm context}$.
\item[(b)] If $\vdash_k \Gamma, x: A\; {\rm context}$, then $\vdash_{k-1} A \; {\rm type}\; (\Gamma)$. 
\item[(c)] If $\vdash_k \Gamma, x: A\; {\rm context}$, then 
$\vdash_{k-1} \Gamma\; {\rm context}$.
\item[(d)] If $\vdash_k x:A \; (\Gamma)$ and $x$ is a variable, then $\vdash_{k-2} A \; {\rm type}\; (\Gamma)$.
\item[(e)] If $\vdash_k f(\bar{b}):A \; (\Gamma)$, then $$\vdash_{u+k-1} A \; {\rm type}\; (\Gamma)$$
where $A=U[\bar{a}/\Delta]$, $(\Delta,f,\bar{i},U) \in \Sigma$ and $\vdash_u U\; {\rm type}\; (\Delta)$, $\vdash_{k-1} \bar{a}: \Gamma \to \Delta$ and $\bar{b} = \bar{a}_{\bar{i}}$.
\item[(f)] Suppose $\Sigma$ is finite. Then there is $u \ge 0$, such that 
if $\vdash_k a:A\; (\Gamma)$, then 
$$\vdash_{k+u} A\; {\rm type} \; (\Gamma).$$
\end{itemize} 
\end{lemma}
{\flushleft \bf Proof.} Cases (a) -- (c) are immediate as the presuppositions are found
in the premisses of the last applied rule. 

Case (d): only rule (R3) can achieve
$\vdash_k x:A \; (\Gamma)$, so $\vdash_{k-1} \Gamma\; {\rm context}$. Using (a) and (b)
at least once we get  $\vdash_{k-2} A \; {\rm type}\; (\Gamma)$.

Case (e): Note that $\vdash_k f(\bar{b}):A \; (\Gamma)$ can only be derived by an application of (R5*). Thus there are $\vdash_{\bar{\ell}} \bar{a}: \Gamma \to \Delta$ and
 $(\Delta,f,\bar{i},U) \in \Sigma$ with  $A=U[\bar{a}/\Delta]$. By the
correctness requirement on signatures
$$\vdash_u U\; {\rm type}\; (\Delta)$$
for some $u$. By the Substitution Lemma we get
$$\vdash_{u+{\rm max}(\bar{\ell})} U[\bar{a}/\Delta]\; {\rm type}\; (\Gamma).$$
Noting that ${\rm max}(\bar{\ell}) \le k-1$ finishes the proof. 

Case (f): this follows from (e) by taking the maximum of the $u$s appearing for
function declarations in $\Sigma$.
$\qed$

\medskip
Introduce the following notation $\Sigma \vdash {\cal U}$ for ${\cal U} \in {\cal J}(\Sigma)$
and $\Sigma \vdash^* {\cal U}$ for ${\cal U} \in {\cal J}^*(\Sigma)$.

\begin{lemma} \label{expbyfunct}
Suppose that $\Sigma$ is a signature such that ${\cal J}^*(\Sigma) = {\cal J}(\Sigma)$. Assume that $D=(\Delta,f,\bar{i},U)$ is a predeclaration, where $f$ is not declared in $\Sigma$, and where $(U \; {\rm type}\; (\Delta)) \in {\cal J}^*(\Sigma)$. Then for
$\Sigma' = \Sigma \cup \{D\}$
$${\cal J}^*(\Sigma') = {\cal J}(\Sigma').$$
\end{lemma}
{\flushleft \bf Proof. }  See Appendix A.  $\qed$

\medskip
\begin{lemma} \label{expbytype} Suppose that $\Sigma$ is a signature such that ${\cal J}^*(\Sigma) = {\cal J}(\Sigma)$. Assume that $D=(\Delta,\bar{i},S)$ is a predeclaration, where $S$ is not declared in $\Sigma$, and where $(\Delta \; {\rm context}) \in {\cal J}^*(\Sigma)$. Then for
$\Sigma' = \Sigma \cup \{D\}$
$${\cal J}^*(\Sigma') = {\cal J}(\Sigma').$$
\end{lemma}
{\flushleft \bf Proof.}  See Appendix A.  $\qed$

\medskip
As a consequence of these lemmas we have: 

\begin{theorem} If $\Sigma$ is a finite, inductive signature, then ${\cal J}(\Sigma) = {\cal J}^*(\Sigma)$.
\end{theorem}

\dontshow{
\section{FOLDS vocabularies  and signatures} \label{FOLDSsec}

A finite FOLDS {\em vocabulary} (or {\em dependent sorts vocabulary,} DSV) (Makkai 1995)  is given by
a finite, skeletal category $K$, which is one-way (has no non-trivial endomorphisms).

\begin{example} \label{F2}
{\em The category with three objects
$$T \two^s_t A \two^d_c O$$
and the six non-identity arrows $s$, $t$, $d$, $c$, $ds=dt$ and $cs=ct$, is a FOLDS vocabulary $K_2$. It is
suitable for describing a 2-category.
}
\end{example}

We show how to translate any FOLDS vocabulary $K$ to a FOLDS-like signature.
First define a binary relation $\le$  on the objects by
$$A \le B \Longleftrightarrow_{\rm def} \text{there is a morphism from $B$ to $A$}.$$
By the properties of $K$ this is a well-founded partial order.  Let $\le^*$ be some linear
order that extends $\le$. Note that $A \le^* B$ implies that there can be no non-identity
morphism $A \to B$. (Otherwise, we have $A=B$ and a nontrivial endomorphism on $A$).
Construct for every object $A$
a non-repeating enumeration 
$$x^A_1,\ldots,x^A_{n(A)}$$ 
of all the non-identity morphisms in $K$
with domain $A$, which is such that
$${\rm cod}(x^A_i) \le^* {\rm cod}(x^A_j) \text{ whenever } i < j.$$
We now define a signature where the objects of $K$ become the type symbols and
where the morphisms of $K$ serve as variables in the declarations.  For this purpose 
we assume that $K \subseteq V$ for some infinite discrete set $V$ with fop 
$(\varphi_{\infty}, {\sf fr}_{\infty})$.

If $\bar{u}=u_1,\ldots,u_n$ is a sequence of elements in $K$ and $v\in K$, 
with ${\rm dom}(u_i) = {\rm cod}(v)$, for all $i$,  we write $\bar{u}^v$ for
$$u_1v,\ldots,u_nv.$$

For each object $A$ of $K$ we define a type pre-declaration (regular)
$${\cal D}_A =_{\rm def} (\Gamma_A,A)$$
by induction on $\le$.
Let
\begin{equation} \label{ctxA}
\Gamma_A = x_1: A_1,\ldots, x_n: A_n
\end{equation}
where $x_i = x_i^{A}$,  and $$A_i = C_i({\rm OV}( \Gamma_{C_i})^{x_i})$$ with $C_i = {\rm cod}(x_i)$ and $n=n(A)$.

If $$B_1 <^* B_2 <^* \cdots <^* B_N$$ are all the objects of $K$, then the pre-signature corresponding to the vocabulary $K$ is
$$\Sigma_K = \{{\cal D}_{B_1},\ldots,{\cal D}_{B_N}\}.$$

Returning to Example \ref{F2} above, we start by choosing order and enumerations
\begin{itemize}
\item $O <^* A <^* T$
\item $x_1^T, x_2^T, x_3^T, x_4^T= ds, cs,  s, t$
\item $x_1^A, x_2^A= c, d$
\item $n(O)=0$.
\end{itemize}
We get the following contexts:
\begin{itemize}
\item $\Gamma_O = \langle \rangle$
\item $\Gamma_A = \langle c:O, d: O\rangle = \langle x_1:O, x_2:O \rangle$
\item \begin{eqnarray*} 
     \Gamma_T &=&  \langle ds:O, cs:O, s:A({\rm OV}(\Gamma_A)^s), t:A({\rm OV}(\Gamma_A)^t) \rangle \\
  &=& \langle ds:O, cs:O, s:A(cs,ds), t:A(ct,dt) \rangle \\
  &=& \langle ds:O, cs:O, s:A(cs,ds), t:A(cs,ds) \rangle \\
   &=& \langle x_1:O, x_2:O, x_3:A(x_2,x_1), x_4:A(x_2,x_1) \rangle
         \end{eqnarray*}
\end{itemize}
The resulting signature is 
$$\Sigma_{K_2}=\{(\Gamma_O,O),(\Gamma_A,A),(\Gamma_T,T) \}.$$
Note that with an alternate choice of the enumerations the order of the arguments
can be different. 

It can be proved that $\Sigma_K$ is a signature for an arbitrary $K$:

\begin{theorem} \label{F2signature}
For any FOLDS vocabulary $K$, the pre-signature $\Sigma_K$ is a signature. 
\end{theorem}
{\flushleft \bf Proof.} See Appendix A. $\qed$

\medskip 
FOLDS vocabularies give a canonical way of maximally hiding variables
for the corresponding signature.
A morphism $f$ of a FOLDS vocabulary $K$ is called {\em reducible} there are
two non-identity arrows $g$ and $h$ in $K$ such that $f=gh$.

\begin{proposition} For a FOLDS vocabulary $K$, and each object $A$ of $K$,
the top variables ${\rm TV}(\Gamma_A)$ are exactly the non-reducible arrows in $K$ 
with domain $A$.
\end{proposition}
{\flushleft \bf Proof.} With reference to the context (\ref{ctxA}), note that $x_i$ is not a top variable if, and only if ,$x_i$ occurs in some $A_j$ with $j>i$.  If $x_i$ occurs in $A_j$, $j>i$, then $x_i$ is reducible. Conversely, if $x_i$ is reducible, then  $x_i=ux_j$, for some $x_j$ and some non-identity $u$. Now ${\rm cod}(x_j) = C_j$, so $u$ is variable in $\Gamma_{C_j}$ by construction. Thus
$ux_j$ occurs in $({\rm OV}(\Gamma_{C_j}))^{x_j}$, and indeed $x_i$ occurs in $A_j$. Now since
$x_i=ux_j$ we must have $j\ne i$, and $j<i$ would imply ${\rm cod}(x_j) \le^* {\rm cod}(x_i)$,
saying that there can be no non-identity map from  ${\rm cod}(x_j)$ to ${\rm cod}(x_i)$, which
however $u$ is an example of. Hence $j>i$.  $\qed$

\medskip
Consider the other direction, from signatures to vocabularies.
Suppose a FOLDS-like inductive signature $\Sigma=\{{\cal D}_1,\ldots,{\cal D}_N\}$ on regular form is given. Thus one can write
$${\cal D}_i =(\Gamma_i,S_i),$$
$n(i) = {\rm length}({\rm OV}(\Gamma_i))$,
and without any loss of generality one can assume that variables are chosen as follows
$$\Gamma_i = x^i_1:A_1^i,\ldots,x^i_{n(i)}:A_{n(i)}^i,$$
and so that $x^i_j = x^{i'}_{j'}$ implies $i=i'$ and $j=j'$, and
 moreover that for $j=1,\ldots, n(i)$,
 $$A_j^i = S_{f(i,j)}(x^i_{g(i,j,1)},\ldots,x^i_{g(i,j, n(f(i,j)))}).$$
Here ${\rm dom}(f) = \{(i,j) : i=1,\ldots,N, j =1,\ldots, n(i)\}$.
By the inductiveness of the signature we can assume
that $1 \le f(i,j) < i$, so that $S_{f(i,j)}$ has been declared before $S_i$. 
We also have that  $1 \le g(i,j,k) < i$.

Define a finite category $K_\Sigma$ as follows. Its objects are $\{S_1,\ldots,S_N\}$, the type symbols declared in the signature.
Define first a directed acyclic graph with multiple edges. Its vertices are the objects above. For every $i=1,\ldots,N$ and 
$j=1,\ldots,n(i)$ we declare 
$x^i_j$ to be an arc from $S_i$ to $S_{f(i,j)}$. Now generate the free category on this graph. 
A typical path starting at $S_i$ looks like this
$$x_{j_n}^{f( \cdots f(f(f(i,j_1),j_2),j_3), \ldots, j_{n-1})} \cdots x_{j_3}^{f(f(i,j_1),j_2)} x_{j_2}^{f(i,j_1)}x_{j_1}^i$$
The following relations then generate the equivalence relation on paths
$$x_k^{f(i,j)}x_j^i = x^i_{g(i,j,k)} \qquad (1 \le j \le n(i), 1 \le k \le n(f(i,j))).$$
This defines the category $K_\Sigma$. It is skeletal and contains no nontrivial endomorphisms since the free category has these properties.

We summarize this as a

\begin{theorem} \label{Signature2F}
For any finite inductive FOLDS-like signature $\Sigma$, the category $K_\Sigma$ is
a FOLDS vocabulary. $\qed$
\end{theorem}
}

\section{The category of contexts} \label{catcontsec}

\subsection{Categories with families}

Categories with families  (Dybjer 1996) is one of several equivalent ways of defining categorical semantics for dependent type theories (Hofmann 1997). We recall:

\begin{definition} \label{cwfdef}
{\em 
A {\em category with families (cwf) }consists of the following data
\begin{itemize}
\item[(a)] A category ${\cal C}$ with a terminal object $\top$. This is thought of as  the category of contexts and substitutions.  For $\Gamma \in {\cal C}$ denote by $\epsilon_\Gamma$ the unique morphism $\Gamma \to \top$.
\item[(b)] For each  object $\Gamma$ of ${\cal C}$, a class $\Ty(\Gamma)$ and for each morphism $f: \Delta \to \Gamma$, 
a class function $\Ty(f): \Ty(\Gamma) \to \Ty(\Delta)$, for which use the notation $A\{f\}$ for $\Ty(f)(A)$, to
suggest that it is the result of performing the substitution $f$ in the type $A$.
These functions should satisfy, for all $A \in \Ty(\Gamma)$, and 
$g:\Theta \to \Delta$, $f:\Delta \to \Gamma$:
\begin{itemize}
\item[(i)] $A\{1_\Gamma\} = A$,
\item[(ii)] $A\{f \circ g\}=A\{f\}\{g\}$.
\end{itemize} 
(Thus $\Ty$ may be regarded as a functor ${\cal C}^{\rm op} \to {\rm Class}$.)
\item[(c)] For each $A \in \Ty(\Gamma)$, an object $\Gamma.A$ in ${\cal C}$ and a
morphism ${\rm p}(A) = {\rm p}_\Gamma(A): \Gamma.A \to \Gamma$.  This tells that each context can be extended
by a type in the context, and that there is a projection from the extended context to the original one.
\item[(d)] For each $A \in \Ty(\Gamma)$, there is a class $\Tm(\Gamma,A)$ --- thought of
as the terms of type $A$. It should be 
 such that for $f: \Delta \to \Gamma$ there is
a class function $\Tm(f): \Tm(\Gamma,A) \rightarrow \Tm(\Delta,A\{f\})$,  where we write $a\{f\}$ for $\Tm(f)(a)$. It should satisfy the following
\begin{itemize}
\item[(i)] $a\{1_{\Gamma} \} = a$ for $a \in \Tm(\Gamma,A) $ $(= \Tm(\Gamma, A\{1_\Gamma\})).$
\item[(ii)] $a\{f \circ g\} = a\{f\}\{g\}$ for $a \in \Tm(\Gamma,A)$

(Note: $\Tm(\Theta, A\{f \circ g\}) = \Tm(\Theta,A\{f\}\{g\})$.)
\end{itemize}
\item[(e)]  For each $A \in \Ty(\Delta)$ there is an element ${\rm v}_{A} \in \Tm(\Delta.A, A\{{\rm p}(A)\})$.
\item[(f)] For any morphism $f:\Gamma \to \Delta$ and $a \in \Tm(\Gamma, A\{f\})$, there is 
$$\langle f,a\rangle_A: \Gamma \to \Delta.A.$$
This construction should satisfy
\begin{itemize}
\item[(i)] ${\rm p}(A)\circ \langle f,a\rangle_A = f$,
\item[(ii)] ${\rm v}_A\{\langle f,a\rangle_A\} = a$,  

(Note: ${\rm v}_A\{\langle f,a\rangle_A\}  \in \Tm(\Gamma, A\{{\rm p}(A)\}\{\langle f,a\rangle_A\})= \Tm(\Gamma, A\{f\}))$
\item[(iii)] $\langle {\rm p}(A) \circ h,{\rm v}_A\{h\}\rangle_A = h$ for any $h:\Gamma \to \Delta.A$.
\item[(iv)] for any $g:\Theta \to \Gamma$,
\begin{equation} \label{substdist}
\langle f,a\rangle_A  \circ g=  \langle f \circ g,  a\{g\}\rangle_A
\end{equation}
\end{itemize}
(Remark: $a\{g\} \in \Tm(\Theta, A\{f\}\{g\} )=\Tm(\Theta,A\{f \circ g\})$.)
$\qed$
\end{itemize}
}
 \end{definition}

A cwf ${\cal C}$ is {\em separated} if for any $\Gamma, \Delta \in {\cal C}$ and $A \in {\rm  Ty}(\Gamma)$ and $B \in \Ty(\Delta)$, 
$$\Gamma.A = \Delta.B \Longrightarrow \text{$\Gamma = \Delta$ and $A=B$}.$$
A cwf ${\cal C}$ is {\em contextual} if for every $\Gamma \in {\cal C}$ there is a unique sequence $A_1 \in \Ty(\top)$, $A_2 \in \Ty(\top.A_1)$,\ldots, 
$A_n \in \Ty(\top.A_1. \ldots.A_{n-1})$, $n\ge 0$ such that
\begin{equation} \label{finctx}
\Gamma = \top.A_1. \ldots.A_n.
\end{equation}

\subsection{Structure of morphisms in cwfs}

Let $\Theta \in {\cal C}$. There is a unique morphism $\varepsilon_\Theta: \Theta \to \top$.
For $a_1 \in \Tm(\Theta, A_1\{\varepsilon_\Theta\})$, we have
$$\langle \varepsilon_\Theta, a_1 \rangle_{A_1} : \Theta \to \top.A_1.$$
For $a_2 \in \Tm(\Theta, A_2\{\langle \varepsilon_\Theta, a_1 \rangle_{A_1}\})$, we get
$$\langle \langle \varepsilon_\Theta, a_1 \rangle_{A_1}, a_2\rangle_{A_2} : \Theta \to \top.A_1.A_2.$$
For $a_3 \in \Tm(\Theta, A_3\{\langle \langle \varepsilon_\Theta, a_1 \rangle_{A_1}, a_2\rangle_{A_2}\})$, we get
$$\langle \langle \langle \varepsilon_\Theta, a_1 \rangle_{A_1}, a_2\rangle_{A_2} a_3\rangle_{A_3}: \Theta \to \top.A_1.A_2.A_3$$
We introduce the notation 
$$[a_1, \ldots, a_n]_{\Theta;A_1,\ldots,A_n}= \langle \cdots \langle \langle \varepsilon_\Theta, a_1 \rangle_{A_1}, a_2\rangle_{A_2}, \ldots, a_n\rangle_{A_n},$$
which satisfies the recursive equations
\begin{eqnarray} \label{req}
 []_{\Theta} &=&  \varepsilon_{\Theta}, \\
{[a_1, \ldots, a_{n+1}]}_{\Theta;A_1,\ldots,A_{n+1}} 
&= &
\langle [a_1, \ldots, a_n ]_{\Theta;A_1,\ldots,A_n}, a_{n+1}\rangle_{A_{n+1}}.
\end{eqnarray}
Dropping subscripts and curly braces in applications of substitution
 we can write
 $$B[a_1,\ldots,a_n] = B\{[a_1, \ldots, a_n]_{\Theta;A_1,\ldots,A_n}\}$$
 and
 $$b[a_1,\ldots,a_n] =  b\{[a_1, \ldots, a_n]_{\Theta;A_1,\ldots,A_n}\}.$$
The above is then more succinctly rendered as:
For $a_1 \in \Tm(\Theta, A_1[])$, 
$$[a_1]: \Theta \to \top.A_1.$$
For $a_2 \in \Tm(\Theta, A_2[a_1])$, we get
$$[a_1, a_2] : \Theta \to \top.A_1.A_2.$$
For $a_3 \in \Tm(\Theta, A_3[a_1, a_2])$, we get
$$[a_1, a_2, a_3] : \Theta \to \top.A_1.A_2.A_3$$
and in general  for $a_k \in \Tm(\Theta, A_k[a_1, a_2,\ldots, a_{k-1}])$,
\begin{equation} \label{gentupl}
[a_1,a_2, \ldots, a_k] : \Theta \to \top.A_1.A_2. \ldots . A_k.
\end{equation}


For $f: \Delta \to \Theta$ we get by repeated uses of (\ref{substdist}) 
$$[a_1,\ldots,a_n] \circ f = [a_1\{f\},\ldots,a_n\{f\}].$$
If $f=[b_1,\ldots,b_m]$, we obtain
\begin{equation} \label{substcomp}
[a_1,\ldots,a_n] \circ [b_1,\ldots,b_m] = [a_1[b_1,\ldots,b_m],\ldots,a_n[b_1,\ldots,b_m]].
\end{equation}
For contexts on the form (\ref{finctx}) we introduce the notation
\begin{eqnarray*}
{\rm p}_{\Gamma}^{(0)} &= & {\rm id}_\Gamma \\
{\rm p}_{\Gamma.A}^{(k+1)} &=& {\rm p}_{\Gamma}^{(k)} \circ {\rm p}_{\Gamma}(A).
\end{eqnarray*}
Then $${\rm p}_{\top.A_1. \ldots.A_n}^{(i)} : \top.A_1. \ldots.A_n \to \top.A_1. \ldots.A_{n-i}$$
and $${\rm p}_{\top.A_1. \ldots.A_n}^{(i)} = 
{\rm p}(A_{n-i+1}) \cdots {\rm p}(A_n).$$
For an arbitrary $f: \Theta \to\top.A_1. \ldots.A_n$ we can then write using (f.iii)
$$
\begin{aligned}
f & =  \langle {\rm p}(A_n) \circ f, {\rm v}_{A_n}\{f\} \rangle \\ 
& =  \langle \langle {\rm p}(A_{n-1}) {\rm p}(A_n)f, {\rm v}_{A_{n-1}}\{{\rm p}(A_n) f\} \rangle, {\rm v}_{A_n}\{f\} \rangle \\ 
&=  \langle \langle \langle   {\rm p}(A_{n-2}){\rm p}(A_{n-1}) {\rm p}(A_n)f, {\rm v}_{A_{n-2}}\{{\rm p}(A_{n-1}){\rm p}(A_n) f\} \rangle, {\rm v}_{A_{n-1}}\{{\rm p}(A_n)f\} \rangle, {\rm v}_{A_n}\{f\} \rangle\\
& \vdots \\
& = [{\rm v}\{{\rm p}^{(n-1)}f\}, {\rm v}\{{\rm p}^{(n-2)}f\},\ldots,{\rm v}\{{\rm p}^{(1)}f\}, {\rm v}\{{\rm p}^{(0)}f\}].
\end{aligned}
$$
Introducing the special notation
$$f_{(i)} = {\rm v}_{A_i}\{{\rm p}^{(n-i)}f\},$$
we have $$f = [f_{(1)},\ldots,f_{(n)}].$$
The identity ${\rm id}_\Gamma={\rm id}$ can be regarded as a generic element of $\Gamma = \top.A_1. \ldots.A_n$. If we introduce suggestively
 the $i$th variable as
\begin{equation} \label{xiproj}
{\sf x}_i =_{\rm def}  {\rm id}_{(i)} = {\rm v}\{{\rm p}^{(n-i)}\} \in \Tm(\top.A_1. \ldots.A_n, A_i\{{\rm p}^{(n-i+1)}\})
\end{equation}
we get
$${\rm id} = [{\sf x}_1 ,\ldots, {\sf x}_n].$$
Note that 
$${\sf x}_i\{f\}= f_{(i)}.$$
Moreover
$${\rm p}_{\Gamma}^{(i)} = [{\sf x}_1,\ldots, {\sf x}_{n-i}].$$
With the notation $${\rm v}_i = {\rm v}\{{\rm p}^{(i-1)}\}$$ we get
de Bruijn indices, and 
$${\rm v}_i = {\sf x}_{n-i+1},$$
so
$${\rm id} = [{\rm v}_n ,\ldots, {\rm v}_1]$$
and 
$${\rm p}_{\Gamma}^{(i)} = [{\rm v}_n,\ldots, {\rm v}_{i+1}].$$

\subsection{Pullbacks of canonical projections}

Pullbacks along canonical projections are particularly well behaved in cwfs, and
occur in substitutive properties of quantifiers. We define the following. For $f: \Delta \to \Gamma$ and $S\in {\rm Ty}(\Gamma)$ a morphism ${\rm q}(f,S): \Delta.S\{f\} \to \Gamma.S$
set
$${\rm q}(f, S)=_{\rm def} \langle f\circ {\rm p}_{\Delta}(S\{f\}), {\rm v}_{S\{f\}}\rangle_{\rm S}.$$
We use the suggestive notation 
$$f.S =_{\rm def} {\rm q}(f,S).$$
As is shown in (Hofmann 1997) one obtains a pullback square which behaves functorially in $f$:
\begin{equation} \label{pqpb}
\bfig\square<1000,600>[\Delta.S\{f\}`\Gamma.S`\Delta`\Gamma;f.S`{\rm p}_{\Delta}(S\{f\})`{\rm p}_{\Gamma}(S)`f]\efig
\end{equation}
These functoriality properties are
\begin{eqnarray*}
1_\Gamma.S &=& 1_{\Gamma.S} \\
(f \circ g).S &= &(f.S)\circ (g.S\{f\}).
\end{eqnarray*}
By the definition of ${\rm q}$ it follows that
\begin{equation}
{\rm v}_S\{f.S\}={\rm v}_{S\{f\}}.
\end{equation}
Further we have the equation
\begin{equation}
f.S\circ \langle g,N\rangle = \langle f \circ g, N \rangle.
\end{equation}

By composing pullback squares we have the following pullback for repeated
extensions of contexts.
\begin{equation} \label{iterpqpb}
\bfig\square<1800,600>[\Delta.S_1\{f\}.S_2\{f.S_1\}.\ldots.S_n\{f.S_1.\ldots.S_{n-1}\}`\Gamma.S_1.\ldots.S_n`\Delta`\Gamma;f.S_1.\ldots.S_n`{\rm p}^{(n)}`{\rm p}^{(n)}`f]\efig
\end{equation}

We denote a repeated extension $\Gamma.S_1.\ldots.S_n$ by $\Gamma.\Sigma$ where
$\Sigma= S_1.\ldots.S_n$. Note that $\Sigma$ is in general not a context by itself. Introduce
the notation 
$$\Sigma\{f\} = S_1\{f\}.S_2\{f.S_1\}.\ldots.S_n\{f.S_1.\ldots.S_{n-1}\}.$$
We also write for iterated projections
$${\rm p}_\Gamma(\Sigma) =  {\rm p}_\Gamma(S_1){\rm p}_{\Gamma.S_1}(S_2) \cdots {\rm p}_{\Gamma.S_1.\ldots.S_{n-1}}(S_n).$$
Note that for $n=1$ this the usual projection, and for $n=0$ it is the identity.
Then (\ref{iterpqpb}) may be written
\begin{equation} \label{iterpqpb2}
\bfig\square<1000,600>[\Delta.\Sigma\{f\} `\Gamma.\Sigma`\Delta`\Gamma;f.\Sigma`{\rm p}_\Delta(\Sigma\{f\})`{\rm p}_\Gamma(\Sigma)`f]\efig
\end{equation}
The functoriality properties also extend to this notation
\begin{itemize}
\item $1_\Gamma.\Sigma = 1_{\Gamma.\Sigma}$
\item $(f \circ g).\Sigma =(f.\Sigma)\circ (g.\Sigma\{f\})$.
\end{itemize}
Moreover introducing the notation
$${\rm v}_{A.\Sigma} =_{\rm def} {\rm v}_A\{{\rm p}(\Sigma)\}$$
we have
$${\rm v}_{A.\Sigma} \in \Tm(\Delta.A.\Sigma, A\{{\rm p}(A.\Sigma)\}).$$

\subsection{Cwf morphisms}

The cwfs may themselves be organised as a (large) category, using the notion of morphism between cwfs introduced by Dybjer (1996), which suffices here. Clairambault and Dybjer (2014) use 2-categories to establish a biequivalence between cwfs extended with certain type constructions and left exact categories extended with corresponding universal constructions.

\begin{definition} \label{defcwf}
{\em 
A {\em cwf morphism} from a  cwf ${\cal C}=({\cal C},\Ty,\Tm)$ to another cwf ${\cal C'} = ({\cal C}',\Ty',\Tm')$ is a triple $(F,\sigma, \theta)$ consisting of a functor $F:{\cal C} \rightarrow {\cal C}'$ 
such that 
\begin{equation} \label{termpres}
F(\top_{\cal C}) =\top_{{\cal C}'}
\end{equation}
and a family of functions $\sigma_\Gamma: \Ty(\Gamma) \rightarrow \Ty'(F(\Gamma))$ satisfying the naturality condition for $f: \Delta \rightarrow \Gamma$,
and $A \in \Ty(\Gamma)$,
\begin{equation} \label{cwfmorph1}
\bfig\square<1000,600>[\Ty(\Gamma)`\Ty'(F(\Gamma))`\Ty(\Delta)`\Ty'(F(\Delta)).;\sigma_{\Gamma}`\_\{f\}`\_\{F(f)\}`\sigma_{\Delta}]\efig
\end{equation}
commutes (that is $\sigma:\Ty \to \Ty'\circ F$ is a natural transformation)
and such that
$$F(\Gamma.A) = F(\Gamma).\sigma_{\Gamma}(A)$$ and 
\begin{equation} \label{ppres}
F({\rm p}_{\Gamma}(A))={\rm p}_{F(\Gamma)}(\sigma_{\Gamma}(A))
\end{equation}
and as third component a family of functions
$$\theta_{\Gamma,A}: \Tm(\Gamma,A) \rightarrow \Tm'(F(\Gamma),\sigma_{\Gamma}(A))$$
such that for $f: \Delta \rightarrow \Gamma$, $A\in \Ty(\Gamma)$ the following diagram commutes
\begin{equation} \label{cwfmorph2}
\bfig\square<1400,600>[\Tm(\Gamma,A)`\Tm'(F(\Gamma),\sigma_\Gamma(A))`\Tm(\Delta,A\{f\})`\Tm'(F(\Delta),\sigma_{\Delta}(A\{f\})).;\theta_{\Gamma,A}`\_\{f\}`\_\{F(f)\}`\theta_{\Delta,A\{f\}}]\efig
\end{equation}
It is required that for $A \in \Ty(\Gamma)$,
\begin{equation} \label{thetav}
\theta_{\Gamma.A, A\{{\rm p}(A)\}}({\rm v}_A) =  {\rm v}_{\sigma_{\Gamma}(A)}
\end{equation}
(the left hand side is in 
\begin{eqnarray*}
\Tm'(F\Gamma.FA, \sigma_{\Gamma.A}(A\{{\rm p}(A)\}))
&=& \Tm'(F\Gamma.FA, \sigma_{\Gamma}(A)\{F({\rm p}(A))\})\\
&=&  \Tm'(F\Gamma.FA, \sigma_{\Gamma}(A)\{{\rm p}(\sigma_{\Gamma}(A))\})
\end{eqnarray*}
so the right hand side has correct type)
and moreover it is required that for
$f:\Delta \to \Gamma$ and $a \in \Tm(\Delta, A\{f\})$,  
\begin{equation} \label{eq12}
F(\langle f,a\rangle_A) = \langle F(f),\theta_{\Delta, A\{f\}}(a) \rangle_{\sigma_\Gamma(A)}. 
\end{equation}
(Here the left hand side is in $F(\Delta) \to F(\Gamma).F(A)$. As $F(f): F(\Delta) \to F(\Gamma)$ and  $\theta_{\Delta, A\{f\}}(a)$ is a member of 
$$\Tm'(F(\Delta),\sigma_{\Delta}(A\{f\})) = 
\Tm'(F(\Delta),\sigma_{\Gamma}(A)\{F(f)\}).$$
So the right hand side of (\ref{eq12}) has the correct type as well.) $\qed$
}
\end{definition}

\begin{lemma} \label{ffa_lm}
Suppose that $F: {\cal C} \to {\cal C}'$ is a cwf morphism. For $A \in \Ty_{\cal C}(\Gamma)$ and $f:\Delta \to \Gamma$ in ${\cal C}$,
$$F(f.A) = F(f).\sigma_\Gamma(A).$$
\end{lemma}
{\flushleft \bf Proof.} A straightforward calculation gives
\begin{eqnarray*}
F(f.A) & = & F(\langle f \circ {\rm p}_\Delta(A\{f\}), {\rm v}_{A\{f\}} \rangle) \\
   & = & \langle  F(f \circ {\rm p}_\Delta(A\{f\}), \theta_{\Delta.A\{f\}, A\{ f \circ {\rm p}_\Delta(A\{f\})\}}({\rm v}_{A\{f\}}) \rangle_{\sigma_\Gamma(A)} \\
    & = & \langle  F(f) \circ F({\rm p}_\Delta(A\{f\}), \theta_{\Delta.A\{f\}, A\{ f\}\{{\rm p}_\Delta(A\{f\})\}}({\rm v}_{A\{f\}}) \rangle_{\sigma_\Gamma(A)} \\
    & = & \langle  F(f) \circ  {\rm p}_{F\Delta}(\sigma_\Delta(A\{f\})), {\rm v}_{\sigma_\Delta(A\{f\})} \rangle_{\sigma_\Gamma(A)} \\
    & = & \langle  F(f) \circ  {\rm p}_{F\Delta}(\sigma_\Gamma(A)\{F(f)\}), {\rm v}_{\sigma_\Gamma(A)\{F(f)\}} \rangle_{\sigma_\Gamma(A)} \\
    & = & F(f).\sigma_\Gamma(A).
\end{eqnarray*}
$\qed$

The identity morphism on the cwf ${\cal C}=({\cal C},\Ty,\Tm)$ is given by the 
triple $(I_{\cal C},\iota^{\cal C}, \varepsilon^{\cal C})$ where $I_{\cal C}$ is the identity functor on ${\cal C}$, $\iota^{\cal C}$ and 
$\varepsilon^{\cal C}$ are given by identities
$$\iota^{\cal C}_\Gamma = 1_{\Ty(\Gamma)}$$
and
$$\varepsilon^{\cal C}_{\Gamma,A} = 1_{\Tm(\Gamma,A)}.$$
Suppose that $$({\cal C},\Ty,\Tm) \to^{(F,\sigma, \theta)}({\cal C}',\Ty',\Tm')\to^{(F',\sigma', \theta')} ({\cal C}'',\Ty'',\Tm'')$$ are two cwf morphisms.
Define their composition 
$$(F',\sigma', \theta') \circ (F,\sigma, \theta) = (F' \circ F,\sigma' \circ \sigma,\theta' \circ \theta)$$
where
$$(\sigma' \circ \sigma)_\Gamma =_{\rm def} \sigma'_{F(\Gamma)} \circ \sigma_\Gamma$$
and 
$$(\theta' \circ \theta)_{\Gamma,A} =_{\rm def} \theta'_{F(\Gamma),\sigma_{\Gamma}(A)} \circ \theta_{\Gamma,A}.$$

\begin{theorem} \label{catofcwfs} The cwfs and cwf morphisms form a category.
\end{theorem}
{\flushleft \bf Proof.}  See Appendix A. $\qed$

\medskip 
Note that by (\ref{termpres}) it follows that empty context maps $\varepsilon$ are
preserved
\begin{equation} \label{emptympres}
 F(\varepsilon_{\Theta}) =\varepsilon_{F(\Theta)}.
\end{equation}
For morphisms $[a_1,\ldots,a_n] :\Theta \rightarrow \top.A_1.A_2.\ldots.A_n$, where
$$a_k \in \Tm(\Theta, A_k[a_1, a_2,\ldots, a_{k-1}]),$$
as in (\ref{gentupl})
we have the following useful equation 
\begin{equation} \label{tuplepres}
\begin{aligned}
& F([a_1,\ldots,a_n]_{\Theta;A_1,\ldots, A_n}) =   \\
& [\theta_{\Gamma,A_1}(a_1),\ldots,
\theta_{\Gamma, A_{n}\{[a_1, \ldots, a_{n-1} ]_{\Theta;A_1,\ldots,A_{n-1}}\}}(a_n)
]_{F(\Theta);\sigma_{\top}(A_1),\ldots,
 \sigma_{\top.A_1.\ldots.A_{n-1}}(A_n)},
\end{aligned}
\end{equation}
and dropping subscripts this reads simply
\begin{equation} \label{tuplepres2}
F([a_1,\ldots,a_n]) =  [\theta(a_1),\ldots,\theta(a_n)].
\end{equation}
Indeed, for $n=0$ the equation (\ref{tuplepres}) is just (\ref{emptympres}). The
inductive step is as follows, where (\ref{eq12}) is used in the second equation
\begin{equation}
\begin{aligned}
& F([a_1,\ldots,a_{n+1}]_{\Theta;A_1,\ldots, A_{n+1}}) \\
&=
F(\langle {[a_1, \ldots, a_n ]}_{\Theta;A_1,\ldots,A_n}, a_{n+1}\rangle_{A_{n+1}})  \\
&= \bigl\langle F({[a_1, \ldots, a_n ]}_{\Theta;A_1,\ldots,A_n}),
\theta_{\Gamma, A_{n+1}\{[a_1, \ldots, a_n ]_{\Theta;A_1,\ldots,A_n}\}}(a_{n+1}) \bigr\rangle_{\sigma_{\top.A_1.\ldots.A_n}(A_{n+1})} \\
&= \bigl\langle [\theta_{\Gamma,A_1}(a_1),\ldots,
\theta_{\Gamma, A_{n}\{[a_1, \ldots, a_{n-1} ]_{\Theta;A_1,\ldots,A_{n-1}}\}}(a_n)
]_{F(\Theta);\sigma_{\top}(A_1),\ldots,
 \sigma_{\top.A_1.\ldots.A_{n-1}}(A_n)}, \\
 & \qquad  \qquad
 \theta_{\Gamma, A_{n+1}\{[a_1, \ldots, a_n ]_{\Theta;A_1,\ldots,A_n}\}}(a_{n+1}) \bigr\rangle_{\sigma_{\top.A_1.\ldots.A_n}(A_{n+1})} \\
 &=  \bigl[\theta_{\Gamma,A_1}(a_1),\ldots,
\theta_{\Gamma, A_{n}\{[a_1, \ldots, a_{n-1} ]_{\Theta;A_1,\ldots,A_{n-1}}\}}(a_n)
, \\
 & \qquad  \qquad
 \theta_{\Gamma, A_{n+1}\{[a_1, \ldots, a_n ]_{\Theta;A_1,\ldots,A_n}\}}(a_{n+1}) 
 \bigr]_{F(\Theta);\sigma_{\top}(A_1),\ldots,
 \sigma_{\top.A_1.\ldots.A_n}(A_{n+1})}
\end{aligned}
\end{equation}

\begin{remark}  \label{FunctorRem} {\em
It follows from (\ref{cwfmorph2}), (\ref{ppres}) and (\ref{thetav}) that 
\begin{equation} \label{x2x}
\theta({\sf x}_i) = {\sf x}_i.
\end{equation}
Thus $$F(1_\Gamma) = F([{\sf x}_1,\ldots,{\sf x}_n]) = [{\sf x}_1,\ldots,{\sf x}_n]=1_\Gamma.$$
Furthermore using (\ref{substcomp}) and (\ref{cwfmorph2})
\begin{eqnarray*}
F([a_1,\ldots,a_n] \circ [b_1,\ldots,b_m]) &=&  
  [\theta (a_1[b_1,\ldots,b_m]),\ldots,\theta(a_n[b_1,\ldots,b_m])] \\
  &=& [\theta (a_1)\{F([b_1,\ldots,b_m])\}, \ldots,\theta(a_n)\{F([b_1,\ldots,b_m])\}]  \\
  &=& [\theta (a_1),\ldots, \theta(a_n)] \circ F([b_1,\ldots,b_m]) \\
  &=& F([a_1,\ldots,a_n]) \circ F([b_1,\ldots,b_m]).
\end{eqnarray*}
From this follows that if ${\cal C}$ is a contextual cwf, functoriality of $F$ is a consequence of
(\ref{tuplepres2}) and (\ref{x2x}). }
\end{remark}



\subsection{Free cwf on a signature}

\begin{construction} \label{freecwf}
The contexts and context maps for a fixed signature $\Sigma$ form a category with families ${\cal F}_\Sigma$.
\end{construction}
{\flushleft \bf Proof.} 
Let $\Sigma$ be a fixed signature. We work with derivations of judgements in ${\cal J}= {\cal J}(\Sigma)$.

Let  $\Gamma = x_1:A_1,\ldots,x_n:A_n$ be a context with respect to the signature. By rule (R3) we have for all $i=1,\ldots,n$
$$x_i:A_i\; (\Gamma)$$
Now trivially, $A_i=A_i[x_1,\ldots, x_{i-1}/x_1,\ldots, x_{i-1}]$ so
$$\iota_{\Gamma} =_{\rm def} (x_1,\ldots,x_n): \Gamma \to \Gamma$$
is a context map.

Suppose that $\Delta=y_1:B_1,\ldots,y_m:B_m$ and $\Theta = z_1:C_1,\ldots,z_k:C_k$ are contexts and that 
$$\bar{s}=(s_1,\ldots,s_m): \Gamma \to \Delta \text{ and }\bar{t}=(t_1,\ldots,t_k):\Delta \to \Theta$$ 
are context maps.
Thus 
\begin{equation} \label{map1}
\begin{array}{l}
 s_i:B_i[s_1,\ldots,s_{i-1}/y_1,\ldots,y_{i-1}] \quad (\Gamma)\qquad (i=1,\ldots,m)
\end{array}
\end{equation}
\begin{equation} \label{map2}
\begin{array}{l}
  t_j:C_j[t_1,\ldots,t_{j-1}/z_1,\ldots,z_{j-1}]  \quad (\Delta) \qquad (j=1,\ldots,k)
\end{array}
\end{equation}
By the Substitution Lemma with $\bar{s}: \Gamma \to \Delta$ applied to (\ref{map2}) we get
\begin{equation} \label{map3}
\begin{array}{l}
t_j[s_1,\ldots,s_m/y_1,\ldots,y_m]: \\
\qquad \qquad C_j[t_1,\ldots,t_{j-1}/z_1,\ldots,z_{j-1}][s_1,\ldots,s_m/y_1,\ldots,y_m]  \quad (\Gamma)\qquad (j=1,\ldots,k)
\end{array}
\end{equation}
The types of (\ref{map3}) can be rewritten as
$$\begin{array}{l}
C_j[t_1,\ldots,t_{j-1}/z_1,\ldots,z_{j-1}][s_1,\ldots,s_m/y_1,\ldots,y_m] = \\
\quad C_j[t_1[s_1,\ldots,s_m/y_1,\ldots,y_m],\ldots,t_{j-1}[s_1,\ldots,s_m/y_1,\ldots,y_m]/z_1,\ldots,z_{j-1}]
\end{array}$$
The vector defined by
 $$\bar{t} \circ \bar{s} =_{\rm def}
 \bigl(t_1[s_1,\ldots,s_m/y_1,\ldots,y_m],\ldots,t_k[s_1,\ldots,s_m/y_1,\ldots,y_m]\bigr)
 $$
is thus a context map $\Gamma \to \Theta$, the composition of $\tau$ and $\sigma$. We have
\begin{eqnarray*}
\bar{t} \circ \iota_{\Delta} &=& \bigl(t_1[y_1,\ldots,y_m/y_1,\ldots,y_m],\ldots,t_k[y_1,\ldots,y_m/y_1,\ldots,y_m]\bigr) \\
&=& (t_1,\ldots,t_k) = \bar{t}
\end{eqnarray*}
\begin{eqnarray*}
\iota_{\Delta} \circ \bar{s} &=&   \bigl(y_1[s_1,\ldots,s_m/y_1,\ldots,y_m],\ldots,y_k[s_1,\ldots,s_m/y_1,\ldots,y_m]\bigr)\ \\
&=& (s_1,\ldots,s_m) = \bar{s}
\end{eqnarray*}
Thus the identity rules are verified. Suppose now that $$\Pi=w_1:P_1,\ldots,w_{\ell}:P_{\ell}$$ is yet another context and that $\bar{p}=(p_1,\ldots,p_{\ell}): \Theta \to \Pi$ is a context map.
\begin{eqnarray*}
\bar{p} \circ (\bar{t} \circ \bar{s})&=& \bar{p} \circ\bigl(t_1[s_1,\ldots,s_m/y_1,\ldots,y_m],\ldots,t_k[s_1,\ldots,s_m/y_1,\ldots,y_m]\bigr)
\\
&=& \bigl(p_1[t_1[s_1,\ldots,s_m/y_1,\ldots,y_m],\ldots,t_k[s_1,\ldots,s_m/y_1,\ldots,y_m]/z_1,\ldots.z_k],  \\
&&    \ldots, p_{\ell}[t_1[s_1,\ldots,s_m/y_1,\ldots,y_m],\ldots,t_k[s_1,\ldots,s_m/y_1,\ldots,y_m]/z_1,\ldots.z_k]\bigr)\\
&=& \bigl(p_1[t_1,\ldots,t_k/z_1,\ldots.z_k][s_1,\ldots,s_m/y_1,\ldots,y_m], \\
&&    \ldots, p_{\ell}[t_1,\ldots,t_k/z_1,\ldots.z_k][s_1,\ldots,s_m/y_1,\ldots,y_m]\bigr)\\
&=&  \bigl(p_1[t_1,\ldots,t_k/z_1,\ldots.z_k], \ldots, p_{\ell}[t_1,\ldots,t_k/z_1,\ldots.z_k]\bigr) \circ \bar{s} \\
&=& (\bar{p} \circ \bar{t}) \circ \bar{s}
\end{eqnarray*}
From this we can build a category ${\cal F}_{\Sigma} = {\cal F}$ consisting of contexts and context maps. The objects of
the category is  formally then
\begin{equation} \label{x1}
{\rm Ob} \, {\cal F} = 
\{ \Gamma \in {\rm precontext} : (\Gamma \; {\rm context}) \in {\cal J}\}.
\end{equation}
and arrows
$${\rm Arr} \, {\cal F} = 
\{(\Gamma,\Delta,\bar{a}) : (\bar{a}: \Gamma \to \Delta) \in {\cal J}\}.$$
Composition of arrows is defined as
$$(\Gamma,\Delta,\bar{a}) \circ (\Theta,\Gamma,\bar{b}) = (\Theta,\Delta,\bar{a} \circ \bar{b}).$$
The identity arrow is given by
$$1_\Gamma = (\Gamma, \Gamma, \iota_\Gamma).$$
The empty context $\langle \rangle$ is the terminal object in ${\cal F}$.
For objects $\Gamma, \Delta$ we write as usual ${\cal F}(\Gamma,\Delta)$ for
the arrows $\sigma =(\Gamma', \Delta',\bar{s})$ with $\Gamma' = \Gamma$ and
$\Delta' = \Delta$ (syntactic equalities). We introduce further notation for substitution,
writing 
$$A[\sigma/\Delta] = A[\bar{s}/\Delta] \mbox{ and } a[\sigma/\Delta] = a[\bar{s}/\Delta].$$

We define a functor ${\rm Ty}:{\cal F}^{\rm op} \to {\rm Set}$ which assigns
types to contexts
\begin{equation} \label{x2}
{\rm Ty}(\Gamma) = \{(\Gamma, A) \in {\rm precontext} \times {\rm pretype} :  
 (A \; {\rm type} \; (\Gamma)) \in {\cal J}\},
\end{equation}
and for $\sigma \in {\cal F}(\Delta,\Gamma)$, and ${\rm A} = (\Gamma, A)$  let 
$${\rm Ty}(\sigma)({\rm A}) = (\Delta, A[\sigma/\Gamma]).$$
Note that ${\rm Ty}(\sigma)({\rm A}) \in {\rm Ty}(\Delta)$ by the Substitution Lemma.
Clearly
$${\rm Ty}(1_{\Gamma})({\rm A}) = (\Gamma, A[\iota_{\Gamma}/\Gamma])={\rm A}.$$
For $\sigma = (\Delta, \Gamma, \bar{s}) \in {\cal F}(\Delta,\Gamma)$ and 
 $\tau = (\Theta, \Delta, \bar{t}) \in {\cal F}(\Theta,\Delta)$, where $\Delta=\langle y_1:B_1,\ldots,y_m:B_m \rangle$, we have for ${\rm A} = (\Theta,A)$
\begin{eqnarray*}
{\rm Ty}(\sigma \circ  \theta)({\rm A}) &=& (\Theta, A[\bar{s} \circ \bar{t}/\Gamma]) \\
&=& (\Theta, A[s_1[t_1,\ldots,t_m/y_1,\ldots,y_m],\ldots,  \\
& & \qquad \qquad s_n[t_1,\ldots,t_m/y_1,\ldots,y_m]/x_1,\ldots,x_n]) \\
&=& (\Theta, A[s_1,\ldots, s_n/x_1,\ldots,x_n] [t_1,\ldots,t_m/y_1,\ldots,y_m])\\
&=& (\Theta, A[\bar{s}/\Gamma] [\bar{t}/\Delta])\\
&=&{\rm Ty}(\theta)({\rm Ty}(\sigma)({\rm A})).\\
\end{eqnarray*}
We write as is usual ${\rm A}\{\sigma\}$ for ${\rm Ty}(\sigma)({\rm A})$.

For $\Gamma \in {\cal F}$ and ${\rm A} = (\Gamma, A) \in {\rm Ty}(\Gamma)$, let  the set of terms be
\begin{equation} \label{x3}
\begin{aligned}
{\rm Tm}(\Gamma,{\rm A}) =\{((\Gamma, A), a) & \in ({\rm precontext}\times 
{\rm pretype}) \times {\rm preelement} \,|\,  \\
&\qquad \qquad \qquad(a:A\; (\Gamma)) 
\in {\cal J}\}.
\end{aligned}
\end{equation}

Define for $\sigma \in {\cal F}(\Delta, \Gamma)$ and ${\rm a} = ({\rm A}, a) \in {\rm Tm}(\Gamma,{\rm A})$, 
${\rm A} = (\Gamma, A)$
$${\rm a}\{\sigma\} = ({\rm A}\{\sigma\}, a[\sigma/\Gamma])  
= ((\Delta, A[\sigma/\Gamma]), a[\sigma/\Gamma]).$$ 
It is easily checked that ${\rm a}\{\sigma\} \in {\rm Tm}(\Delta,A\{\sigma\})$, by using the Substitution Lemma, and that 
$${\rm a}\{1_\Gamma\} = {\rm a}$$ 
$${\rm a}\{\sigma \circ \tau\} = {\rm a}\{\sigma\}\{\tau\} \qquad 
(\tau \in {\cal F}(\Theta,\Delta))$$

If $\Gamma = x_1:A_1,\ldots,x_n:A_n \in {\rm Ob}({\cal F})$ and ${\rm S} = (\Gamma,S)
\in \Ty(\Gamma)$, define employing the fresh variable provider
$$\Gamma.{\rm S} =_{\rm def}
 \langle \Gamma, {\rm fresh}(\Gamma):S\rangle$$
and
$${\rm p}_{\Gamma}({\rm S}) =_{\rm def} (\Gamma. {\rm S}, \Gamma, (x_1,\ldots,x_n)).$$
Note the special cases of substitution: for ${\rm p}_{\Gamma}({\rm S}) \in {\cal F}(\Gamma.{\rm S}, \Gamma)$
and  ${\rm A} = (\Gamma, A) \in {\rm Ty}(\Gamma)$,
$${\rm A}\{{\rm p}_{\Gamma}({\rm S})\} = (\Gamma.{\rm S}, A[{\rm p}_{\Gamma}({\rm S})/\Gamma])
= (\Gamma.{\rm S}, A[x_1,\ldots,x_n/x_1,\ldots,x_n]) = (\Gamma.{\rm S}, A).$$
and ${\rm a} = ({\rm A}, a)  \in {\rm Tm}(\Gamma,{\rm A})$,
\begin{equation} \label{eq4_2}
{\rm  a}\{{\rm p}_{\Gamma}({\rm S})\} =  ((\Gamma.{\rm S}, A[{\rm p}_{\Gamma}({\rm S})/\Gamma]), a[{\rm p}_{\Gamma}({\rm S})/\Gamma])
=  ((\Gamma.{\rm S}, A), a)
\end{equation}
We then introduce 
\begin{equation} \label{eq4_3}
\begin{aligned}
{\rm v}_{\Gamma,{\rm S}} = _{\rm def}   ({\rm S}\{{\rm p}_{\Gamma}({\rm S})\},  {\sf fr}(x_1,\ldots,x_n)) &= 
((\Gamma.{\rm S},S[{\rm p}_{\Gamma}(S)/\Gamma]), {\sf fr}(x_1,\ldots,x_n)) \\
&= ((\Gamma.{\rm S},S), {\sf fr}(x_1,\ldots,x_n))
\end{aligned}
\end{equation}
 which is readily proved to be in
$${\rm Tm}(\Gamma.{\rm S},{\rm S}\{{\rm p}_{\Gamma}({\rm S})\}).$$

Then for  $\sigma = (\Delta, \Gamma, (s_1,\ldots,s_n)) \in {\cal F}(\Delta,\Gamma)$,
${\rm S}= (\Gamma, S) \in {\rm Ty}(\Gamma)$, and ${\rm a} =({\rm S},a) \in {\rm Tm}(\Delta, {\rm S}\{\sigma\})$, define
$$\langle \sigma, {\rm a} \rangle_{\rm S} = (\Delta, \Gamma.{\rm S}, (s_1,\ldots,s_n,a)).$$ 
which is a context morphism ${\cal F}(\Delta, \Gamma.{\rm S})$. Then 
for $\tau \in {\cal F}(\Theta, \Delta)$, 
we have
$$\langle \sigma,{\rm a}\rangle_{\rm S} \circ \tau = \langle \sigma \circ \tau, {\rm a}\{\tau\} \rangle_{\rm S}.$$
Moreover
$${\rm p}_{\Gamma}({\rm S}) \circ \langle \sigma, {\rm a}\rangle_{\rm S} = \sigma$$
and
$${\rm v}_{\Gamma,{\rm S}} \{\langle \sigma,{\rm a}\rangle_{\rm S}\} = {\rm a}$$
and finally
$$1_{\Gamma,{\rm S}} = \langle {\rm p}_{\Gamma}({\rm S}), {\rm v}_{\Gamma,{\rm S}}  \rangle_{\rm S}.$$ 
$\qed$

\medskip
We now show that judgements in ${\cal J}(\Sigma)$  can be derived in contexts with a standardized order of variables, and vice versa  (Corollary \ref{judgementstandard}).

\begin{definition} {\em Let $\Sigma = ({\cal V}, F, T)$ where the variable system is ${\cal V} =(V, \varphi, {\sf fr})$. A context 
$$x_1:A_1,\ldots,x_n:A_n$$
is on {\em standard form} with respect to $\Sigma$, if $x_k= {\sf fr}(\{x_1,\ldots, x_{k-1}\})$
for all $k=1,\ldots,n$.
An infinite sequence $\sigma:{\mathbb N} \to V$ of variables is {\em fresh} with respect to $\Sigma$ if
$$\sigma(n) \in \varphi(\{\sigma(0),\ldots,\sigma(n-1)\}),$$
for all $n \in {\mathbb N}$.  The standard ordering of variables given
by $v_n= {\sf fr}(\{v_0,\ldots, v_{n-1}\})$ is an example of a fresh infinite sequence.
}
\end{definition}

We define a context map which puts variables in the order of a fresh sequence $\sigma$. For a precontext $\Gamma$, let $\sigma(\Gamma):\Gamma^\sigma \to \Gamma$ be given by
the simultaneous recursive definition
\begin{eqnarray*}
\sigma(\langle \rangle) &=& () \\
\sigma(\langle \Gamma, x:A \rangle) &=& (\sigma(\Gamma), \sigma(|\Gamma|))\\
\langle \rangle^\sigma &=& \langle \rangle \\
\langle \Gamma, x:A \rangle^\sigma &=&  \langle \Gamma^\sigma, 
\sigma(|\Gamma|):A[\sigma(\Gamma)/\Gamma]\rangle.
\end{eqnarray*}
This is a context map whenever $\Gamma$ is a context:

\begin{lemma} \label{varstandard} Suppose that $\Sigma$ is a signature, and that
$\sigma: {\mathbb N} \to V$ is a fresh sequence of variables for $\Sigma$. If $(\Gamma\; {\rm context}) \in {\cal J}(\Sigma)$, then variable ordering $\sigma$ for $\Gamma$ forms
a context map
$$(\sigma(\Gamma): \Gamma^\sigma \to \Gamma) \in {\cal J}(\Sigma),$$
whose inverse 
$(\sigma^{-1}(\Gamma):\Gamma \to \Gamma^\sigma) \in {\cal J}(\Sigma)$, 
is given by $\sigma^{-1}(\Gamma) = {\rm OV}(\Gamma)$.
\end{lemma}
{\flushleft \bf Proof.} See Appendix A. $\qed$

\medskip
We can now prove the result about variable standardization for judgements, as
a corollary to Lemma \ref{varstandard}.

\begin{corollary}  \label{judgementstandard} Suppose that $\Sigma$ is a signature,
and that $\sigma: {\mathbb N} \to V$ is a fresh sequence of variables for $\Sigma$.
Then
\begin{itemize}
\item[(a)] if $(\Gamma \; {\rm context}) \in {\cal J}(\Sigma)$, then
$(\Gamma^{\sigma} \; {\rm context}) \in {\cal J}(\Sigma)$,
\item[(b)] if $(\Gamma \; {\rm context}) \in {\cal J}(\Sigma)$, then  
$$(A \; {\rm type} \; (\Gamma)) \in {\cal J}(\Sigma) \Longleftrightarrow (A[\sigma(\Gamma)/\Gamma] \; {\rm type} \; (\Gamma^{\sigma}) \in {\cal J}(\Sigma),$$
\item[(c)] if $(\Gamma \; {\rm context}) \in {\cal J}(\Sigma)$, then 
$$(a: A  \; (\Gamma)) \in {\cal J}(\Sigma)   \Longleftrightarrow (a[\sigma(\Gamma)/\Gamma]: A[\sigma(\Gamma)/\Gamma] \; (\Gamma^{\sigma}) \in {\cal J}(\Sigma),$$
\item[(d)] if $\Gamma$ is a precontext with respect to $\Sigma$, and 
$(\Gamma^{\sigma} \; {\rm context}) \in {\cal J}(\Sigma)$,
then $(\Gamma \; {\rm context}) \in {\cal J}(\Sigma)$.
\end{itemize}
\end{corollary}
{\flushleft \bf Proof.} As for (a): If $(\Gamma \; {\rm context}) \in {\cal J}(\Sigma)$, then by Lemma \ref{varstandard},
\begin{equation} \label{js1}
(\sigma(\Gamma):\Gamma^{\sigma} \to \Gamma) \in {\cal J}(\Sigma).
\end{equation}
In particular $(\Gamma^{\sigma} \; {\rm context}) \in {\cal J}(\Sigma)$.

As for (b): By Lemma \ref{varstandard} and  the Substitution lemma we obtain
$$(A[\sigma(\Gamma)/\Gamma] \; {\rm type} \; (\Gamma^{\sigma}) \in {\cal J}(\Sigma),$$
as required. The converse is obtained by applying the inverse substitution
$\sigma^{-1}(\Gamma):\Gamma \to \Gamma^{\sigma}$.

As for (c): By Lemma \ref{varstandard} and the Substitution lemma  again we get
$$(a[\sigma(\Gamma)/\Gamma]: A[\sigma(\Gamma)/\Gamma] \; (\Gamma^{\sigma}) \in {\cal J}(\Sigma).$$
The converse is obtained by applying the inverse substitution. 

As for (d): Induction on derivations. For $\Gamma = \langle \rangle$ the statement is trivial.
Suppose $\Gamma,x:A$ is a precontext. Hence $x \in {\rm Fresh}(\Gamma)$ and
$\Gamma$ is a precontext.
If $((\Gamma,x:A)^\sigma \; {\rm context})\in {\cal J}(\Sigma)$, then with shorter
derivations  $(\Gamma^{\sigma} \; {\rm context})\in {\cal J}(\Sigma)$ and
$(A[\sigma(\Gamma)/\Gamma] \; {\rm type}\; (\Gamma^{\sigma})) \in {\cal J}(\Sigma)$.
By inductive hypothesis, $(\Gamma \; {\rm context}) \in {\cal J}(\Sigma)$.
We can apply Lemma \ref{varstandard}, and by (b) obtain,
$(A\; {\rm type}\; (\Gamma)) \in {\cal J}(\Sigma)$. Now as $x \in {\rm Fresh}(\Gamma)$,
we get $(\Gamma, x:A \; {\rm context}) \in {\cal J}(\Sigma)$ as required.
$\qed$

\medskip
\begin{remark}  \label{freerem}
{\em
The cwf ${\cal F}_{\Sigma}$  defined in Construction \ref{freecwf} 
is in general not contextual, unless we restrict
the fresh variable provider. Consider a signature with a constant type $T$. Then $\langle y:T \rangle$ is a context for any variable 
$y$ when $\varphi=\varphi_{\infty}$, but unless $y={\sf fr}_{\infty}()$, this context can never be written as an extension of the terminal object in ${\cal F}_{\Sigma}$.

To obtain a contextual cwf we consider the fresh variable provider $\varphi_1$, ${\sf fr}_1$ and 
$V={\mathbb N}_{+}$. In effect this just calling the $n$th variable in the context, $n$.
Indeed note that then any context $\Gamma$ in ${\cal F}_{\Sigma}$ has the form
$$1:A_1,\ldots,n:A_n.$$
The $i$th projection is, by (\ref{eq4_2}) and (\ref{eq4_3})
\begin{equation} \label{eq4_4}
{\sf x}_i = {\rm v}\{{\rm p}^{(n-i)}\} = ((\Gamma, A_i),i).
\end{equation}
}
\end{remark}

By this remark we can easily see that

\begin{theorem} \label{freecwf2}
The contexts and context maps for a fixed signature $\Sigma$, where the fvp is of de Bruijn type, form a contextual cwf ${\cal F}_\Sigma$. $\qed$
\end{theorem}

\medskip
\subsection{Grading the free cwf by proof height}

\medskip
The components of ${\cal F}_\Sigma$ have a natural grading based on
height of derivations.  We introduce graded versions of 
(\ref{x1}), (\ref{x2}) and (\ref{x3}) using (\ref{gradedJ}) as follows. For $n \in {\mathbb N}$,
\begin{itemize}
\item $
{\rm Ob}^\Sigma_n= {\rm Ob}_n \, {\cal F}_{\Sigma} = 
\{ \Gamma \in {\rm precontext} \,|\, (\Gamma \; {\rm context}) \in {\cal J}_n(\Sigma)\}.$
\item For $\Gamma \in {\rm Ob}^\Sigma_n$,
$${\rm Ty}^\Sigma_n(\Gamma) = \{(\Gamma, A) \in {\rm precontext} \times {\rm pretype} \,|\,  (A \; {\rm type} \; (\Gamma)) \in {\cal J}_n(\Sigma)\},$$
\item For $\Gamma \in {\rm Ob}^\Sigma_n$, $(\Gamma,A) \in {\rm Ty}^{\Sigma}_n(\Gamma)$,
$$
\begin{aligned}
{\rm Tm}^\Sigma_n(\Gamma,{\rm A}) =\{((\Gamma, A), a) \in & ({\rm precontext}\times 
{\rm pretype}) \times {\rm preelement} \,|\,   \\
& \qquad (a:A\; (\Gamma)) 
\in {\cal J}_n(\Sigma)\}.
\end{aligned}
$$
\end{itemize}
Each of these sets is monotone in $n$, i.e. for $m<n$ we have
\begin{itemize}
\item ${\rm Ob}^\Sigma_m \subseteq {\rm Ob}^\Sigma_n$,
\item ${\rm Ty}^\Sigma_m(\Gamma) \subseteq {\rm Ty}^\Sigma_n(\Gamma)$ for $\Gamma \in
{\rm Ob}^\Sigma_m$,
\item ${\rm Tm}^\Sigma_m(\Gamma,{\rm A})  \subseteq {\rm Tm}^\Sigma_n(\Gamma,{\rm A})$
for $\Gamma \in {\rm Ob}^\Sigma_m$, and $(\Gamma,{\rm A}) \in 
{\rm Ty}^\Sigma_m(\Gamma)$.
\end{itemize}
We note that these sets have the following presupposition like properties:
\begin{itemize}
\item[(P1)] If $\Gamma \in {\rm Ob}^\Sigma_n$, then either $\Gamma = \langle\rangle$, or
$\Gamma= \langle \Gamma', u:C\rangle$ with $n>0$ and $\Gamma' \in  {\rm Ob}^\Sigma_{n-1}$ and
$(\Gamma',C) \in {\rm Ty}^{\Sigma}_{n-1}(\Gamma')$.
\item[(P2)] If $(\Gamma,A) \in {\rm Ty}^\Sigma_n(\Gamma)$ where $\Gamma \in {\rm Ob}^\Sigma_n$, then $n>0$ and there is a unique declaration $(\Delta,\bar{i},S) \in \Sigma$ and
a unique $(\bar{t}: \Gamma \to \Delta) \in {\cal J}_{n-1}(\Sigma)$, so that $A=S(\bar{t})$.
Thus $\Gamma  \in {\rm Ob}^\Sigma_{n-1}$, $\Delta  \in {\rm Ob}^\Sigma_{n-1}$ and if $\Delta = y_1:B_1,\ldots, y_r:B_r$,  then for $k=1,\ldots,r$:
$$t_k: B_k[\bar{t}^{k-1}/\Delta^{k-1}]\; (\Gamma) \in {\cal J}_{n-1}(\Sigma).$$
Then by Theorem \ref{presup}, $
B_k[\bar{t}^{k-1}/\Delta^{k-1}]\; {\rm type}\; (\Gamma) \in {\cal J}_{n-1}(\Sigma)$.
\item[(P3)] If $((\Gamma,A),a) \in {\rm Tm}^\Sigma_n(\Gamma,{\rm A})$ where  
$(\Gamma,A) \in {\rm Ty}^\Sigma_n(\Gamma)$ and $\Gamma \in {\rm Ob}^\Sigma_n$, then
$n>0$ and there is a unique declaration $(\Delta,f,\bar{i},U) \in \Sigma$ and
unique $(\bar{t}: \Gamma \to \Delta) \in {\cal J}_{n-1}(\Sigma)$, so that
$a=f(\bar{t})$ and $A=U[\bar{t}/\Delta]$. Thus $\Gamma  \in {\rm Ob}^\Sigma_{n-1}$, $\Delta  \in {\rm Ob}^\Sigma_{n-1}$ and if $\Delta = y_1:B_1,\ldots, y_r:B_r$,  then for $k=1,\ldots,r$:
$$t_k: B_k[\bar{t}^{k-1}/\Delta^{k-1}]\; (\Gamma) \in {\cal J}_{n-1}(\Sigma).$$
Again by Theorem \ref{presup}, 
$B_k[\bar{t}^{k-1}/\Delta^{k-1}]\; {\rm type}\; (\Gamma) \in {\cal J}_{n-1}(\Sigma)$.
\end{itemize}

\section{Models for signatures and extension problems} \label{extprob}

We use a functorial definition of a model for a signature, and then show 
how it relates to models defined by specifying the component of the signature.

In this section {\em we consider only regular signatures $\Sigma$, over symbol systems where the fresh variable provider is $(\varphi_1,{\sf fr}_1)$.} Thus according to Theorem \ref{freecwf2}, ${\cal F}_{\Sigma}$ will be a contextual cwf.

\subsection{Models for type systems}

A {\em model of a type system} over the signature $\Sigma$ is a cwf morphism 
$(F,\sigma,\theta):{\cal F}_{\Sigma} \rightarrow {\cal C}$ into some cwf ${\cal C}$.
Such models are determined by their values on the component of the signature
(Lemma \ref{cwfuniq}). 
To see this we first use Definition \ref{defcwf} to spell out the meaning of 
$(F,\sigma,\theta):{\cal F}_{\Sigma} \rightarrow {\cal C}$
being a cwf morphism, by considering judgements in ${\cal J} = {\cal J}(\Sigma)$:
\begin{itemize}
\item[(M1)] If $(\Gamma \; {\rm context})  \in {\cal J}$, then  $F(\Gamma) \in {\rm Ob}\; {\cal C}$.
\item[(M2)] If $( A \; {\rm type}\;   (\Gamma) )  \in {\cal J}$, then $\sigma_{\Gamma}(\Gamma, A) \in \Ty_{{\cal C}}(F(\Gamma))$.
\item[(M3)] If $(a:A \;  (\Gamma) )   \in {\cal J}$, then $\theta_{\Gamma,A}((\Gamma,A),a) \in \Tm_{{\cal C}}(F(\Gamma),\sigma_{\Gamma}(\Gamma, A)).$
\end{itemize}
We drop the subscripts from $\sigma$ and $\theta$ below.

For $(\bar{b}: \Delta \to \Gamma) \in {\cal J}$, we have 
\begin{itemize}
\item[(M4)] if $( A \; {\rm type}\;   (\Gamma) )  \in {\cal J}$,  then
$\sigma(\Delta, A[\bar{b}/\Gamma]) = \sigma(\Gamma, A)\{F(\Delta,\Gamma,\bar{b})\},$
\item[(M5)] if $(a:A \;  (\Gamma) )   \in {\cal J}$, then
$\theta(\Delta, A[\bar{b}/\Gamma], a[\bar{b}/\Gamma]) = 
\theta(\Gamma, A, a)\{F(\Delta,\Gamma,\bar{b})\},$
\item[(M6)] if $\bar{b}=b_1,\ldots,b_n$ and $\Gamma = x_1:B_1,\ldots,x_n:B_n$ then
$$F(\Delta,\Gamma,\bar{b}) = 
[\theta(\Delta, B_1, b_1),\ldots,\theta(\Delta, B_n[\bar{b}^{n-1}/\Gamma^{n-1}], b_n)].$$
\end{itemize}
We introduce the notation
$$\bar{\theta}(\Delta,\Gamma,\bar{b}) =_{\rm def}
[\theta(\Delta, B_1, b_1),\ldots,\theta(\Delta, B_n[\bar{b}^{n-1}/\Gamma^{n-1}], b_n)].$$
Moreover we have
\begin{itemize}
\item[(M7)] if $(x_i:A_i \;  (x_1:A_1,\ldots,x_n:A_n))   \in {\cal J}$,
$$\theta(\langle x_1:A_1,\ldots,x_n:A_n\rangle,A_i, x_i) = {\sf x}_i.$$
\end{itemize}
As particular cases of (M4) and (M5) we have, letting $\bar{x}= {\rm OV}(\Gamma)$,
\begin{itemize}
\item[(M8)] for $(\Gamma,S) \in \Sigma$, 
$\sigma(\Delta,S(\bar{b})) = \sigma(\Gamma, S(\bar{x}))\{F(\Delta,\Gamma,\bar{b})\},$
\item[(M9)] for $(\Gamma,f,U) \in \Sigma$, 
$\theta(\Delta,U[\bar{b}/\Gamma], f(\bar{b})) = 
\theta(\Gamma,U, f(\bar{x}))\{F(\Delta,\Gamma,\bar{b})\}.$
\end{itemize}

\begin{lemma} \label{cwfuniq} Let $\Sigma$ be a regular signature. Suppose that $(F,\sigma,\theta), (F',\sigma',\theta'):{\cal F}_{\Sigma} \rightarrow {\cal C}$ are two cwf morphisms
with 
\begin{eqnarray*}
&F(\Gamma) = F'(\Gamma) \qquad 
&(\text{$(\Gamma,S) \in \Sigma$ or  $(\Gamma,f,U) \in \Sigma$}) \\
&\sigma(\Gamma, S(\check{\Gamma})) = \sigma'(\Gamma, S(\check{\Gamma}))
\qquad &((\Gamma,S) \in \Sigma) \\
&\sigma(\Gamma, U) = \sigma'(\Gamma, U)
\qquad &((\Gamma,f,U) \in \Sigma) \\
& \theta(\Gamma,U,f(\check{\Gamma})) = \theta'(\Gamma,U,f(\check{\Gamma})) \qquad &((\Gamma,f,U) \in \Sigma) 
\end{eqnarray*}
Then $(F,\sigma,\theta) = (F',\sigma',\theta')$.
\end{lemma}
{\flushleft \bf Proof.} Induction on the grading $n$ of  ${\rm Ob}_n^\Sigma$, 
$\Ty_n^\Sigma$ and $\Tm_n^\Sigma$. We prove that for all $n$:
\begin{itemize}
\item[(a)] $F(\Gamma) = F'(\Gamma)$ for $\Gamma \in {\rm Ob}_n^\Sigma$ 
\item[(b)] $\sigma(\Gamma, A) = \sigma'(\Gamma, A)$ for 
$\Gamma \in {\rm Ob}_n^\Sigma$ and $(\Gamma, A) \in 
\Ty_n^\Sigma(\Gamma)$
\item[(c)] $\theta(\Gamma, A, a) = \theta'(\Gamma, A, a)$ for
$\Gamma \in {\rm Ob}_m^\Sigma$, $(\Gamma, A) \in 
\Ty_n^\Sigma(\Gamma)$, and $((\Gamma,A),a) \in
\Tm_n^\Sigma(\Gamma,A)$.
\end{itemize}

$n=0$: We have $F(\langle\rangle) = \top_{\cal C} = F'(\langle \rangle)$ so
since ${\rm Ob}_0^\Sigma = \{\langle \rangle\}$,  $F$ and $F'$ are equal on 
this set. The sets $\Ty_0^\Sigma$ and $\Tm_0^\Sigma$ are empty
so trivially $\sigma$ and $\sigma'$, and $\theta$ and $\theta'$ agree on 
these sets respectively. Thus (a) -- (c) holds for $n=0$.

\medskip
$n=m+1$:  Assume as inductive hypothesis that (a) -- (c) holds for $n=m$.
Suppose that $\Gamma \in {\rm Ob}_{m+1}^\Sigma$. 
Thus either $\Gamma=\langle\rangle$ or  
$\Gamma=\langle\Gamma', x: A\rangle$.
In the first case, $\Gamma \in {\rm Ob}_0^\Sigma$, and 
$F(\Gamma) = F'(\Gamma)$ by the base case.
In the second case we have according to  (P1) above that $\Gamma' \in 
{\rm Ob}_m^\Sigma$ and $(\Gamma',C) \in \Ty_m^\Sigma(\Gamma')$.
Thus by the inductive hypothesis
$$F(\langle\Gamma', x: A\rangle) = 
F(\Gamma).\sigma(\Gamma',A) =
F'(\Gamma).\sigma'(\Gamma',A) =
F'(\langle\Gamma', x: A\rangle).
$$
Thus (a) is proved for $n=m+1$.

Suppose that
$\Gamma \in {\rm Ob}_{m+1}^\Sigma$ and $(\Gamma, A) \in 
\Ty_{m+1}^\Sigma(\Gamma)$. By (P2) above there is a 
unique declaration $(\Delta,\bar{i},S) \in \Sigma$ and
a unique $(\bar{t}: \Gamma \to \Delta) \in {\cal J}_{m}(\Sigma)$, so that $A=S(\bar{t})$.
Thus $\Gamma  \in {\rm Ob}^\Sigma_{m}$, $\Delta  \in {\rm Ob}^\Sigma_{m}$ and writing $\Delta = y_1:B_1,\ldots, y_r:B_r$,  we have for $k=1,\ldots,r$:
$$t_k: B_k[\bar{t}^{k-1}/\Delta^{k-1}]\; (\Gamma) \in {\cal J}_{m}(\Sigma)$$
and $$
B_k[\bar{t}^{k-1}/\Delta^{k-1}]\; {\rm type}\; (\Gamma) \in {\cal J}_{m}(\Sigma).$$
By the induction hypothesis (b) for $m$ we have
\begin{equation} \label{sigmaih}
\sigma(\Gamma,B_k[\bar{t}^{k-1}/\Delta^{k-1}]) =
\sigma'(\Gamma,B_k[\bar{t}^{k-1}/\Delta^{k-1}]) \quad (k=1,\ldots,n).
\end{equation}
By inductive hypothesis (c) for $m$ we then have
$$
[\theta(\Gamma, B_1, t_1),\ldots,\theta(\Gamma, B_n[\bar{t}^{n-1}/\Delta^{n-1}], t_n)] =
[\theta'(\Gamma, B_1, t_1),\ldots,\theta'(\Gamma, B_n[\bar{t}^{n-1}/\Delta^{n-1}], t_n)]
$$
so
\begin{equation} \label{thetaih}
\bar{\theta}(\Gamma,\Delta,\bar{t}) = \bar{\theta'}(\Gamma,\Delta,\bar{t}).
\end{equation}
Now
\begin{eqnarray*}
\sigma(\Gamma,A) & = & \sigma(\Delta,S(\check{\Delta})\{\bar{\theta}(\Gamma,\Delta,\bar{t})\} \qquad (\text{(M6) and (M8)})\\
  & = & \sigma'(\Delta,S(\check{\Delta}))\{\bar{\theta}(\Gamma,\Delta,\bar{t})\}  \qquad (\text{by (ii)})
  \\
  & = & \sigma'(\Delta,S(\check{\Delta}))\{\bar{\theta'}(\Gamma,\Delta,\bar{t})\} 
   \qquad (\ref{thetaih})\\  
  & = & \sigma'(\Gamma,A). \qquad \qquad \qquad \qquad \qquad \quad (\text{(M6) and (M8)})
\end{eqnarray*}
Thus (b) is proved for $n=m+1$.

Suppose that
$\Gamma \in {\rm Ob}_{m+1}^\Sigma$, $(\Gamma, A) \in 
\Ty_{m+1}^\Sigma(\Gamma)$ and $((\Gamma,A),a) \in
\Tm_{m+1}^\Sigma(\Gamma,A)$. By (P3) there is a unique declaration $(\Delta,f,\bar{i},U) \in \Sigma$ and
unique $(\bar{t}: \Gamma \to \Delta) \in {\cal J}_{n-1}(\Sigma)$, so that
$a=f(\bar{t})$ and $A=U[\bar{t}/\Delta]$. Thus $\Gamma  \in {\rm Ob}^\Sigma_{m}$, $\Delta  \in {\rm Ob}^\Sigma_{m}$ and writing $\Delta = y_1:B_1,\ldots, y_r:B_r$,  we have for $k=1,\ldots,r$:
$$t_k: B_k[\bar{t}^{k-1}/\Delta^{k-1}]\; (\Gamma) \in {\cal J}_{m}(\Sigma)$$
and $$
B_k[\bar{t}^{k-1}/\Delta^{k-1}]\; {\rm type}\; (\Gamma) \in {\cal J}_{m}(\Sigma).$$
As above the inductive hypothesis (b) and (c) for $m$ gives  (\ref{sigmaih}) and (\ref{thetaih}).
\begin{eqnarray*}
\theta(\Gamma, A, a) &=&  \theta(\Delta,U, f(\bar{x}))\{\bar{\theta}(\Gamma,\Delta,\bar{t})\} \\
&=& \theta'(\Delta,U, f(\bar{x}))\{\bar{\theta}(\Gamma,\Delta,\bar{t})\}  \qquad (\text{by (iv)}) \\
&=& \theta'(\Delta,U, f(\bar{x}))\{\bar{\theta'}(\Gamma,\Delta,\bar{t})\} \qquad (\ref{thetaih}) \\
&=&  \theta'(\Gamma, A, a). \\
\end{eqnarray*}
This proves (c) for (m+1). $\qed$

\medskip
The following result show that it is enough to specify the values on objects of
 ${\cal F}_{\Sigma}$ to define a cwf morphism ${\cal F}_{\Sigma} \to {\cal C}$.

\begin{theorem}  \label{thm41}
Suppose that $F$, $\sigma$ 
and $\theta$ are assignments  that satisfy (M1) -- (M3), (M7)  and that
for $(\bar{b}: \Delta \to \Gamma) \in {\cal J}(\Sigma)$,  
\begin{itemize}
\item[(M8')] for all $(\Gamma,S) \in \Sigma$, 
$\sigma(\Delta, S(\bar{b})) = 
\sigma(\Gamma, S(\bar{x}))\{
\bar{\theta}(\Delta,\Gamma, \bar{b})\},$
\item[(M9')] for all $(\Gamma,f,U) \in \Sigma$, 
$
\theta(\Delta,U[\bar{b}/\Gamma], f(\bar{b})) =  
\theta(\Gamma,U,f(\bar{x}))\{\bar{\theta}(\Delta,\Gamma, \bar{b})\}.
$
\end{itemize}
Then letting 
\begin{equation} \label{Fdef}
F(\Delta,\Gamma,\bar{b}) =\bar{\theta}(\Delta,\Gamma, \bar{b}),
\end{equation} 
we obtain a cwf morphism $(F,\sigma,\theta): {\cal F}_{\Sigma} \rightarrow {\cal C}$. 
\end{theorem}
{\bf \flushleft Proof.} By induction on the height of derivations one shows that (M8') and (M9') implies (M4) and (M5).  Suppose that for all $m<k$, we have for any $(\bar{b}: \Delta \to \Gamma) \in {\cal J}(\Sigma)$, 
$\Gamma = x_1:B_1,\ldots,x_n:B_n$
\begin{itemize}
\item[(M4')] if $\vdash_m A \; {\rm type}\;   (\Gamma)$, 
$$\sigma(\Delta, A[\bar{b}/\Gamma]) = \sigma(\Gamma, A)\{[\theta(\Delta,B_1,b_1),\ldots,\theta(\Delta,B_n[\bar{b}^{n-1}/\Gamma^{n-1}], b_n)]\},$$
\item[(M5')] if $\vdash_m a:A \;  (\Gamma)$,
$$\theta(\Delta, A[\bar{b}/\Gamma], a[\bar{b}/\Gamma]) = 
\theta(\Gamma, A, a)\{[\theta(\Delta,B_1,b_1),\ldots,\theta(\Delta,B_n[\bar{b}^{n-1}/\Gamma^{n-1}], b_n)]\},$$
\end{itemize}

Ad (M4'): 
Suppose $\vdash_k A \; {\rm type}\;  (\Gamma)$. Thus $A = S(a_1,\ldots,a_p)$, for some $(\Theta,S) \in \Sigma$ with
\begin{equation} \label{avect}
\vdash_{\bar{r}}(a_1,\ldots,a_p): \Gamma \to \Theta
\end{equation}
and $\bar{r} < k$. Assume $\Theta = z_1:C_1,\ldots, z_p:C_p$ and
write $\bar{a} = a_1,\ldots,a_p$, then by (\ref{avect}),
\begin{equation} \label{avect2}
\vdash_{r_{i+2}} a_i: C_i[\bar{a}^{i-1}/\Theta^{i-1}] \; (\Gamma)
\end{equation} 
for $i=1,\ldots, p$. Then applying substitution $\bar{b}$, 
\begin{equation} \label{avect3}
a_i[\bar{b}/\Gamma]: C_i[\bar{a}^{i-1}/\Theta^{i-1}][\bar{b}/\Gamma] \; (\Delta).
\end{equation}
i.e.\
\begin{equation} \label{avect4}
a_i[\bar{b} /\Gamma]: 
C_i[a_1[\bar{b}/\Gamma],\ldots, a_{i-1}[\bar{b}/\Gamma]/z_1,\ldots,z_{i-1}] \; (\Delta).
\end{equation}
So 
\begin{equation} \label{avect5}
(a_1[\bar{b}/\Gamma],\ldots, a_p[\bar{b}/\Gamma]):  \Delta \to \Theta
\end{equation}

We have 
$$\sigma(\Delta, A[\bar{b}/\Gamma])  
 = \sigma(\Delta,  S(a_1,\ldots,a_p)[\bar{b}/\Gamma]) 
= \sigma(\Delta,  S(a_1[\bar{b}/\Gamma],\ldots,a_p[\bar{b}/\Gamma]))
$$
By (M8') and  (\ref{avect5}) the right hand side is
$$\sigma(\Theta,S(\bar{z}))\{[\theta(\Delta, C_1, a_1[\bar{b}/\Gamma]),\ldots, \theta(\Delta, C_p[a_1[\bar{b}/\Gamma],\ldots,a_{p-1}[\bar{b}/\Gamma]/z_1,\ldots,z_{p-1}],a_p[\bar{b}/\Gamma])]\}$$
which is
$$\sigma(\Theta,S(\bar{z}))\{[\theta(\Delta, C_1, a_1[\bar{b}/\Gamma]),\ldots, \theta(\Delta, C_p[a_1,\ldots,a_{p-1}/z_1,\ldots,z_{p-1}][\bar{b}/\Gamma],a_p[\bar{b}/\Gamma])]\}$$
By (\ref{avect}) using $C_1 = C_1[\bar{b}/\Gamma]$ we can apply the inductive hypothesis 
to get
$$
\begin{aligned}
& \sigma(\Theta,S(\bar{z}))\Bigl\{[\theta(\Gamma, C_1, a_1) \{[\theta(\Delta,B_1,b_1),\ldots,\theta(\Delta,B_n[\bar{b}^{n-1}/\Gamma^{n-1}], b_n)]\}, \\
& \ldots, \\
& \theta(\Gamma, C_p[a_1,\ldots,a_{p-1}/z_1,\ldots,z_{p-1}],a_p)
\{[\theta(\Delta,B_1,b_1),\ldots,\theta(\Delta,B_n[\bar{b}^{n-1}/\Gamma^{n-1}], b_n)]\}
]\Bigr\}
\end{aligned}
$$
By distributivity and functoriality we obtain  this expression
$$
\begin{aligned}
& \sigma(\Theta,S(\bar{z}))\Bigl\{[\theta(\Gamma, C_1, a_1) , \ldots, \theta(\Gamma, C_p[a_1,\ldots,a_{p-1}/z_1,\ldots,z_{p-1}],a_p)]\Bigr\} \\
& \qquad \qquad \quad
\Bigl\{[\theta(\Delta,B_1,b_1),\ldots,\theta(\Delta,B_n[\bar{b}^{n-1}/\Gamma^{n-1}], b_n)]\Bigr\}
\end{aligned}
$$
Then using (M8') again we obtain
$$\sigma(\Delta, A[\bar{b}/\Gamma])  = \sigma(\Gamma,S(a_1,\ldots,a_p))\Bigl\{[\theta(\Delta,B_1,b_1),\ldots,\theta(\Delta,B_n[\bar{b}^{n-1}/\Gamma^{n-1}], b_n)]\Bigr\}
$$
But since $A=S(a_1,\ldots,a_p)$ this is what was required. We have proved (M4') for $k$.

\medskip
Ad (M5'): Suppose $\vdash_k a:A \;  (\Gamma)$,
 Thus $a = f(a_1,\ldots,a_p)$ and $A =U[a_1,\ldots,a_p/z_1,\ldots,z_p]$ for some $(\Theta,f,U) \in \Sigma$ with
\begin{equation} \label{avectb}
\vdash_{\bar{r}}(a_1,\ldots,a_p): \Gamma \to \Theta
\end{equation}
and $\bar{r} < k$,  where  $\Theta = z_1:C_1,\ldots, z_p:C_p$.
Write $\bar{a} = a_1,\ldots,a_p$, then by (\ref{avectb}),
\begin{equation} \label{avect2b}
\vdash_{r_{i+2}} a_i: C_i[\bar{a}^{i-1}/\Theta^{i-1}] \; (\Gamma)
\end{equation} 
for $i=1,\ldots, p$. Then applying substitution $\bar{b}$, 
\begin{equation} \label{avect3b}
a_i[\bar{b}/\Gamma]: 
C_i[\bar{a}^{i-1}/\Theta^{i-1}][\bar{b}/\Gamma] \; (\Delta).
\end{equation}
i.e.\
\begin{equation} \label{avect4b}
a_i[\bar{b}/\Gamma]: C_i[a_1[\bar{b}/\Gamma],\ldots, a_{i-1}[\bar{b}/\Gamma]/z_1,\ldots,z_{i-1}] \; (\Delta).
\end{equation}
So 
\begin{equation} \label{avect5b}
(a_1[\bar{b}/\Gamma],\ldots, a_p[\bar{b}/\Gamma]):  \Delta \to \Theta
\end{equation}
We have 
$$
\begin{aligned}
& \theta(\Delta, U[a_1[\bar{b}/\Gamma],\ldots,a_p[\bar{b}/\Gamma]/z_1,\ldots,z_{p-1}], a[\bar{b}/\Gamma])  \\
& = \theta(\Delta,  U[a_1[\bar{b}/\Gamma],\ldots,a_p[\bar{b}/\Gamma]/z_1,\ldots,z_{p-1}], f(a_1,\ldots,a_p)[\bar{b}/\Gamma]) \\
& = \theta(\Delta, U[a_1[\bar{b}/\Gamma],\ldots,a_p[\bar{b}/\Gamma]/z_1,\ldots,z_{p-1}],  f(a_1[\bar{b}/\Gamma],\ldots,a_p[\bar{b}/\Gamma])).
\end{aligned}
$$

By (M9') and (\ref{avect5b}) the last expression equals
$$\theta(\Theta,U, f(\bar{z}))\{[\theta(\Delta, C_1, a_1[\bar{b}/\Gamma]),\ldots, \theta(\Delta, C_p[a_1[\bar{b}/\Gamma],\ldots,a_{p-1}[\bar{b}/\Gamma]/z_1,\ldots,z_{p-1}],a_p[\bar{b}/\Gamma])]\}$$
By (\ref{avectb}) using $C_1 = C_1[\bar{b}/\Gamma]$ we can apply the inductive hypothesis 
to get
$$
\begin{aligned}
& \theta(\Theta,U,f(\bar{z}))\Bigl\{[\theta(\Gamma, C_1, a_1) \{[\theta(\Delta,B_1,b_1),\ldots,\theta(\Delta,B_n[\bar{b}^{n-1}/\Gamma^{n-1}], b_n)]\}, \\
& \ldots, \\
& \theta(\Gamma, C_p[a_1,\ldots,a_{p-1}/z_1,\ldots,z_{p-1}],a_p)
\{[\theta(\Delta,B_1,b_1),\ldots,\theta(\Delta,B_n[\bar{b}^{n-1}/\Gamma^{n-1}], b_n)]\}
]\Bigr\}
\end{aligned}
$$
By distributivity and functoriality we obtain  this expression
$$
\begin{aligned}
& \theta(\Theta,U,f(\bar{z}))\Bigl\{[\theta(\Gamma, C_1, a_1) , \ldots, \theta(\Gamma, C_p[a_1,\ldots,a_{p-1}/z_1,\ldots,z_{p-1}],a_p)]\Bigr\} \\
& \qquad \qquad \quad
\Bigl\{[\theta(\Delta,B_1,b_1),\ldots,\theta(\Delta,B_n[\bar{b}^{n-1}/\Gamma^{n-1}], b_n)]\Bigr\}
\end{aligned}
$$
Then using (M9') again we obtain
$$\begin{aligned}
&\theta(\Delta, A[\bar{b}/\Gamma], a[\bar{b}/\Gamma])   \\
&= \theta(\Gamma,U[a_1,\ldots,a_p/z_1,\ldots,z_p], f(a_1,\ldots,a_p))\Bigl\{[\theta(\Delta,B_1,b_1),\ldots,\theta(\Delta,B_n[\bar{b}^{n-1}/\Gamma^{n-1}], b_n)]\Bigr\}.
\end{aligned}
$$
But since $a=f(a_1,\ldots,a_p)$ and $A =U[a_1,\ldots,a_p/z_1,\ldots,z_p]$ 
this is what was required. We have proved (M5') for $k$.

\medskip
 The definition (\ref{Fdef}) then gives a functor according to
Remark \ref{FunctorRem}. $\qed$

\medskip
For two signatures with $\Sigma \subseteq \Sigma'$ there is a canonical
embedding $E: {\cal F}_{\Sigma} \rightarrow {\cal F}_{\Sigma'}$ given
by letting  $E = (E,\sigma, \theta)$ be the identity on all input
\begin{eqnarray*}
E(\Gamma) &=&\Gamma, \\
E(\Delta,\Gamma,\bar{b}) &=&  (\Delta,\Gamma,\bar{b}),\\
\sigma_\Gamma({\rm A}) &=& {\rm A}, \\
\theta_{\Gamma,A}({\rm a}) &=& {\rm a}.
\end{eqnarray*}
The conditions follows readily since by 
Lemma \ref{Jmono}, ${\cal J}(\Sigma) \subseteq {\cal J}(\Sigma')$.

\subsection{Extension problem for signatures}

The {\em extension problem} is given a cwf 
morphism $G = (G, \sigma, \theta): {\cal F}_{\Sigma} \rightarrow {\cal C}$,
a signature $\Sigma'$ including $\Sigma$, and some suitable data in ${\cal C}$, find
a cwf morphism $G' = (G', \sigma', \theta'): {\cal F}_{\Sigma'} \rightarrow {\cal C}$ 
mapping the new syntactic entities to this data and such that $G' \circ E = G$,
where $E$ is the canonical cwf embedding ${\cal F}_{\Sigma} \rightarrow {\cal F}_{\Sigma'}$.
 
\medskip

We have the following general main cases of the extension problem:

\begin{theorem} \label{extbytype}
Given a signature $\Sigma$ and a cwf 
morphism $G = (G, \sigma, \theta): {\cal F}_{\Sigma} \rightarrow {\cal C}$. Moreover suppose
that $S$ is a type symbol not declared in $\Sigma$ and 
$(\Gamma_S \; {\rm context}) \in {\cal J}(\Sigma)$, so that $\Sigma'= \Sigma \cup \{(\Gamma_S, S)\}$ is a signature.
For $A \in \Ty_{\cal C}(G(\Gamma_S))$, there is a unique  cwf 
morphism $G' = (G', \sigma', \theta'): {\cal F}_{\Sigma'} \rightarrow {\cal C}$ 
with  $G' \circ E = G$ and 
$$\sigma'(\Gamma_S, S(\bar{x})) = A,$$
where $\bar{x} = {\rm OV}(\Gamma_S)$.
\end{theorem}
{\bf \flushleft Proof.}  First we observe that uniqueness is direct by Lemma \ref{cwfuniq}.
If $G''$ is another cwf morphism ${\cal F}_{\Sigma'} \rightarrow {\cal C}$ with 
$G'' \circ E = G$ and $\sigma'(\Gamma_S, S(\bar{x})) = A$, then $G''$ and $G'$ agree 
on the signature, so by this lemma $G'' = G'$.

We are going to define approximations of $(G', \sigma', \theta')$ for 
each maximum proof height $n \in {\mathbb N}$:
\begin{itemize}
\item[] $G'_n : {\rm Ob}^{\Sigma'}_n  \to {\rm Ob} \, {\cal C}$, 
\item[] $\sigma'_n : {\rm Ty}^{\Sigma'}_n(\Gamma) \to {\rm Ty}^{\cal C}(G'_n(\Gamma))$ 
for $\Gamma \in {\rm Ob}^{\Sigma'}_n$,
\item[] $\theta'_n: {\rm Tm}^{\Sigma'}_n(\Gamma,{\rm A}) \to {\rm Tm}^{\cal C}(G'_n(\Gamma),\sigma'_n(\Gamma,A))$
for $\Gamma \in {\rm Ob}^{\Sigma'}_n$,  $(\Gamma,A) \in {\rm Ty}^{\Sigma'}_n(\Gamma)$.
\end{itemize}
Suppose that these
 $G'_k$, $\sigma'_k$ and $\theta'_k$ have the properties that for any $p<k$, 
\begin{itemize}
\item[(I1)] $G'_p(\Gamma) = G'_k(\Gamma)$ for all $\Gamma \in {\rm Ob}^{\Sigma'}_p$,
\item[(I2)] $\sigma'_p(\Gamma, A) = \sigma'_k(\Gamma, A)$ for all $\Gamma \in {\rm Ob}^{\Sigma'}_p$,
$(\Gamma,A) \in {\rm Ty}^{\Sigma'}_p(\Gamma)$,
\item[(I3)] $\theta'_p((\Gamma, A), a) = \theta'_k((\Gamma, A), a)$ for all $\Gamma \in {\rm Ob}^{\Sigma'}_p$, $(\Gamma,A) \in {\rm Ty}^{\Sigma'}_p(\Gamma)$, $((\Gamma, A), a) \in  
{\rm Tm}^{\Sigma'}_p(\Gamma,{\rm A})$.
\end{itemize}
Then these independence properties imply that the following are well defined:
\begin{itemize}
\item[] $G'(\Gamma) = G'_k(\Gamma)$, where $\Gamma \in {\rm Ob}^{\Sigma'}_k$,
\item[] $\sigma'(\Gamma, A) = \sigma'_k(\Gamma, A)$, where $\Gamma \in {\rm Ob}^{\Sigma'}_k$,  $(\Gamma,A) \in {\rm Ty}^{\Sigma'}_k(\Gamma)$,
\item[] $\theta'((\Gamma, A), a)= \theta'_k((\Gamma, A), a)$, where $\Gamma \in {\rm Ob}^{\Sigma'}_k$, $(\Gamma,A) \in {\rm Ty}^{\Sigma'}_k(\Gamma)$, $((\Gamma, A), a) \in  
{\rm Tm}^{\Sigma'}_k(\Gamma,{\rm A})$.
\end{itemize}

\medskip

Note that ${\rm Ob}^{\Sigma'}_0 = \{\langle \rangle\}$ but 
${\rm Ty}^{\Sigma'}_n(\Gamma) = \emptyset$ and 
${\rm Tm}^{\Sigma'}_n(\Gamma,{\rm A})= \emptyset$. We
define  $G'_0(\langle\rangle) = \top$
and let $\sigma'_0$ and $\theta_0$ be the unique 
functions with empty domain.

Now define functions on level $n+1$:
\begin{itemize}
\item $G'_{n+1}(\langle\rangle) = \top$,
\item $G'_{n+1}(\Gamma,u:C) = G'_n(\Gamma).\sigma'_n(\Gamma,C)$ for 
$(\Gamma, u:C\; {\rm context}) \in {\cal J}_{n+1}(\Sigma')$,
\item $\sigma'_{n+1}(\Delta, S(\bar{t})) = A\Bigl\{\bar{\theta'}_n(\Delta,\Gamma_S,\bar{t})\Bigr\}$ where 
$(\bar{t}: \Delta \to \Gamma_S) \in {\cal J}_n(\Sigma')$,
\item $\sigma'_{n+1}(\Delta, T(\bar{t})) = \sigma(\Gamma_T, T(\bar{y}))\Bigl\{\bar{\theta'}_n(\Delta,\Gamma_T,\bar{t})\Bigr\}$ if
$T \ne S$ and $(\bar{t}: \Delta \to \Gamma_T) \in {\cal J}_n(\Sigma')$  and $(\Gamma_T,T) \in \Sigma$, and
where $\bar{y} = {\rm OV}(\Gamma_T)$,
\item $\theta'_{n+1}(\Delta, C_i, z_i) =  {\sf x}_i$ if $\Delta = z_1:C_1,\ldots,z_m:C_m$, if
$(\Delta \; {\rm context}) \in {\cal J}_n(\Sigma')$,
\item $\theta'_{n+1}(\Delta, U_f[\bar{t}/\Gamma_f], f(\bar{t})) = \theta(\Gamma_f,U_f,f(\bar{y}))\Bigl\{\bar{\theta'}_n(\Delta,\Gamma_f,\bar{t})\Bigr\}$ if 
$(\bar{t}: \Delta \to \Gamma_f) \in {\cal J}_n(\Sigma')$, $(U_f[\bar{t}/\Gamma_f]\; {\rm type}\; (\Delta)) \in {\cal J}_n(\Sigma')$  and $(\Gamma_f,f,U_f) \in \Sigma$, and
where $\bar{y} = {\rm OV}(\Gamma_f)$.
\end{itemize}
Note that by properties (P1) -- (P3) the above recursive definition is well-defined. The independence properties (I1) -- (I3) are now verified by induction on $k$.
By the indepedence properties we have now
\begin{itemize}
\item[(G1)] $G'(\langle\rangle) = \top$,
\item[(G2)] $G'(\Gamma,u:C) = G'(\Gamma).\sigma'(\Gamma,C)$ for 
$(\Gamma, u:C\; {\rm context}) \in {\cal J}(\Sigma')$,
\item[(G3)] $\sigma'(\Delta, S(\bar{t})) = A\Bigl\{\bar{\theta'}(\Delta,\Gamma_S,\bar{t})\Bigr\}$ where 
$(\bar{t}: \Delta \to \Gamma_S) \in {\cal J}(\Sigma')$,
\item[(G4)] $\sigma'(\Delta, T(\bar{t})) = \sigma(\Gamma_T, T(\bar{y}))\Bigl\{\bar{\theta'}(\Delta,\Gamma_T,\bar{t})\Bigr\}$ if
$T \ne S$ and $(\bar{t}: \Delta \to \Gamma_T) \in {\cal J}(\Sigma')$  and $(\Gamma_T,T) \in \Sigma$, and
where $\bar{y} = {\rm OV}(\Gamma_T)$,
\item[(G5)] $\theta'(\Delta, C_i, z_i) =  {\sf x}_i$ if $\Delta = z_1:C_1,\ldots,z_m:C_m$, and
$(\Delta \; {\rm context}) \in {\cal J}(\Sigma')$,
\item[(G6)] $\theta'(\Delta, U_f[\bar{t}/\Gamma_f], f(\bar{t})) = \theta(\Gamma_f,U_f,f(\bar{y})) \Bigl\{\bar{\theta'}(\Delta,\Gamma_f,\bar{t})\Bigr\}$ if 
$(\bar{t}: \Delta \to \Gamma_f) \in {\cal J}(\Sigma')$, $(U_f[\bar{t}/\Gamma_f]\; {\rm type}\; (\Delta)) \in {\cal J}(\Sigma')$  and $(\Gamma_f,f,U_f) \in \Sigma$, and
where $\bar{y} = {\rm OV}(\Gamma_f)$.
\end{itemize}
Note that using identity maps  in the above we get
\begin{equation} \label{eq56}
\sigma'(\Gamma_S, S(\bar{x})) = A
\end{equation}
\begin{equation}  \label{eq57}
\sigma'(\Gamma_T, T(\bar{y})) = \sigma(\Gamma_T, T(\bar{y})) \qquad T \ne S
\end{equation}
\begin{equation}  \label{eq58}
\theta'(\Gamma_f, U_f, f(\bar{y})) = \theta(\Gamma_f,U_f,f(\bar{y}))
\end{equation}

We now show that this extends to a cwf morphism $(G',\sigma',\theta'): {\cal F}_{\Sigma} \rightarrow {\cal C}$ by checking the conditions of Theorem \ref{thm41}. Conditions (M1) -- (M3) are clear by the above, condition (M7) is true by definition. We need to check the remaining (M8') and (M9').

(M8'): This follows from (\ref{eq56}) and (G3) for the case $S$ and from (\ref{eq57}) and (G4) for the general case.

(M9'): This follows from (\ref{eq58}) and (G6).

Hence by Theorem \ref{thm41}, there is a cwf morphism $G': {\cal F}_{\Sigma'} \to {\cal C}$.
We need to check that $G' \circ E = G$ and $\sigma'(\Gamma_S, S(\bar{x})) = A$. The latter follows by (\ref{eq56}) and the former follows by Lemma \ref{cwfuniq} and (\ref{eq57}) and
(\ref{eq58}). $\qed$

\begin{theorem} \label{extbyfun}
Given a signature $\Sigma$ and a cwf 
morphism $G = (G, \sigma, \theta): {\cal F}_{\Sigma} \rightarrow {\cal C}$. Moreover suppose
that $f$ is a function symbol not declared in $\Sigma$, and
$(U_f \; {\rm type}\; (\Gamma_f)) \in {\cal J}(\Sigma_f)$   so that $\Sigma'= \Sigma \cup \{(\Gamma_f, f, U_f)\}$ is a signature.
For $a \in \Tm_{\cal C}(G(\Gamma_f), \sigma_{\Gamma_f}(\Gamma_f, U_f))$, there is a unique cwf 
morphism $G' = (G', \sigma', \theta'): {\cal F}_{\Sigma'} \rightarrow {\cal C}$ with  $G' \circ E = G$
and 
$$\theta'(\Gamma_f, U_f, f(\bar{x})) = a,$$
where  $\bar{x} = {\rm OV}(\Gamma_f)$.
\end{theorem}
{\bf \flushleft Proof.}  Uniqueness follows as in  Theorem \ref{extbytype}
by Lemma \ref{cwfuniq}. 

As in Theorem \ref{extbytype}, we
define  $G'_0(\langle\rangle) = \top$
and let $\sigma'_0$ and $\theta_0$ be the unique 
functions with empty domain.
Define functions on level $n+1$:
\begin{itemize}
\item $G'_{n+1}(\langle\rangle) = \top$,
\item $G'_{n+1}(\Gamma,u:C) = G'_n(\Gamma).\sigma'_n(\Gamma,C)$ for 
$(\Gamma, u:C\; {\rm context}) \in {\cal J}_{n+1}(\Sigma')$,
\item $\sigma'_{n+1}(\Delta, T(\bar{t})) = \sigma(\Gamma_T, T(\bar{y}))\Bigl\{\bar{\theta'}_n(\Delta,\Gamma_T,\bar{t})\Bigr\}$ if $(\bar{t}: \Delta \to \Gamma_T) \in {\cal J}_n(\Sigma')$  and $(\Gamma_T,T) \in \Sigma$, and
where $\bar{y} = {\rm OV}(\Gamma_T)$,
\item $\theta'_{n+1}(\Delta, C_i, z_i) =  {\sf x}_i$ if $\Delta = z_1:C_1,\ldots,z_m:C_m$, and
$(\Delta \; {\rm context}) \in {\cal J}_n(\Sigma')$,
\item $\theta'_{n+1}(\Delta, U_f[\bar{t}/\Gamma_f], f(\bar{t})) = a\Bigl\{\bar{\theta'}_n(\Delta,\Gamma_f,\bar{t})\Bigr\}$ if $(\bar{t}: \Delta \to \Gamma_f) \in {\cal J}_n(\Sigma')$
\item $\theta'_{n+1}(\Delta, U_g[\bar{t}/\Gamma_g], g(\bar{t})) = \theta(\Gamma_g,U_g,g(\bar{y}))\Bigl\{\bar{\theta'}_n(\Delta,\Gamma_g,\bar{t})\Bigr\}$ if 
$(\bar{t}: \Delta \to \Gamma_f) \in {\cal J}_n(\Sigma')$, $(U_g[\bar{t}/\Gamma_g]\; {\rm type}\; (\Delta)) \in {\cal J}_n(\Sigma')$  and $(\Gamma_g,g,U_g) \in \Sigma$, $g \ne f$ and
where $\bar{y} = {\rm OV}(\Gamma_g)$.
\end{itemize}

Similarly to Theorem \ref{extbytype} we construct $(G',\sigma',\theta')$ so that
\begin{itemize}
\item[(G1)] $G'(\langle\rangle) = \top$
\item[(G2)] $G'(\Gamma,u:C) = G'(\Gamma).\sigma'(\Gamma,C)$, for $(\Gamma, u:C\; {\rm context}) \in {\cal J}(\Sigma')$,
\item[(G3)] $\sigma'(\Delta, T(\bar{t})) = \sigma(\Gamma_T, T(\bar{y}))\Bigl\{
\bar{\theta'}(\Delta,\Gamma_T,\bar{t})\Bigr\}$ if
$(\bar{t}: \Delta \to \Gamma_T) \in {\cal J}(\Sigma')$  and $(\Gamma_T,T) \in \Sigma$,
 and where $\bar{y} = {\rm OV}(\Gamma_T)$.
\item[(G4)] $\theta'(\Delta, C_i, z_i) =  {\sf x}_i$ if $\Delta = z_1:C_1,\ldots,z_m:C_m$, and
$(\Delta \; {\rm context}) \in {\cal J}(\Sigma')$,
\item[(G5)] 
$\theta'(\Delta, U_f[\bar{t}/\Gamma_f], f(\bar{t})) = a \Bigl\{\bar{\theta'}(\Delta,\Gamma_f,\bar{t})\Bigr\}$, 
if $(\bar{t}: \Delta \to \Gamma_f) \in {\cal J}(\Sigma')$,
\item[(G6)] $\theta'(\Delta, U_g[\bar{t}/\Gamma_g], g(\bar{t})) = \theta(\Gamma_g,U_g,g(\bar{y}))\Bigl\{
\bar{\theta'}(\Delta,\Gamma_g,\bar{t})\Bigr\}$ if 
$(\bar{t}: \Delta \to \Gamma_g) \in {\cal J}(\Sigma')$ and $(\Gamma_g,g,U_g) \in \Sigma$, $g \ne f$ and
where $\bar{y}= {\rm OV}(\Gamma_g)$.
\end{itemize}

The verification of the properties of $G'$ is also similar to Theorem \ref{extbytype}.
Using identity maps  in the above we get
\begin{equation} \label{eq55}
\sigma'(\Gamma_T, T(\bar{y})) = \sigma(\Gamma_T, T(\bar{y}))
\end{equation}
\begin{equation} \label{eq56}
\theta'(\Gamma_f, U_f, f(\bar{x})) = a
\end{equation}
\begin{equation} \label{eq57}
\theta'(\Gamma_g, U_g, g(\bar{y})) = \theta(\Gamma_g,U_g,g(\bar{y})) \qquad \text{ for } g\ne f
\end{equation}
We  show that this extends to a cwf morphism $(G',\sigma',\theta'): {\cal F}_{\Sigma} \rightarrow {\cal C}$ by checking the conditions of Theorem \ref{thm41}. Conditions (M1) -- (M3) are clear by the above, condition (M7) is true by definition. We need to check the remaining (M8') and (M9').

(M8'): this follows from (\ref{eq55}) and (G3).

(M9'): this is a consequence of (\ref{eq56}) and (\ref{eq57}) using respectively (G5) and (G6).

Hence by Theorem \ref{thm41}, there is a cwf morphism $G': {\cal F}_{\Sigma'} \to {\cal C}$.
We check that $G' \circ E = G$ and  $\theta'(\Gamma_f, U_f, f(\bar{x})) = a$. The latter follows by (\ref{eq56}), and the former is a consequence  of Lemma \ref{cwfuniq} using
(\ref{eq55}) and (\ref{eq57}). $\qed$

\section{Hyperdoctrines over cwfs} \label{semcomplsec}
 
 We introduce a logical structure of propositions over the type structure provided by cwfs
 using the technique of hyperdoctrines (cf.\ Seely 1983).
 
 Recall that a {\em Heyting algebra}  is a reflexive, transitive order $\le$ with operations $\land,\lor,\rightarrow$ and constants $\bot,\top$ such that for all $x,y,z$
 \begin{itemize}
 \item[(a)] $\bot \le x \le \top$,
 \item[(b)] $z \le x \land y$ iff $z \le x$ and $z \le y$,
 \item[(c)] $x \lor y \le z$ iff $x \le z$ and $y \le z$,
 \item[(d)] $z \le (x \rightarrow y)$ iff $z \land x \le y$.
 \end{itemize}
 Note that we are not requiring $\le$ to be antisymmetric, so we should probably
 say Heyting prealgebras rather than Heyting algebras.  A morphism between Heyting algebras is an order preserving homomorphism which respects the operations and constants. 
  
We suggest a generalization of hyperdoctrines to dependent types:
semantics of first-order logic  over a cwf can be given by the following data
\begin{enumerate}
\item A category with families ${\cal C}$.
\item A functor ${\rm Pr}: {\cal C}^{op} \to {\rm Heyting}$ into the category of Heyting algebras. For $f: \Delta \to \Gamma$ we write $$R\{f\}=_{\rm def } {\rm Pr}(f)(R).$$

\item For any $\Gamma \in {\rm Ob}({\cal C})$ and $S \in \Ty(\Gamma)$ monotone operations 
$\forall_{S},\exists_{S}: {\rm Pr}(\Gamma.S) \to {\rm Pr}(\Gamma)$ such that
\begin{enumerate}
\item For $Q \in {\rm Pr}(\Gamma)$, $R \in {\rm Pr}(\Gamma.S)$,
 $$Q \le \forall_{S}(R) \Longleftrightarrow  Q\{{\rm p}(S)\} \le R$$
\item For $Q \in {\rm Pr}(\Gamma)$, $R \in {\rm Pr}(\Gamma.S)$,
 $$\exists_{S}(R) \le Q \Longleftrightarrow R \le  Q\{{\rm p}(S)\}.$$
\end{enumerate}
\item For the pullback square
 \begin{equation} \label{BCC}
\bfig\square<1000,600>[\Delta.S\{f\}`\Gamma.S`\Delta`\Gamma;f.S`{\rm p}(S\{f\})`{\rm p}(S)`f]\efig
\end{equation}
we have for $R \in {\rm Pr}(\Gamma.S)$, the following {\em Beck-Chevalley conditions:}
\begin{enumerate}
\item $\forall_{S}(R)\{f\} =  \forall_{S\{f\}}(R\{f.S\})$,
\item  $\exists_{S}(R)\{f\}  = \exists_{S\{f\}}(R\{f.S\})$.
\end{enumerate}
\end{enumerate}
 
The resulting structure ${\cal H}=({\cal C},{\rm Pr}, \forall, \exists)$ is called a {\em first-order hyperdoctrine over $\cal C$}. Write ${\cal C}^{\cal H}={\cal C}$ and ${\rm Pr}^{\cal H}={\rm Pr}$. It is also meaningful to consider fragments of first-order logic: 
\begin{itemize}
\item  {\em cwf with Horn doctrine}: ${\rm Pr}$ is a contravariant functor into the category of lower semi lattices and we drop conditions (3) and (4).
\item  {\em cwf with regular doctrine}: ${\rm Pr}$ is a contravariant functor into the category of lower semi lattices and we drop conditions (3a) and (4a) and add the {\em Frobenius condition:}
For $Q \in {\rm Pr}(\Gamma)$, $R \in {\rm Pr}(\Gamma.S)$,
\begin{equation}\label{frob}
 Q \land \exists_{S}(R) \le 
 \exists_{S}(Q\{{\rm p}(S)\} \land R)
\end{equation}
(cf.\ Johnstone (2002)).
\item {\em cwf with coherent doctrine}: this is a cwf with regular doctrine, but ${\rm Pr}$ is a contravariant functor into the category of {\em distributive lattices}.
\end{itemize}
 
 \medskip
 Suppose that $F: {\cal C} \to {\cal C}'$ is a cwf morphism. Assume that
 ${\cal H}= ({\cal C}, {\rm Pr}, \forall,\exists)$
 and ${\cal H}'=({\cal C}', {\rm Pr}', \forall',\exists')$
 are two hyperdoctrines over the respective cwfs. An {\em $F$-based morphism of hyperdoctrines} 
 $$G: {\cal H} \to {\cal H}'$$ is a natural transformation $G: {\rm Pr} \to {\rm Pr}' \circ F$
 $$\bfig
 \square<800,500>[{\rm Pr}(\Gamma)`{\rm Pr}'(F\Gamma)`{\rm Pr}(\Delta)`{\rm Pr}'(F\Delta);G_\Gamma`\_\{f\}`\_\{Ff\}`G_\Delta]
 \efig
 $$
 such that for $R \in {\rm Pr}(\Gamma.S)$,
 \begin{enumerate}
\item $G_\Gamma(\forall_S(R))=\forall'_{\sigma_\Gamma(S)}(G_{\Gamma.S}(R))$,
\item $G_\Gamma(\exists_S(R))=\exists'_{\sigma_\Gamma(S)}(G_{\Gamma.S}(R))$.
\end{enumerate}

\subsection{Horn logic in any cwf} \label{horncwf}

We note that a canonical Horn doctrine arises from every cwf by a propositions-as-types interpretation. Let ${\cal C}$ be a cwf. Define for $\Gamma \in {\cal C}$,
$${\rm Pr}^\land_{\cal C}(\Gamma) = ({\rm F}(\Gamma),\le)$$
where ${\rm F}(\Gamma)$ consists of all finite sequences 
$\langle A_1,\ldots,A_n\rangle$ of types in $\Ty(\Gamma)$, which 
are ordered by $\le_\Gamma$ as follows
$$
\begin{aligned}
&\langle A_1,\ldots,A_n\rangle \le_\Gamma \langle B_1,\ldots,B_m\rangle \Longleftrightarrow \\
&\forall k=1,\ldots m, \text{ $\Tm(\Gamma. A_1. A_2\{{\rm p}\}.\ldots . A_n\{{\rm p}^{(n-1)}\}, B_k\{{\rm p}^{(n)}\})$
is inhabited}.
\end{aligned}$$
For $f: \Delta \to \Gamma$, let
\begin{equation} \label{AAf}
\langle A_1,\ldots,A_n\rangle\{f\} = \langle A_1\{f\},\ldots,A_n\{f\}\rangle.
\end{equation}
The semilattice structure is given by $\top_{\Gamma} = \langle\rangle$ and 
$$\langle A_1,\ldots,A_n\rangle \wedge_{\Gamma} \langle B_1,\ldots,B_m\rangle =
\langle A_1,\ldots,A_n, B_1,\ldots,B_m\rangle.$$

\begin{theorem} For any cwf ${\cal C}$, the pair $({\cal C}, {\rm Pr}^\land_{\cal C})$ constructed above is a cwf with Horn doctrine. 
\end{theorem}
{\flushleft \bf Proof.} We show that ${\rm Pr}^\land_{\cal C}$ is a lower semilattice. The element $\top_\Gamma$ is obviously the top element. To see that the relation is reflexive we note
that 
$$ {\rm v}_{A_k\{{\rm p}^{(k-1)}\}}\{{\rm p}^{(n-k)}\}\in \Tm(\Gamma. A_1. A_2\{{\rm p}\}.\ldots . A_n\{{\rm p}^{(n-1)}\}, A_k\{{\rm p}^{(n)}\})$$
for all $k=1,\ldots,n$. This also proves that for each $k$
\begin{equation} \label{A1AnAk}
\langle A_1,\ldots,A_n\rangle \le_\Gamma \langle A_k\rangle.
\end{equation}
Suppose that
$\langle A_1,\ldots,A_n\rangle \le_\Gamma \langle B_1,\ldots,B_m\rangle$
and 
$\langle B_1,\ldots,B_m\rangle \le_\Gamma \langle C_1,\ldots,C_p\rangle$.
Suppose that these relations are witnessed by
$$b_k \in \Tm(\Gamma. A_1. A_2\{{\rm p}\}.\ldots . A_n\{{\rm p}^{(n-1)}\}, B_k\{{\rm p}^{(n)}\})$$
$$c_k \in \Tm(\Gamma. B_1. B_2\{{\rm p}\}.\ldots . B_m\{{\rm p}^{(m-1)}\}, C_{\ell}\{{\rm p}^{(m)}\}).$$
Now 
$$f=\langle {\rm p}^{(n)}, b_1,\ldots, b_m\rangle:
\Gamma. A_1. A_2\{{\rm p}\}.\ldots . A_n\{{\rm p}^{(n-1)}\} \to 
B_1. B_2\{{\rm p}\}.\ldots . B_m\{{\rm p}^{(m-1)}\},$$
so
$$c_k\{f\} \in \Tm(\Gamma. A_1. A_2\{{\rm p}\}.\ldots . A_n\{{\rm p}^{(n-1)}\}, C_{\ell}\{{\rm p}^{(m)}\}\{f\}).$$
But $C_{\ell}\{{\rm p}^{(m)}\}\{f\} = C_{\ell}\{{\rm p}^{(n)}\}$ which proves
$\langle A_1,\ldots,A_n\rangle \le_\Gamma  \langle C_1,\ldots,C_p\rangle$.
By (\ref{A1AnAk}) and the definition of $\le_\Gamma$ it is easy to see that
$$\langle A_1,\ldots,A_n\rangle \wedge_{\Gamma} \langle B_1,\ldots,B_m\rangle
 \le_\Gamma \langle A_1,\ldots,A_n\rangle, \langle B_1,\ldots,B_m\rangle.$$
Suppose that 
$$\langle C_1,\ldots,C_p\rangle \le_\Gamma  \langle A_1,\ldots,A_n\rangle, \langle B_1,\ldots,B_m\rangle.$$
But then by the definition of the order
$$\langle C_1,\ldots,C_p\rangle \le_\Gamma  \langle A_1,\ldots,A_n, 
B_1,\ldots,B_m\rangle =  \langle A_1,\ldots,A_n\rangle \wedge_{\Gamma} \langle B_1,\ldots,B_m\rangle$$
as required.

The required properties of ${\rm Pr}^{\land}_{\cal C}$ follows directly from (\ref{AAf}) and the functoriality
of $\Ty$. $\qed$

\dontshow{
\medskip
Suppose that ${\cal C}$ is a contextual cwf with a Horn doctrine ${\rm Pr},\le$. 
From this define a regular doctrine ${\rm Pr}^*_{\cal C}$ as follows:
For $\Gamma \in {\rm Ob}\, {\cal C}$,  let ${\rm Pr}^*(\Gamma)$ consist of
pairs $(\Gamma.A_1.\ldots.A_n, P)$ such that 
$\Gamma.\Sigma \in {\rm Ob}\, {\cal C}$ and $P \in {\rm Pr}(\Gamma.\Sigma)$. Introduce the following relation:
$$(\Gamma.\Sigma, P) \le_\Gamma^* (\Gamma.\Theta, Q)$$
if, and only if, there is a context morphism $f: \Gamma.\Sigma 
\to \Gamma.\Theta$ with ${\rm p}^{(|\Theta|)} \circ f = {\rm p}^{(|\Sigma|)}$ and
$P \le_{\Gamma.\Sigma} Q\{f\}$. It is easily verified to be a preorder using the
functoriality of $\Ty$.

For $g: \Delta \to \Gamma$, define the action by $g$:
$$(\Gamma.\Sigma, P)\{g\}= (\Delta.\Sigma\{g\}, P\{g.\Sigma\})$$
which belongs to ${\rm Pr}^*(\Delta)$.

The semi lattice structure on ${\rm Pr}^*(\Gamma)$ is given by
$$\top^*_\Gamma = (\Gamma, \top_\Gamma)$$
and
$$(\Gamma.\Sigma, P) \land_\Gamma^* (\Gamma.\Theta, Q) = (\Gamma.\Sigma.\Theta\{{\rm p}^{(|\Sigma|)}\},  P\{{\rm p}^{(|\Theta|)}\} \land Q\{{\rm p}^{(|\Sigma|)}.\Theta\}).$$
Note that  $(\Gamma.\Sigma, P) \le_\Gamma^* \top^*_\Gamma$ follows taking 
$f= {\rm p}^{(|\Sigma|)}$ and by the property that $\top_\Gamma\{ {\rm p}^{(|\Sigma|)}\} =
\top_{\Gamma.\Sigma}$. Next we need
to check that
\begin{equation} \label{inf1}
(\Gamma.\Sigma.\Theta\{{\rm p}^{(|\Sigma|)}\},  P\{{\rm p}^{(|\Theta|)}\} \land Q\{{\rm p}^{(|\Sigma|)}.\Theta\})  \le^*_{\Gamma} (\Gamma.\Sigma, P)
\end{equation}
and 
\begin{equation} \label{inf2}
(\Gamma.\Sigma.\Theta\{{\rm p}^{(|\Sigma|)}\},  P\{{\rm p}^{(|\Theta|)}\} \land Q\{{\rm p}^{(|\Sigma|)}.\Theta\})  \le^*_{\Gamma} (\Gamma.\Theta, Q).
\end{equation}
The inequality (\ref{inf1}) is witnessed by $f={\rm p}^{(|\Theta|)}$. As for (\ref{inf2}) taking
$f={\rm p}^{(|\Sigma|)}.\Theta$ yields the equality.
Now suppose that 
$$(\Gamma.\Xi, R) \le^*_{\Gamma} (\Gamma.\Sigma, P)  \qquad   
(\Gamma.\Xi, R)  \le^*_{\Gamma} (\Gamma.\Theta, Q),$$
are witnessed by $g: \Gamma.\Xi \to \Gamma.\Sigma$ and
$h: \Gamma.\Xi \to \Gamma.\Theta$ respectively, with $R \le P\{g\}$ and $R\le Q\{h\}$,
and  moreover ${\rm p}^{(|\Sigma|)} \circ g= {\rm p}^{(|\Xi|)}$ and
${\rm p}^{(|\Theta|)} \circ h= {\rm p}^{(|\Xi|)}$. These two equations
and the
 pullback square
\begin{equation} \label{iterpqpb3}
\bfig\square<1000,600>[\Gamma.\Sigma.\Theta\{{\rm p}^{(|\Sigma|)}\} `\Gamma.\Theta`\Gamma.\Sigma`\Gamma;{\rm p}^{(|\Sigma|)}.\Theta`{\rm p}^{(|\Theta|)}`{\rm p}^{(|\Theta|)}`{\rm p}^{(|\Sigma|)}]\efig
\end{equation} 
yields a unique $f: \Gamma.\Xi \to \Gamma.\Sigma.\Theta\{{\rm p}^{(|\Sigma|)}\}$ with 
 ${\rm p}^{(|\Theta|)}  \circ f = g$ and  $({\rm p}^{(|\Sigma|)}.\Theta) \circ f = h$. This gives
 by substitution, $R \le P\{{\rm p}^{(|\Theta|)}\}\{f\}$ and  
 $R\le Q\{({\rm p}^{(|\Sigma|)}.\Theta)\}\{f\}$. Moreover 
 $$ {\rm p}^{(|\Sigma|+|\Theta|)} \circ f= 
 {\rm p}^{(|\Sigma|)}\circ {\rm p}^{(|\Theta|)} \circ f = {\rm p}^{(|\Sigma|)}\circ g= {\rm p}^{(|\Xi|)}.$$ 
 Hence 
 $$(\Gamma.\Xi, R)  \le^*_{\Gamma} (\Gamma.\Sigma, P) \land_\Gamma^* (\Gamma.\Theta, Q).$$
 
\medskip
For $S \in \Ty(\Gamma)$ and $(\Gamma.S.\Sigma, P) \in {\rm Pr}^*(\Gamma.S)$ 
let 
$$\exists_{S}(\Gamma.S.\Sigma, P) = (\Gamma.S.\Sigma, P)$$
which belongs to ${\rm Pr}^*(\Gamma)$. Now 
\begin{equation}\label{lhexists}
\exists_{S}(\Gamma.S.\Sigma, P) \le^*_{\Gamma} (\Gamma.\Theta, Q)
\end{equation}
is equivalent to the existence of $f: \Gamma.S.\Sigma \to \Gamma.\Theta$
with ${\rm p}^{(|\Theta|)} \circ f = {\rm p}^{(|\Sigma|+1)}$ and 
\begin{equation} \label{pqfineq}
P \le_{\Gamma.S.\Sigma} Q\{f\}.
\end{equation}
On the other hand
$$
(\Gamma.\Theta, Q)\{{\rm p}(S)\} = 
  (\Gamma.S. \Theta\{{\rm p}(S)\}, Q\{{\rm p}(S).\Theta\}$$
so
\begin{equation}\label{rhexists}
(\Gamma.S.\Sigma, P) \le^*_{\Gamma.S} (\Gamma.\Theta, Q)\{{\rm p}(S)\}
\end{equation}
amounts to the existence of some
$$g: \Gamma.S.\Sigma \to 
\Gamma.S. \Theta\{{\rm p}(S)\}
$$
with
${\rm p}^{(|\Theta|)} \circ g= {\rm p}^{(|\Sigma|)}$ and
\begin{equation} \label{ex60}
P \le_{\Gamma.S.\Sigma} Q\{{\rm p}(S).\Theta\}\{g\}=
Q\{({\rm p}(S).\Theta) \circ g\}.
\end{equation}
We show that (\ref{lhexists}) and (\ref{rhexists}) are equivalent. Assuming (\ref{rhexists})
we can take $f= ({\rm p}(S).\Theta) \circ g$ in (\ref{ex60}) and note that
\begin{eqnarray*}
{\rm p}^{(|\Theta|)} \circ f &=& {\rm p}^{(|\Theta|)} \circ  ({\rm p}(S).\Theta) \circ g \\
&=&  {\rm p}(S)  \circ {\rm p}^{(|\Theta|)}  \circ g\\
&=&  {\rm p}(S)  \circ {\rm p}^{(|\Sigma|)} = {\rm p}^{(|\Sigma|+1)}.
\end{eqnarray*}
This establishes (\ref{lhexists}). Conversely, assume (\ref{lhexists}). We use the fact that
the following is a pullback 
\begin{equation}
\bfig\square<1000,600>[\Gamma.S.\Theta\{{\rm p}(S)\}`\Gamma.\Theta`\Gamma.S`\Gamma;{\rm p}(S).\Theta`{\rm p}^{(|\Theta|)}`{\rm p}^{(|\Theta|)}`{\rm p}(S)]\efig
\end{equation}
Now by (\ref{lhexists}), ${\rm p}(S) \circ {\rm p}^{(|\Sigma|)} = {\rm p}^{(|\Sigma|+1)}={\rm p}^{(|\Theta|)} \circ f$,
so by the pullback above there is a unique
$$g: \Gamma.S.\Sigma \to 
\Gamma.S.\Theta\{{\rm p}(S)\}
$$
such that ${\rm p}^{(|\Theta|)} \circ g = {\rm p}^{(|\Sigma|)}$ and 
$({\rm p}(S).\Theta) \circ g =f$.
Then from (\ref{pqfineq}) follows (\ref{ex60}) so (\ref{lhexists}) holds as required.

Suppose $g: \Delta \to \Gamma$. Then
$$
\begin{aligned}
&\bigl((\Gamma.\Sigma, P) \land_\Gamma^* (\Gamma.\Theta, Q)\bigr)\{g\} \\
& \qquad \qquad = (\Gamma.\Sigma.\Theta\{{\rm p}^{(|\Sigma|)}\},  P\{{\rm p}^{(|\Theta|)}\} \land Q\{{\rm p}^{(|\Sigma|)}.\Theta\})\{g\} \\
& \qquad \qquad = (\Delta.(\Sigma.\Theta\{{\rm p}^{(|\Sigma|)}\})\{g\}, 
(P\{{\rm p}^{(|\Theta|)}\} \land Q\{{\rm p}^{(|\Sigma|)}.\Theta\})\{g.\Sigma.\Theta\{{\rm p}^{(|\Sigma|)}\}\})\\
& \qquad \qquad = (\Delta.\Sigma\{g\}.\Theta\{{\rm p}^{(|\Sigma|)} \circ g.\Sigma\},  P\{{\rm p}^{(|\Theta|)} \circ g.\Sigma.\Theta\{{\rm p}^{(|\Sigma|)}\}\} \land Q\{{\rm p}^{(|\Sigma|)}.\Theta \circ g.\Sigma.\Theta\{{\rm p}^{(|\Sigma|)}\}\}) \\
& \qquad \qquad = (\Delta.\Sigma\{g\}.\Theta\{g \circ {\rm p}^{(|\Sigma\{g\}|)}\}, P\{g.\Sigma \circ {\rm p}^{(|\Theta\{g\}|)}\}
\land Q\{g.\Theta \circ {\rm p}^{(|\Sigma\{g\})}.\Theta\{g\}\})\\
& \qquad \qquad = (\Delta.\Sigma\{g\}.\Theta\{g\}\{{\rm p}^{(|\Sigma\{g\}|)}\}, P\{g.\Sigma\}\{{\rm p}^{(|\Theta\{g\}|)}\}
\land Q\{g.\Theta\}\{{\rm p}^{(|\Sigma\{g\})}.\Theta\{g\}\}) \\
& \qquad \qquad = (\Delta.\Sigma\{g\}, P\{g.\Sigma\}) \land_\Delta^*(\Delta.\Theta\{g\}, Q\{g.\Theta\}) \\
& \qquad \qquad = (\Gamma.\Sigma, P)\{g\} \land_\Delta^* (\Gamma.\Theta, Q)\{g\}.
\end{aligned}
$$
The fourth equality follows by using the below pullback squares:

\begin{equation} \label{pb71}
\bfig\square<1000,600>[\Delta.\Sigma\{g\}`\Gamma.\Sigma`\Delta`\Gamma;g.\Sigma`{\rm p}^{(|\Sigma\{g\}|)}`{\rm p}^{(|\Sigma|)}`g]\efig
\end{equation}
\begin{equation} \label{pb72}
\bfig\square<1400,600>[\Delta.\Sigma\{g\}.\Theta\{{\rm p}^{(|\Sigma|)}\}\{g.\Sigma\}`\Gamma.\Sigma.\Theta\{{\rm p}^{(|\Sigma|)}\}`\Delta.\Sigma\{g\}`\Gamma.\Sigma;
g.\Sigma.\Theta\{{\rm p}^{(|\Sigma|)}\}`{\rm p}^{(|\Theta|)}`{\rm p}^{(|\Theta|)}`g.\Sigma]\efig
\end{equation}
\begin{equation} \label{pb73}
\bfig\square<1400,600>[\Gamma.\Sigma.\Theta\{{\rm p}^{(|\Sigma|)}\}`\Gamma.\Theta`\Gamma.\Sigma`\Gamma;
{\rm p}^{(|\Sigma|)}.\Theta`{\rm p}^{(|\Theta|)}`{\rm p}^{(|\Theta|)}`{\rm p}^{(|\Sigma|)}]\efig
\end{equation}
Combining the two squares (\ref{pb72}) and (\ref{pb73}) and functoriality we have
$${\rm p}^{(|\Sigma|)}.\Theta \circ g.\Sigma.\Theta\{{\rm p}^{(|\Sigma|)}\}
= ({\rm p}^{(|\Sigma|)} \circ g.\Sigma).\Theta$$
By (\ref{pb71}) 
$$({\rm p}^{(|\Sigma|)} \circ g.\Sigma).\Theta = (g \circ {\rm p}^{(|\Sigma\{g\}|)}).\Theta.$$
We have two further pullback squares
\begin{equation} \label{pb74}
\bfig\square<1400,600>[\Delta.\Sigma\{g\}.\Theta\{g\}.\{{\rm p}^{(|\Sigma\{g\}|)}\}`
\Delta.\Theta\{g\}`
\Delta.\Sigma\{g\}`
\Delta;{\rm p}^{(|\Sigma\{g\}|)}.\Theta\{g\}`{\rm p}^{(|\Theta\{g\}|)}`{\rm p}^{(|\Theta\{g\}|)}`{\rm p}^{(|\Sigma\{g\}|)}]
\efig
\end{equation}
\begin{equation} \label{pb75}
\bfig\square<1000,600>[\Delta.\Theta\{g\}`\Gamma.\Theta`\Delta`\Gamma;g.\Theta`{\rm p}^{(|\Theta\{g\}|)}`{\rm p}^{(|\Theta|)}`g]\efig
\end{equation}
Combining the two squares (\ref{pb74}) and (\ref{pb75}) and functoriality we have
$$(g \circ {\rm p}^{(|\Sigma\{g\}|)}).\Theta = g.\Theta \circ {\rm p}^{(|\Sigma\{g\}|)}.\Theta\{g\}.$$
Thus
$${\rm p}^{(|\Sigma|)}.\Theta \circ g.\Sigma.\Theta\{{\rm p}^{(|\Sigma|)}\}
=  g.\Theta \circ {\rm p}^{(|\Sigma\{g\}|)}.\Theta\{g\}$$
as required.

The Beck-Chevalley condition 4(b) follows by evaluating both sides of the equation
using the associated commutative square.

The Frobenius condition:  Suppose that 
$Q= (\Gamma.\Sigma,Q) \in {\rm Pr}(\Gamma)$ and
$R = (\Gamma.S.\Theta,R) \in {\rm Pr}(\Gamma.S)$.
We have to show that
\begin{equation} \label{frob2}
Q \land^* \exists_{S}(R) \le 
 \exists_{S}(Q\{{\rm p}(S)\} \land^* R).
\end{equation}
Now
$$\begin{aligned}
Q \land^*_{\Gamma} \exists_{S}(R) &=  
(\Gamma.\Sigma,Q) \land^*_{\Gamma} (\Gamma.S.\Theta,R) \\
&=
\bigl(\Gamma.\Sigma.(S.\Theta)\{{\rm p}^{(|\Sigma|)}\}, 
Q\{{\rm p}^{(|S.\Theta|)}\} \land R\{{\rm p}^{(|\Sigma|)}.S.\Theta\}\bigr)
\end{aligned}
$$
and 
$$\begin{aligned}
&\exists_{S}(Q\{{\rm p}(S)\} \land^* R) \\
&= 
(\Gamma.S.\Sigma\{{\rm p}(S)\}, Q\{{\rm p}(S).\Sigma\})\land^*_{\Gamma} (\Gamma.S.\Theta,R) \\
&= \bigl(\Gamma.S.\Sigma\{{\rm p}(S)\}.(S.\Theta)\{{\rm p}^{(|S.\Sigma\{{\rm p}(S)|)}.S.\Theta\},
Q\{{\rm p}(S).\Sigma\}\{{\rm p}^{(|S.\Theta|)}\} \land R\{{\rm p}^{(|S.\Sigma\{{\rm p}(S)|)}.S.\Theta\}
\bigr).
\\
\end{aligned}
$$
Let 
$$A=_{\rm def} \Gamma.\Sigma.(S.\Theta)\{{\rm p}^{(|\Sigma|)}\} = 
\Gamma.\Sigma.S\{{\rm p}^{(|\Sigma|)}\}.\Theta\{{\rm p}^{(|\Sigma|)}.S\}$$ 
and
$$B=_{\rm def} \Gamma.S.\Sigma\{{\rm p}(S)\}.(S.\Theta)\{{\rm p}^{(|S.\Sigma\{{\rm p}(S)|)}.S.\Theta\}.$$
To prove (\ref{frob2}) we use the following pullback squares:
\begin{equation} \label{pbsq01}
\bfig\square<1200,600>[B`\Gamma.S.\Theta`\Gamma.S.\Sigma\{{\rm p}\}`\Gamma;
{\rm p}^{(|S.\Sigma\{{\rm p}\}|)}.S.\Theta`{\rm p}^{(|S.\Theta|)}`{\rm p}^{(|S.\Theta|)}`{\rm p}^{(|S.\Sigma\{{\rm p}\}|)}]\efig
\end{equation}
\begin{equation} \label{pbsq02}
\bfig\square<1200,600>[\Gamma.S.\Sigma\{{\rm p}\}`\Gamma.\Sigma`\Gamma.S`\Gamma;
{\rm p}.\Sigma`{\rm p}^{(|\Sigma|)}`{\rm p}^{(|\Sigma|)}`{\rm p}]\efig
\end{equation}
\begin{equation} \label{pbsq03}
\bfig\square<1200,600>[\Gamma.\Sigma.S\{{\rm p}^{(|\Sigma|)}\}`\Gamma.S`\Gamma.\Sigma`\Gamma;
{\rm p}^{(|\Sigma|)}.S`{\rm p}`{\rm p}`{\rm p}^{(|\Sigma|)}]\efig
\end{equation}
By the square (\ref{pbsq03})  and the pullback (\ref{pbsq02}) we get a unique 
$t: \Gamma.\Sigma.S\{{\rm p}^{(|\Sigma|)}\} \to \Gamma.S.\Sigma\{{\rm p}\}$ such that
\begin{equation} \label{ptpS}
{\rm p}^{(|\Sigma|)} \circ t = {\rm p}^{(|\Sigma|)}.S
\qquad
{\rm p}.\Sigma \circ t = {\rm p}.
\end{equation}
Further we have the pullback square
\begin{equation} \label{pbsq04}
\bfig\square<1200,600>[A`\Gamma.S.\Theta`\Gamma.\Sigma`\Gamma;
{\rm p}^{(|\Sigma|)}.S.\Theta`{\rm p}^{(|S.\Theta|)}`{\rm p}^{(|S.\Theta|)}`{\rm p}^{(|\Sigma|)}]\efig
\end{equation}
We show that the square 
\begin{equation} \label{pbsq05}
\bfig\square<1200,600>[A`\Gamma.S.\Theta`\Gamma.S.\Sigma\{{\rm p}\}`\Gamma;
{\rm p}^{(|\Sigma|)}.S.\Theta`t{\rm p}^{(|\Theta|)}`{\rm p}^{(|S.\Theta|)}`{\rm p}^{(|S.\Sigma|)}]\efig
\end{equation}
commutes. Indeed,
\begin{eqnarray*}
{\rm p}^{(|S.\Theta|)} \circ {\rm p}^{(|\Sigma|)}.S.\Theta&=& {\rm p}^{(|\Sigma|)} \circ  {\rm p}^{(|S.\Theta|)} \qquad \text{by (\ref{pbsq04})}\\
&=& {\rm p}^{(|\Sigma|)} \circ {\rm p} \circ {\rm p}^{(|\Theta|)} \\
&=& {\rm p} \circ {\rm p}^{|\Sigma|}.S  \circ  {\rm p}^{(|\Theta|)} \qquad \text{by (\ref{pbsq03})} \\
&=& {\rm p} \circ {\rm p}^{(|\Sigma|)} \circ t  \circ  {\rm p}^{(|\Theta|)} \qquad \text{by (\ref{ptpS})} \\
&=& {\rm p}^{(|S.\Sigma|)} \circ t  \circ  {\rm p}^{(|\Theta|)} \\
\end{eqnarray*}
as required. Thus by the pullback square (\ref{pbsq01}) there is a unique $r: A \to B$ with
\begin{equation}
{\rm p}^{|S.\Theta|}  \circ r = t \circ  {\rm p}^{|\Theta|}
\end{equation}
\begin{equation}
{\rm p}^{(|S.\Sigma\{{\rm p}\}|)}.S.\Theta \circ r =  {\rm p}^{(|\Sigma|)}.S.\Theta\\
\end{equation}
Thus
$$
\begin{aligned}
& \bigl(Q\{{\rm p}(S).\Sigma\}\{{\rm p}^{(|S.\Theta|)}\} \land R\{{\rm p}^{(|S.\Sigma\{{\rm p}(S)|)}.S.\Theta\}\bigr)\{r\} \\
&\qquad \qquad \qquad = 
Q\{{\rm p}(S).\Sigma\}\{{\rm p}^{(|S.\Theta|)}\}\{r\} \land R\{{\rm p}^{(|S.\Sigma\{{\rm p}(S)|)}.S.\Theta\}\{r\}
  \\
&\qquad \qquad \qquad = 
Q\{{\rm p}(S).\Sigma\}\{t \circ  {\rm p}^{|\Theta|}\} \land R\{{\rm p}^{(|\Sigma|)}.S.\Theta\}
  \\
&\qquad \qquad \qquad = 
Q\{{\rm p}(S) \circ {\rm p}^{|\Theta|}\} \land R\{{\rm p}^{(|\Sigma|)}.S.\Theta\} \\
 &\qquad \qquad \qquad = 
Q\{{\rm p}^{|S.\Theta|}\} \land R\{{\rm p}^{(|\Sigma|)}.S.\Theta\}
  \end{aligned}
$$
as required. Furthermore
\begin{eqnarray*}
{\rm p}^{|S.\Sigma\{{\rm p}(S)\}.(S.\Theta)\{{\rm p}^{(|S.\Sigma\{{\rm p}(S)|)}.S.\Theta\}|}\circ r &=& {\rm p}^{|S.\Sigma\{{\rm p}(S)\}|} \circ t \circ {\bf p}^{|\Theta|} \\
&=& {\rm p} \circ {\rm p}^{(|\Sigma|)}.S  \circ {\bf p}^{|\Theta|} \\
&=& {\rm p}^{|\Sigma|} \circ {\rm p} \circ {\bf p}^{|\Theta|} \\
&=& {\rm p}^{|\Sigma.(S.\Theta)\{{\rm p}^{(|\Sigma|)}\}|}.
\end{eqnarray*}

This establishes the Frobenius condition.

\medskip
We summarize the above as

\begin{theorem} Let   $({\cal C}, {\rm Pr}_{\cal C})$ be any contextual cwf
with Horn doctrine. Then
${\rm Pr}^*_{\cal C}$ as defined above is a regular doctrine. $\qed$
\end{theorem}
}

\subsection{Interpretation into a category with families and type constructions} \label{intcwfsec}

\medskip
Let ${\cal C}$ be a cwf. Define for $\Gamma \in {\cal C}$,
$${\rm Pr}_{\cal C}(\Gamma) = (\Ty(\Gamma),\le)$$ where
$A \le B \Longleftrightarrow_{\rm def}  \text{$\Tm(\Gamma.A, B\{{\rm p}(A)\})$ is inhabited}.$
For $\sigma: \Delta \to \Gamma$, define
$${\rm Pr}_{\cal C}(f)(A) =_{\rm def} \Ty(f)(A).$$

\begin{lemma} For a cwf ${\cal C}$, ${\rm Pr}_{\cal C}$ is a functor ${\cal C}^{\rm op} \rightarrow {\rm Preorder}$.
\end{lemma}
{\flushleft \bf Proof.} Let $\Gamma$ be fixed. We show that ${\rm Pr}_{\cal C}(\Gamma)$ is a preorder. We have
${\rm v}_A \in \Tm(\Gamma.A, A\{{\rm p}(A)\})$, so $A \le A$.
Suppose $A \le B$ and $B \le C$, so that for some $s \in \Tm(\Gamma.A, B\{{\rm p}(A)\})$
and $t \in \Tm(\Gamma.B, C\{{\rm p}(B)\})$. Thus $\langle {\rm p}(A), s \rangle : \Gamma.A \rightarrow \Gamma.B$,
and hence by applying the $\Ty$ functor with this map
$$t\{\langle {\rm p}(A), s \rangle\}\in \Tm(\Gamma.A, C\{{\rm p}(B)\}\{\langle {\rm p}(A), s \rangle\})
= \Tm(\Gamma.A, C\{{\rm p}(A)\}).$$
Hence $A \le C$.

Suppose that $f: \Delta \rightarrow \Gamma$ and $A \le B$ 
in ${\rm Pr}_{\cal C}(\Gamma)$. Then $f.A: \Delta.A\{f\} \to \Gamma.A$
and $s \in \Tm(\Gamma.A, B\{{\rm p}(A)\})$, so 
$$s\{f.A\} \in \Tm(\Delta.A\{f\}, B\{{\rm p}(A)\}\{f.A\}) = \Tm(\Delta.A\{f\}, B\{f\}\{{\rm p}(A\{f\}) \rangle\}),
$$
that is $A\{f\} \le B\{f\}$. Functoriality of  ${\rm Pr}_{\cal C}$ follows immediately by the functoriality of $\Ty$. $\qed$

\begin{theorem} {\label{MLTThdoc}}  If ${\cal C}$ is a CwF which admits the type constructions $\Sigma, \Pi, +, {\rm N}_0$ and ${\rm N}_1$ then $({\cal C},{\rm Pr}_{\cal C})$ is a CwF with first-order doctrine.
\end{theorem}
{\flushleft \bf Proof.} The proof is essentially the standard proposition-as-types intepretation of logic in Martin-L\"of Type Theory. For details, see e.g.\ Appendix B of (Palmgren 2016). $\qed$

\section{Dependently typed first-order logic} \label{DFOL}

Here we introduce a version of the dependently typed first-order logic of (Belo 2007)
which generalizes the standard presentation of many-sorted first-order logic (Johnstone 2002). The first version DFOL uses capture avoiding substitutions and standardized contexts,
and is suitable for an easy completeness proof.
The second version DFOL* is closer to (Belo 2007) and is introduced and related to DFOL
in Section \ref{dfolstar}.

The predicate symbols are given by a discrete set $P$.
Suppose that $\Sigma$ is a term signature over a given symbol system $({\cal V},T,F)$ and 
$P$ is disjoint from $T \cup F$. A {\em predicate declaration relative to $\Sigma$} is a  triple $D=(\Gamma,({i_1},\ldots,{i_k}), R)$ where $\Gamma \in {\cal J}(\Sigma)$, ${i_1},\ldots,{i_k} \in {\rm DS}(\Gamma)$ and $R \in P$.  The symbol declared is
${\rm decl}(D)=_{\rm def} R$. If $\Gamma = x_1 : A_1,\ldots,x_n : A_n$, we also write
the declaration as
$$R(x_{i_1},\ldots,x_{i_k}) \quad (\Gamma).$$
The declaration is {\em regular} if $i_1,\ldots, i_k = 1, \ldots, n$.
A {\em predicate signature relative to $\Sigma$} is a set $\Pi$ of predicate declarations relative to $\Sigma$ satisfying the uniqueness condition $D=D'$ whenever ${\rm decl}(D) = {\rm decl}(D')$ for $D,D' \in \Pi$.
A first order signature is a pair $(\Sigma, \Pi)$ consisting of a term signature $\Sigma$ and
a predicate signature $\Pi$ relative to $\Sigma$.

Given this we can form the set of formulas for a signature $(\Sigma,\Pi)$. 
Note that we do not assume equality as a standard predicate. 
Let ${\rm Form}(\Sigma,\Pi)$ be the smallest set of judgement expressions 
$$\phi \; {\rm form}\; (\Gamma)$$ formed by
the rules: 
\begin{itemize}
\item[(F1)] For each predicate declaration $(\Delta, \bar{i},R)$ in $\Pi$,  and 
 $(\bar{a}: \Gamma \to \Delta) \in {\cal J}(\Sigma)$
we have
$$
\infer{R(\bar{a}_{i}) \;{\rm form} \; (\Gamma)}{}
$$ 
\item[(F2)] For $(\Gamma \;{\rm context}) \in {\cal J}(\Sigma)$,
$$\infer{ \bot \;{\rm form}  \; (\Gamma)}{}\qquad
\infer{ \top \;{\rm form}  \; (\Gamma)}{}$$
\item[(F3)] For  the connectives $\bigcirc= \land, \lor, \rightarrow$:
$$\infer{ (\phi \bigcirc \psi) \;{\rm form}  \; (\Gamma)}{ \phi  \;{\rm form}\; (\Gamma)  &
\psi  \;{\rm form}  \; (\Gamma) }
$$
\item[(F4)] For the quantifiers $Q= \forall,\exists$, and 
for $(A \; {\rm type} \; (\Gamma)) \in {\cal J}(\Sigma)$,
$$\infer{ (Q x:A)\phi  \;{\rm form}  \; (\Gamma)}{\phi  \;{\rm form}  \; (\Gamma,x:A)} 
$$
\end{itemize}
The free variables of formulas are defined by induction in the usual way, but where the
quantifier case is modified in order to take the variables of the type into account:
$$
{\rm FV}((Q x:A)\phi) = {\rm FV}(A) \cup ({\rm FV}(\phi)  \setminus \{x\}).
$$

\begin{remark} \label{boundfreedist}
{\em The free variables of the formula $\phi$ in context $\Gamma$ will all be declared in $\Gamma$.
Observe that by the formation rule for quantifiers, the bound and free variables of the formula $\phi$ will always be distinct. Also if  $((Q x:A)\phi  \;{\rm form}  \; (\Gamma)) \in {\rm Form}(\Sigma,\Pi)$, then $x$ cannot be bound in $\phi$. Thus for instance
$$(\forall x:A)R(x)\; {\rm form} \; (x:A)$$
is not well-formed.} 
\end{remark}

We have the following result about presuppositions, which is proved by induction on 
formation of formulas.

\begin{lemma}  \label{formpresup} If $(\phi\; {\rm form} \; (\Gamma)) \in {\rm Form}(\Sigma,\Pi)$, then 
$(\Gamma \; {\rm context}) \in {\cal J}(\Sigma)$. $\qed$
\end{lemma} 

We define inductively the notion of a formula $(\phi \; {\rm form} \; (\Gamma)) \in {\rm Form}(\Sigma, \Pi)$ being on {\em standard form (with respect to $\Sigma$)} as follows
\begin{itemize} 
\item for $\phi$ atomic, $(\phi \; {\rm form} \; (\Gamma))$ is on standard form, whenever
$\Gamma$ is a context on standard form with respect to $\Sigma$ (cf.\ Remark \ref{ctxstandard}).
\item if $(\phi \; {\rm form} \; (\Gamma))$ and $(\psi \; {\rm form} \; (\Gamma))$ are on standard form, then so is $(\phi \bigcirc \psi \; {\rm form} \; (\Gamma))$, for $\bigcirc = \land,\lor, \rightarrow$,
\item if $(\phi \; {\rm form} \; (\Gamma,x:A))$ is on standard form, then so is 
$((Qx:A) \phi \; {\rm form} \; (\Gamma))$, for $Q = \forall, \exists$.
\end{itemize}

The following is evident:

\begin{lemma} If $(\phi\; {\rm form} \; (\Gamma)) \in {\rm Form}(\Sigma,\Pi)$,  and $\Sigma$ has the de Bruijn property, then $(\phi\; {\rm form} \; (\Gamma))$ is on standard form
with respect to $\Sigma$. $\qed$
\end{lemma}

\medskip
{\em Capture avoiding substitution:} We define for $(\phi\; {\rm form} \; (\Gamma))
\in {\rm Form}(\Sigma,\Pi)$ and $(\bar{a}: \Delta \to \Gamma) \in {\cal J}(\Sigma)$ the capture avoiding substitution instance 
$$\phi\{(\Delta,\Gamma,\bar{a})\}$$ 
by induction on formulas:
\begin{itemize}
  \item if $(\phi\; {\rm form}\; (\Gamma)) \in {\rm Form}(\Sigma,\Pi)$ and $\phi$ is atomic, 
  then $$\phi\{(\Delta,\Gamma,\bar{a})\} = \phi[\bar{a}/\Gamma]$$ (syntactic substitution).
  \item if $(\phi\; {\rm form}\; (\Gamma)) \in {\rm Form}(\Sigma,\Pi)$ and 
  $(\psi\; {\rm form}\; (\Gamma)) \in {\rm Form}(\Sigma,\Pi)$, then
  $$(\phi \bigcirc \psi)\{(\Delta,\Gamma,\bar{a})\} =
      \phi\{(\Delta,\Gamma,\bar{a})\} \bigcirc \psi\{(\Delta,\Gamma,\bar{a})\}$$
      for $\bigcirc= \land, \lor, \rightarrow$.
  \item if $(\theta\; {\rm form}\; (\Gamma,x:A)) \in  {\rm Form}(\Sigma,\Pi)$, then
   for $Q=\forall, \exists$,
  $$((Qx:A)\theta)\{(\Delta,\Gamma,\bar{a})\} =
   (Qy:A[\bar{a}/\Gamma])\bigl(\theta\{(\langle\Delta,y:A[\bar{a}/\Gamma]\rangle,\langle \Gamma,x:A\rangle,(\bar{a},y))\}\bigr)$$
   where $y= {\rm fresh}(\Delta)$.
\end{itemize}

\begin{lemma} If $(\phi\; {\rm form}\; (\Gamma)) \in {\rm Form}(\Sigma,\Pi)$  and 
$(\bar{a}: \Delta \to \Gamma) \in {\cal J}(\Sigma)$, then
$(\phi\{(\Delta,\Gamma,\bar{a})\} \; {\rm form} \; (\Delta)) \in {\rm Form}(\Sigma,\Pi)$.
Note that the height of formulas are preserved under substitution.
\end{lemma}
{\flushleft \bf Proof.} By induction on formulas. $\qed$

\medskip
\begin{definition} \label{sequent} {\em 
For $(\phi\; {\rm form} \; (\Gamma)), (\psi \; {\rm form} \; (\Gamma)) 
\in {\rm Form}(\Sigma,\Pi)$, we have a judgement  expression called {\em sequent}
$$ \phi \sequent \psi \quad (\Gamma)$$
also written
$$ \phi \sequentunder{\Gamma} \psi$$
intended to mean that in the context $\Gamma$, the formula $\psi$ is true whenever $\phi$ is true.
A {\em theory} with respect to the signature $(\Pi,\Sigma)$ is a set of sequents over
 the same signature.}
\end{definition}

Let $T$ be a theory with respect to the signature $(\Pi,\Sigma)$. Let ${\rm Thm}(\Pi,\Sigma,T)$ denote the smallest set of sequents containing $T$ and closed under
the propositional and quantificational rules, and the substitution rules below.

{\em Propositional rules.} For $(\phi\; {\rm form} \; (\Gamma)), (\psi \; {\rm form} \; (\Gamma)),
 (\theta \; {\rm form} \; (\Gamma)) \in {\rm Form}(\Sigma,\Pi)$
\begin{itemize} 
\item[(ref)]  $$\infer{ \phi \sequentunder{\Gamma} \phi}{}.$$
\item[(cut)] $$\infer{ \phi \sequentunder{\Gamma} \psi}{\phi \sequentunder{\Gamma} \theta &\theta \sequentunder{\Gamma} \psi}$$
\item[(conj)] 
$$\infer{ \theta \land \psi \sequentunder{\Gamma} \theta}{}
\quad
\infer{\theta \land \psi \sequentunder{\Gamma} \psi}{}
\quad
\infer{ \phi \sequentunder{\Gamma} \theta \land \psi}{\phi \sequentunder{\Gamma} \theta &\phi \sequentunder{\Gamma} \psi}
\qquad
\infer{ \phi \sequentunder{\Gamma} \top}{}
$$
\item[(disj)] 
$$\infer{\theta \sequentunder{\Gamma} \theta \lor \psi }{}
\quad
\infer{\psi \sequentunder{\Gamma}\theta \lor \psi }{}
\quad
\infer{ \theta \lor \psi  \sequentunder{\Gamma} \phi}{\theta \sequentunder{\Gamma}  \phi &
\psi  \sequentunder{\Gamma} \phi}
\qquad
\infer{ \bot \sequentunder{\Gamma} \phi}{}$$
\item[(impl)] 
$$ \infer{ \theta   \sequentunder{\Gamma} \psi \rightarrow \phi}{\theta \land \psi  \sequentunder{\Gamma} \phi} \qquad 
\infer{\theta \land \psi  \sequentunder{\Gamma} \phi}{ \theta   \sequentunder{\Gamma} \psi \rightarrow \phi}
$$
\end{itemize}

Write for the canonical projection
\begin{equation} \label{canonpr}
{\bf p}_\Gamma(x:A) = (\langle \Gamma,x:A\rangle, \Gamma, {\rm OV}(\Gamma)).
\end{equation}

{\em Quantificational rules.} For $(\phi\; {\rm form} \; (\Gamma)), 
(\psi\; {\rm form} \; (\Gamma, x:A)) \in {\rm Form}(\Sigma,\Pi)$:

\begin{itemize}
\item[(univ)] $$\infer{\phi \sequentunder{\Gamma} (\forall x:A)\psi}{ \phi\{{\bf p}_\Gamma(x:A)\} \sequentunder{\Gamma,x:A}\psi} \qquad
\infer{\phi\{{\bf p}_\Gamma(x:A)\} \sequentunder{\Gamma,x:A}\psi}{\phi \sequentunder{\Gamma} (\forall x:A)\psi}$$
\item[(exis)] $$\infer{(\exists x:A)\psi \sequentunder{\Gamma} \phi }{\psi\sequentunder{\Gamma,x:A}\phi\{{\bf p}_\Gamma(x:A)\}} \qquad
\infer{\psi  \sequentunder{\Gamma,x:A}\phi\{{\bf p}_\Gamma(x:A)\}}{(\exists x:A)\psi \sequentunder{\Gamma} \phi}$$
\end{itemize}

{\em Substitution rule.}  For $(\phi\; {\rm form} \; (\Gamma)), 
(\psi\; {\rm form} \; (\Gamma)) \in {\rm Form}(\Sigma,\Pi)$,  and 
$(\bar{a}: \Delta \to \Gamma) \in {\cal J}(\Sigma)$,

\begin{itemize}
\item[(subs)] $$\infer{\phi\{(\Delta,\Gamma,\bar{a})\} \sequentunder{\Delta} \psi\{(\Delta,\Gamma,\bar{a})\}}{\phi \sequentunder{\Gamma} \psi}.$$
\end{itemize}

\subsection{Soundness and completeness}

In this subsection we assume that the variable system of $\Sigma$ has the de Bruijn property.
We show that the Lindenbaum-Tarski algebra $T$ over $(\Sigma,\Pi)$ forms a first-order
hyperdoctrine over ${\cal F}_\Sigma$.

For $\Gamma  \in {\cal F}_\Sigma$, define 
$${\rm Pr}_{\Sigma,\Pi,T}(\Gamma)= {\rm Pr}(\Gamma) =_{\rm def} \{(\Gamma, \phi) : (\phi\; {\rm form}\; (\Gamma)) 
\in {\rm Form}(\Sigma, \Pi)\}.$$
and define its order relation $(\le_{{\rm Pr}(\Gamma)}) = (\le)$ as follows
$$(\Gamma, \phi) \le (\Gamma, \psi)  \text{ iff }
(\phi \sequentunder{\Gamma} \psi) \in {\rm Thm}(\Sigma,\Pi,T)$$
It is clear that $({\rm Pr}(\Gamma),\le)$ is a Heyting algebra.
For an arrow $(\Delta,\Gamma,\bar{a})$ from $\Delta$ to $\Gamma$ in ${\cal F}_\Sigma$, define for $(\Gamma,\phi) \in {\rm Pr}(\Gamma)$,
$${\rm Pr}((\Delta,\Gamma,\bar{a}))(\Gamma,\phi) = (\Delta, \phi\{(\Delta,\Gamma,\bar{a})\})$$
By the substitution rule this an order preserving map, and by the definition of capture avoiding substitution it is clear that it preserves $\top,\land, \bot,\lor, \rightarrow$.
We consider the axioms for quantifiers. For 
$(\langle \Gamma,x:A\rangle, \psi) \in {\rm Pr}(\langle \Gamma,x:A\rangle)$ define
$$\forall_{(\Gamma, A)}(\langle \Gamma,x:A\rangle, \psi) =_{\rm def} (\Gamma, (\forall x:A)\psi)$$
$$\exists_{(\Gamma, A)}(\langle \Gamma,x:A\rangle, \psi) =_{\rm def} (\Gamma, (\exists x:A)\psi).$$
Consider $S=(\Gamma, A)$, $Q=(\Gamma, \phi) \in {\rm Pr}(\Gamma)$ and
$R=(\langle \Gamma,x:A\rangle, \psi) \in {\rm Pr}(\langle \Gamma,x:A\rangle)$.
Then by the quantifier rules,
\begin{eqnarray*}
Q \le \forall_S(R) &\text{ iff }&  (\phi \sequentunder{\Gamma} (\forall x:A)\psi) \in {\rm Thm}(\Sigma,\Pi,T) \\
&\text{ iff }& (\phi\{{\bf p}_\Gamma(x:A)\} \sequentunder{\Gamma,x:A}\psi) \in {\rm Thm}(\Sigma,\Pi,T)
\end{eqnarray*}
and
\begin{eqnarray*}
\exists_S(R) \le Q &\text{ iff }&  ((\exists x:A)\psi \sequentunder{\Gamma} \phi) \in {\rm Thm}(\Sigma,\Pi,T) 
\\
&\text{ iff }& (\psi  \sequentunder{\Gamma,x:A}\phi\{{\bf p}_\Gamma(x:A)\}) \in {\rm Thm}(\Sigma,\Pi,T).
\end{eqnarray*}
Now, since ${\rm fresh}(\Gamma)=x$,
\begin{eqnarray*}
Q\{{\bf p}(S)\} &=& Q\{(\Gamma. {\rm S}, \Gamma,{\rm OV}(\Gamma))\} \\
&=& Q\{( \langle \Gamma, {\rm fresh}(\Gamma):A\rangle, \Gamma,{\rm OV}(\Gamma))\} \\
&=&Q\{(\langle \Gamma,x:A\rangle, \Gamma, {\rm OV}(\Gamma))\} \\
&=&\phi\{{\bf p}_\Gamma(x:A)\},
\end{eqnarray*}
conditions 3(a) and 3(b) are fulfilled.

We check condition 4. Suppose $S=(\Gamma, A)$,  and
$R=(\langle \Gamma,x:A\rangle, \psi) \in {\rm Pr}(\langle \Gamma,x:A\rangle)$.
Let  $f=(\Delta,\Gamma,\bar{a}) \in {\rm Arr}({\cal F}_\Sigma)$.
Then
\begin{eqnarray*}
\forall_{S}(R)\{f\} &= & ((\forall x:A)\psi)\{(\Delta,\Gamma,\bar{a})\}  \\
&= &
(\forall y:A[\bar{a}/\Gamma])\bigl(\psi\{(\langle\Delta,y:A[\bar{a}/\Gamma]\rangle,\langle \Gamma,x:A\rangle,(\bar{a},y))\}\bigr)
\end{eqnarray*}
On the other hand, since  
$$S\{f\} = (\Gamma, A)\{(\Delta,\Gamma,\bar{a})\} 
        = (\Delta, A[\bar{a}/\Gamma])$$
$${\rm p}_{\Delta}(S\{f\}) 
        = (\Delta.S\{f\}, \Delta, {\rm OV}(\Delta))
        = (\langle \Delta, {\rm fresh}(\Delta):S\{f\}\rangle, \Delta, {\rm OV}(\Delta))$$
and  ${\rm fresh}(\Delta) =y$, we have
\begin{eqnarray*}
{\rm v}_{\Delta,S\{f\}} &= & (S\{f\}\{{\rm p}_{\Delta}(S\{f\})\}, {\rm fresh}(\Delta)) \\
 &= &((\Delta, A[\bar{a}/\Gamma])\{(\langle \Delta, y:S\{f\}\rangle, \Delta, {\rm OV}(\Delta))\}, y) \\
 &= &((\Delta, A[\bar{a}/\Gamma][{\rm OV}(\Delta)/\Delta]), y) \\
 &= &((\Delta, A[\bar{a}/\Gamma]), y)
\end{eqnarray*}
\begin{eqnarray*}
f.S &= & \langle f\circ {\rm p}_{\Delta}(S\{f\}), {\rm v}_{S\{f\}}\rangle_{\rm S} \\
&= & \langle (\Delta,\Gamma,\bar{a}) \circ (\langle \Delta, y:S\{f\}\rangle, \Delta, {\rm OV}(\Delta)), 
((\Delta, A[\bar{a}/\Gamma]), y) \rangle \\
&= & \langle (\langle\Delta, y:S\{f\}\rangle,\Gamma,\bar{a} \circ {\rm OV}(\Delta)) , 
((\Delta, A[\bar{a}/\Gamma]), y) \rangle \\
&= & \langle (\langle\Delta, y:A[\bar{a}/\Gamma]\rangle,\Gamma,\bar{a}) , 
((\Delta, A[\bar{a}/\Gamma]), y) \rangle \\
&=&  (\langle \Delta, y:A[\bar{a}/\Gamma] \rangle, \langle \Gamma, x:A\rangle, (\bar{a},y))
\end{eqnarray*}
Hence
\begin{eqnarray*}
\forall_{S\{f\}}(R\{f.S\})&=& \forall_{S\{f\}}((\langle \Gamma,x:A\rangle, \psi)\{f.S\}) \\
&=&(\forall y:A[\bar{a}/\Gamma])\bigl(\psi\{(\langle\Delta,y:A[\bar{a}/\Gamma]\rangle,\langle \Gamma,x:A\rangle,(\bar{a},y))\}\bigr).
\end{eqnarray*}
Thus $\forall_{S}(R)\{f\}  = \forall_{S\{f\}}(R\{f.S\})$ as required. The verification of
$$\exists_{S}(R)\{f\}  = \exists_{S\{f\}}(R\{f.S\})$$
is the same as above, just replacing $\forall$ with $\exists$.

We conclude that

\begin{theorem} (Completeness: Universal model) The Lindenbaum-Tarski algebra for a theory $T$ over $(\Sigma,\Pi)$, 
$${\cal H}_{\Sigma,\Pi,T}=_{\rm def} ({\cal F}_\Sigma, {\rm Pr}_{\Sigma,\Pi,T}, \forall,\exists)$$ is a hyperdoctrine over ${\cal F}_\Sigma$. 
It has the universal model property, i.e.\ that
$$(\Gamma,\phi) \le (\Gamma, \psi) \text{ if, and only if, } (\phi \sequentunder{\Gamma} \psi) 
\in {\rm Thm}(\Sigma,\Pi,T).$$
$\qed$
\end{theorem}

\medskip 
Let $T$ be a dependent first-order theory over $(\Sigma,\Pi)$. Let
${\cal H}$ be the hyperdoctrine  ${\cal H}_{\Sigma,\Pi,T} =({\cal F}_\Sigma, {\rm Pr}_{\Sigma,\Pi,T}, \forall,\exists)$.
A {\em dependent first-order model of $T$} consists of a $\Sigma$-structure
given by a cwf morphism $F: {\cal F}_{\Sigma} \to {\cal C}$ and a hyperdoctrine
${\cal D} = ({\cal C}, {\rm Pr}^{\cal D},\forall',\exists')$ with an $F$-based hyperdoctrine morphism $G: {\cal H} \to {\cal D}$.  
This means that $G$ is a natural transformation $G: {\rm Pr}^{\cal H} \to {\rm Pr}^{\cal D} \circ F$
 such that for $R \in {\rm Pr}^{\cal H}(\Gamma.S)$,
 \begin{enumerate}
\item $G_\Gamma(\forall_S(R))=\forall'_{\sigma_\Gamma(S)}(G_{\Gamma.S}(R))$,
\item $G_\Gamma(\exists_S(R))=\exists'_{\sigma_\Gamma(S)}(G_{\Gamma.S}(R))$.
\end{enumerate}
That $G$ is natural 
\begin{equation} 
\bfig\square<1000,600>[{\rm Pr}^{\cal H}(\Gamma)`{\rm Pr}^{\cal D}(F\Gamma)`{\rm Pr}^{\cal H}(\Delta)`{\rm Pr}^{\cal D}(F\Delta);G_\Gamma`{\rm Pr}^{\cal H}(f)`{\rm Pr}^{\cal D}(Ff)`G_\Delta]\efig
\end{equation}
amounts to the following condition on substitution: for 
$f=(\Delta,\Gamma,\bar{t})$,
$$G_\Delta(\Delta, \phi\{(\Delta,\Gamma,\bar{t})\}) =
{\rm Pr}^{\cal D}(F(\Delta,\Gamma,\bar{t}))(G_\Gamma(\Gamma,\phi)) 
=G_\Gamma(\Gamma,\phi)\{F(\Delta,\Gamma,\bar{t})\}.$$

\begin{example} {\em Let $T$ and $T'$ be dependent first-order theories over $(\Sigma,\Pi)$,
with $T \subseteq T'$. Then there is a hyperdoctrine morphism $G: {\cal H}_{\Sigma,\Pi,T} \to  
{\cal H}_{\Sigma,\Pi,T'}$
given by letting $F$ be the identity cwf morphism on ${\cal F}_{\Sigma}$,  and defining
$G_\Gamma(\Gamma,\phi)= (\Gamma, \phi)$. It is order preserving,  since
$$(\phi \sequentunder{\Gamma} \psi) \in {\rm Thm}(\Sigma,\Pi,T)
\text{ implies }
(\phi \sequentunder{\Gamma} \psi) \in {\rm Thm}(\Sigma,\Pi,T').$$
}
\end{example}

\medskip
\begin{theorem} \label{determbyrel}
Let $T$ be a dependent first-order theory over $(\Sigma,\Pi)$, where the signatures are regular.
Let $G$ and $G'$ be two $F$-based hyperdoctrine morphisms ${\cal H}_{\Sigma,\Pi,T} \to {\cal D}$.
If for every $(\Gamma, R) \in \Pi$,
$$G_{\Gamma}(\Gamma, R(\check{\Gamma}))  = G'_{\Gamma}(\Gamma, R(\check{\Gamma}))$$
holds in ${\rm Pr}^{\cal D}(F(\Gamma))$,  then $G=G'$.
\end{theorem}
{\flushleft \bf Proof.} Write ${\rm Pr}= {\rm Pr}^{{\cal H}_{\Sigma,\Pi,T}}$. It suffices to show that  for each $\Delta \in {\cal F}_\Sigma$ and
each $P \in {\rm Pr}(\Delta)$,
$$G_{\Delta}(P)  = G'_{\Delta}(P).$$
We prove this by induction on derivations in ${\rm Form}(\Sigma,\Pi)$. Let $(\Gamma, R) \in \Pi$ and 
$(\bar{a}: \Delta \to \Gamma) \in {\cal J}(\Sigma)$.
Then by naturality of $G$ and $G'$ and the assumption
\begin{eqnarray*}
G_{\Delta}(\Delta, R(\bar{a})) &= &G_{\Delta}((\Gamma, R(\check{\Gamma})\{(\Delta, \Gamma, \bar{a})\}\bigr) \\
&= &  G_{\Gamma}(\Gamma, R(\check{\Gamma})\{ F(\Delta, \Gamma, \bar{a})\} \\
&= &  G'_{\Gamma}(\Gamma, R(\check{\Gamma})\{ F(\Delta, \Gamma, \bar{a})\} \\
&= &G'_{\Delta}((\Gamma, R(\check{\Gamma})\{(\Delta, \Gamma, \bar{a})\}\bigr) \\
&= &G'_{\Delta}(\Delta, R(\bar{a}))
\end{eqnarray*}
Moreover since $G_{\Delta}$ and $G'_{\Delta}$ are Heyting algebra morphisms, the cases for propostional
operations are immediate induction steps.
Suppose $R \in {\rm Pr}(\Gamma.S)$. Then using the inductive hypothesis,
\begin{eqnarray*}
G_{\Gamma}(\forall_S(R)) &=& \forall_{\sigma_\Gamma(S)}(G_{\Gamma.S}(R)) \\
&=& \forall_{\sigma_\Gamma(S)}(G'_{\Gamma.S}(R))  \\
&=& G'_{\Gamma}(\forall_S(R)).
\end{eqnarray*}
The case for $\exists$ is proved similarly. $\qed$

\begin{theorem} Let $\Sigma$ be a regular signature and $F: {\cal F}_{\Sigma} \to {\cal C}$ a cwf morphism. Suppose that $\Pi$ is a regular predicate signature 
over $\Sigma$. Let ${\cal D}$ be a hyperdoctrine based on ${\cal C}$, and suppose 
that for each $R=(\Gamma_R,P_R) \in \Pi$, there is $R^* \in {\rm Pr}_{\cal C}(\Gamma_R)$.
Then there is a unique $F$-based hyperdoctrine morphism 
$G: {\cal H}_{\Sigma,\Pi,\emptyset} \to {\cal D}$ with
$$G_{\Gamma_R}(P_R({\rm OV}(\Gamma_R))) = R^*$$
for each $R \in \Pi$.
\end{theorem}
{\flushleft \bf Proof.} The uniqueness follows by Theorem \ref{determbyrel}.
Suppose ${\cal D}=({\cal C},{\rm Pr}^{\cal D},\forall',\exists')$ and write 
$F= (F,\sigma,\theta)$.
Define $G$ by induction on the height of derivations in ${\rm Form}(\Sigma,\Pi)$.
We define
for each $R \in \Pi$  and $(\bar{a}: \Delta \to \Gamma_R) \in {\cal J}(\Sigma)$ define
$$G_{\Delta}(\Delta, R(\bar{a})) = R^*\{F(\Delta, \Gamma_R, \bar{a})\}$$
Furthermore let 
$$G_{\Delta}(\Delta, \top) =\top' \qquad  G_{\Delta}(\Delta, \bot) =\bot'.$$
For the propositional operations $\bigcirc = \land, \lor, \rightarrow$,
$$G_{\Delta}(\Delta, \phi \bigcirc \psi) = G_{\Delta}(\Delta, \phi) \bigcirc' 
G_{\Delta}(\Delta, \psi).$$
For the quantifiers $Q= \forall,\exists$, 
$$G_{\Delta}(\Delta, (Qx:A)\phi) = Q'_{\sigma_{\Delta}(A)}(G_{\Delta,x:A}(\langle \Delta,x:A\rangle, \phi)).$$
Next, we need show that $G$ is a natural transformation, i.e.\
$$G_\Delta(\Delta,\phi)\{F(\Theta,\Delta,\bar{t})\} =G_\Theta(\Theta, \phi\{(\Theta,\Delta,\bar{t})\})$$
for all $(\bar{t}: \Theta \to \Delta) \in {\cal J}(\Sigma)$.
This is done by induction on the formula $\phi$:

As for the atomic predicate case: 
\begin{eqnarray*}
G_{\Delta}(\Delta, R(\bar{a}))\{F(\Theta,\Delta,\bar{t})\} &=&
R^*\{F(\Delta, \Gamma_R, \bar{a})\}\{F(\Theta,\Delta,\bar{t})\} \\
&=& R^*\{F(\Theta, \Gamma_R, \bar{a} \circ \bar{t})\} \\
&=& G_{\Theta}(\Theta, R(\bar{a}[\bar{t}/\Delta])) \\
&=& G_{\Theta}(\Theta, R(\bar{a})\{(\Theta, \Delta, \bar{t})\}).
\end{eqnarray*}

As for the constants $K= \bot,\top$: 
\begin{eqnarray*}
G_{\Delta}(\Delta, K)\{F(\Theta,\Delta,\bar{t})\} &=& K'\{F(\Theta,\Delta,\bar{t})\}  \\
&=& K' \\
&=&  G_{\Theta}(\Theta, K)\\
&=&  G_{\Theta}(\Theta, K\{(\Theta,\Delta,\bar{t})\})\\
\end{eqnarray*}

As for the connectives $\bigcirc = \land, \lor, \rightarrow$,
\begin{eqnarray*}
G_{\Delta}(\Delta, \phi \bigcirc \psi)\{F(\Theta,\Delta,\bar{t})\} &=& 
(G_{\Delta}(\Delta, \phi)  \bigcirc'  G_{\Delta}(\Delta,\psi))\{F(\Theta,\Delta,\bar{t})\}\\
&=& 
(G_{\Delta}(\Delta, \phi)\{F(\Theta,\Delta,\bar{t})\})  \bigcirc'  (G_{\Delta}(\Delta,\psi))\{F(\Theta,\Delta,\bar{t})\}\\
&=& G_\Theta(\Theta, \phi\{(\Theta,\Delta,\bar{t})\}) \bigcirc'
G_\Theta(\Theta, \psi\{(\Theta,\Delta,\bar{t})\}) \\
&=& G_\Theta(\Theta, \phi\{(\Theta,\Delta,\bar{t})\} \bigcirc  \psi\{(\Theta,\Delta,\bar{t})\}) \\
&=& G_\Theta(\Theta, (\phi \bigcirc \psi)\{(\Theta,\Delta,\bar{t})\}).
\end{eqnarray*}

As for the quantifiers $Q= \forall,\exists$, we have the following equations, where step 1 
is by definition of $G$,  step 2 by the Beck-Chevalley conditions (\ref{BCC}) and step 3 is by naturality of $\sigma$,
\begin{eqnarray*}
\lefteqn{G_{\Delta}(\Delta, (Qx:A)\phi)\{F(\Theta,\Delta,\bar{t})\}} \\
 &=& 
Q'_{\sigma_{\Delta}(A)}(G_{\Delta,x:A}(\langle \Delta,x:A\rangle, \phi))\{F(\Theta,\Delta,\bar{t})\}
\\
&=& Q'_{\sigma_{\Delta}(A)\{F(\Theta,\Delta,\bar{t})\}}(G_{\Delta,x:A}(\langle \Delta,x:A\rangle, \phi)\{F(\Theta,\Delta,\bar{t}).\sigma_{\Delta}(A)\})\\
&=& Q'_{\sigma_{\Theta}(A[\bar{t}/\Delta])}(G_{\Delta,x:A}(\langle \Delta,x:A\rangle, \phi)\{F(\Theta,\Delta,\bar{t}).\sigma_{\Delta}(A)\})\\
\end{eqnarray*}
On the other hand we have the equations below, where step 1 is by definition of substitution,
step 2 by the definition of $G$ and step 3 the inductive hypothesis,
\begin{eqnarray*}
\lefteqn{G_{\Delta}(\Delta, ((Qx:A)\phi)\{(\Theta,\Delta,\bar{t})\})} \\
&=&
G_{\Theta}(\Theta, (Qy:A[\bar{t}/\Delta])\phi\{(\langle \Theta, y:A[\bar{t}/\Delta]\rangle,\langle\Delta,x:A\rangle,(\bar{t},y))\}) \\
&=& Q'_{\sigma_\Theta(A[\bar{t}/\Delta])}(G_{\Theta,y:A[\bar{t}/\Delta]}(\langle \Theta, y:A[\bar{t}/\Delta] \rangle,
\phi\{(\langle \Theta, y:A[\bar{t}/\Delta]\rangle,\langle\Delta,x:A\rangle,(\bar{t},y))\})) \\
&=& Q'_{\sigma_\Theta(A[\bar{t}/\Delta])}(G_{\Delta, x:A}(\langle \Delta, x:A \rangle,\phi)
\{F(\langle \Theta, y:A[\bar{t}/\Delta]\rangle,\langle\Delta,x:A\rangle,(\bar{t},y))\}) \\
\end{eqnarray*}
where $y={\rm fresh}(\Theta)$. Now it only remains to prove:
\begin{eqnarray*}
F(\Theta,\Delta,\bar{t}).\sigma_{\Delta}(A)
& = & F(\langle \Theta, y:A[\bar{t}/\Delta]\rangle,\langle\Delta,x:A\rangle,(\bar{t},y))
\end{eqnarray*}
But this follows from Lemma  \ref{ffa_lm}. $\qed$

\begin{theorem} (Soundness) Let $\Sigma$ be a regular signature and $F: {\cal F}_{\Sigma} \to {\cal C}$ a cwf morphism. Suppose that $\Pi$ is a regular predicate signature 
over $\Sigma$. Let ${\cal D}$ be a hyperdoctrine based on ${\cal C}$.
Suppose that  $G: {\cal H}_{\Sigma,\Pi,\emptyset} \to {\cal D}$ is an
$F$-based hyperdoctrine morphism. Let $T$ be a theory over $\Sigma,\Pi$ and assume 
that for every sequent $(\phi \sequentunder{\Gamma} \psi) \in T$, 
$$G_\Gamma(\Gamma,\phi) \le G_\Gamma(\Gamma,\psi).$$
Then $G: {\cal H}_{\Sigma,\Pi,T} \to {\cal D}$ is also a model of $T$.
\end{theorem}
{\flushleft \bf Proof.} 
We prove by induction on proofs that
\begin{equation} \label{sound}
(\phi \sequentunder{\Gamma} \psi) \in {\rm Thm}(\Sigma,\Pi,T) \text{ implies }
G_\Gamma(\Gamma,\phi) \le G_\Gamma(\Gamma,\psi).
\end{equation}
The basic case is the assumption. The (ref) rule is trivial, and (cut) rule follows from the inductive hypotheses.
The propositional rules follows from the inductive hypothesis noting that $G_\Gamma$ preserves the propositional
operations.

Quantifier rules.

(univ, left): Suppose the last rule applied is 
$$\infer{\phi \sequentunder{\Gamma} (\forall x:A)\psi}{ \phi\{{\bf p}_\Gamma(x:A)\} \sequentunder{\Gamma,x:A}\psi}$$
By inductive hypothesis
$$G_{\Gamma,x:A}(\langle \Gamma,x:A\rangle,\phi\{{\bf p}_\Gamma(x:A)\} ) \le G_{\Gamma,x:A}(\langle \Gamma,x:A\rangle,\psi).$$
By naturality of $G$ and since $F$ preserves projections (\ref{ppres}):
\begin{eqnarray*}
G_{\Gamma,x:A}(\langle \Gamma,x:A\rangle,\phi\{{\bf p}_\Gamma(x:A)\}) &=&
(G_{\Gamma}(\Gamma,\phi)) \{F({\bf p}_\Gamma(x:A))\}\\
&=&
(G_{\Gamma}(\Gamma,\phi))\{{\rm p}_{F(\Gamma)}(\sigma_{\Gamma}(x:A))\}
\end{eqnarray*}
Thus we have
$$(G_{\Gamma}(\Gamma,\phi))\{{\rm p}_{F(\Gamma)}(\sigma_{\Gamma}(x:A))\} \le 
G_{\Gamma,x:A}(\langle \Gamma,x:A\rangle,\psi).$$
By property 3(a) of the hyperdoctrine,  and using that $G$ is morphism of hyperdoctrines
(property 1):
$$G_{\Gamma}(\Gamma,\phi) \le \forall_{\sigma_{\Gamma}(x:A)}(G_{\Gamma,x:A}(\langle \Gamma,x:A\rangle,\psi))= G_\Gamma(\Gamma, (\forall x:A)\psi)$$

(univ,right): similar.

(exis, left): analogous using property 2 of $G$. 

(exis, right): analogous using property 2 of $G$. 

\begin{itemize}
\item[(univ)] $$\infer{\phi \sequentunder{\Gamma} (\forall x:A)\psi}{ \phi\{{\bf p}_\Gamma(x:A)\} \sequentunder{\Gamma,x:A}\psi} \qquad
\infer{\phi\{{\bf p}_\Gamma(x:A)\} \sequentunder{\Gamma,x:A}\psi}{\phi \sequentunder{\Gamma} (\forall x:A)\psi}$$
\item[(exis)] $$\infer{(\exists x:A)\psi \sequentunder{\Gamma} \phi }{\psi\sequentunder{\Gamma,x:A}\phi\{{\bf p}_\Gamma(x:A)\}} \qquad
\infer{\psi  \sequentunder{\Gamma,x:A}\phi\{{\bf p}_\Gamma(x:A)\}}{(\exists x:A)\psi \sequentunder{\Gamma} \phi}$$
\end{itemize}

Subsititution rule:  Suppose the last rule used is (subs). I.e. for $(\bar{a}: \Delta \to \Gamma) \in {\cal J}(\Sigma)$,
\begin{itemize}
\item[(subs)] $$\infer{\phi\{(\Delta,\Gamma,\bar{a})\} \sequentunder{\Delta} \psi\{(\Delta,\Gamma,\bar{a})\}}{\phi \sequentunder{\Gamma} \psi}.$$
\end{itemize}
By inductive hypothesis we have
$$G_\Gamma(\Gamma,\phi) \le G_\Gamma(\Gamma,\psi).$$
Thus applying the substitution $F(\Delta,\Gamma,\bar{a}))$,
$$(G_\Gamma(\Gamma,\phi))\{F(\Delta,\Gamma,\bar{a}))\} \le (G_\Gamma(\Gamma,\psi))\{F(\Delta,\Gamma,\bar{a}))\}$$
The naturality of $G$ gives
$$G_\Delta(\Delta,\phi\{(\Delta,\Gamma,\bar{a})\}) \le G_\Delta(\Delta,\psi\{(\Delta,\Gamma,\bar{a})\})$$
as required.
 $\qed$

\subsection{A variant with syntactic substitution} \label{dfolstar}

We briefly consider a version of DFOL with ordinary syntactic substitution
and unrestricted variable systems. This version will be called DFOL* and
is close to the system of (Belo 2007).  In that paper, raw formulas are defined from
 the symbols from our signatures $(\Sigma,\Pi)$. The notion of simultaneous
syntactic substitution $\phi[\bar{t}/\bar{x}]$ into formulas is defined. To make this operation totally defined we use $\alpha$-equivalence between formulas. Taking equivalence classes of formulas allows substitution to be defined on representatives where no illegitmate binding of free variables in $\bar{t}$ by bound variables of $\phi$ at the location of substitution takes place; see (Prop.\ 3 and 4, Belo 2007). We can now reconsider the inductive definition
of ${\rm Form}(\Sigma,\Pi)$  in (F1) -- (F4) as taking the corresponding $\alpha$-equivalence classes of formulas instead of only formulas.  We call the resulting set ${\rm Form}^*(\Sigma,\Pi)$. (Note that the previous example $(\forall x:A)R(x) \; {\rm form}\; (x:A)$  (of Remark \ref{boundfreedist}) is actually in this set, since it is $\alpha$-equivalent to $(\forall y:A)R(y) \; {\rm form}\; (x:A)$.)

\medskip
The rules (ref), (cut) and rules for the propositional connectives are the same as
for DFOL. The substitution rule is now

{\flushleft \em Syntactic substitution rule.}  For $(\phi\; {\rm form} \; (\Gamma)), 
(\psi\; {\rm form} \; (\Gamma)) \in {\rm Form}^*(\Sigma,\Pi)$,  and any
$(\bar{a}: \Delta \to \Gamma) \in {\cal J}(\Sigma)$
\begin{itemize}
\item[(sub*)] $$\infer{\phi[\bar{a}/\Gamma] \sequentunder{\Delta} \psi[\bar{a}/\Gamma]}{\phi \sequentunder{\Gamma} \psi}.$$
\end{itemize}

The quantifier rules can now be simplified to align with those of (Johnstone 2002).

{\em Quantificational rules.} For $(\phi\; {\rm form} \; (\Gamma)), 
(\psi\; {\rm form} \; (\Gamma, x:A)) \in {\rm Form}^*(\Sigma,\Pi)$:

\begin{itemize}
\item[(univ*)] $$\infer{\phi \sequentunder{\Gamma} (\forall x:A)\psi}{ \phi \sequentunder{\Gamma,x:A}\psi} \qquad
\infer{\phi \sequentunder{\Gamma,x:A}\psi}{\phi \sequentunder{\Gamma} (\forall x:A)\psi}$$
\item[(exis*)] $$\infer{(\exists x:A)\psi \sequentunder{\Gamma} \phi }{\psi\sequentunder{\Gamma,x:A}\phi} \qquad
\infer{\psi  \sequentunder{\Gamma,x:A}\phi}{(\exists x:A)\psi \sequentunder{\Gamma} \phi}$$
\end{itemize}

For a theory $T$ over $(\Sigma,\Pi)$ let ${\rm Thm}^*(\Sigma,\Pi,T)$ denote the smallest set of sequents derivable from $T$ using the 
above mentioned rules.

\medskip
 A signature $(\Sigma,\Pi)$ is called {\em normal} if it is regular and has the de Bruijn property. For a normal signature $\Sigma$ with the variable system ${\cal V}$ we let $\Sigma^{(\infty)}$ denote the same signature with unrestricted variable system ${\cal V}^{(\infty)}$. The formulas of
 ${\rm Form}(\Sigma^{(\infty)},\Pi)$ thus have no restrictions on the bound variable names possible.
  
\begin{lemma} \label{subvssub} Let $(\Sigma,\Pi)$ be a normal signature.
For every $(\phi\; {\rm form} \; (\Gamma)) \in {\rm Form}(\Sigma,\Pi)$, and
$(\bar{t}: \Delta \to \Gamma) \in {\cal J}(\Sigma)$ on standard form,
$$(\phi\{(\Delta,\Gamma,\bar{t})\}\; {\rm form} \; (\Gamma)) \equiv_\alpha 
(\phi[\bar{t}/\Gamma]\; {\rm form} \; (\Gamma))$$
in ${\rm Form}(\Sigma^{(\infty)},\Pi)$.
\end{lemma}
{\flushleft \bf Proof.} The atomic base case is clear, as well as the inductive cases for propositional connectives. Suppose $\phi=(Qx:A)\theta$, ($Q=\forall$ or $\exists$) where 
 $(\theta\; {\rm form}\; (\Gamma,x:A)) \in  {\rm Form}(\Sigma,\Pi)$. Then
 $$((Qx:A)\theta)\{(\Delta,\Gamma,\bar{a})\} =
   (Qy:A[\bar{a}/\Gamma])\bigl(\theta\{(\langle\Delta,y:A[\bar{a}/\Gamma]\rangle,\langle \Gamma,x:A\rangle,(\bar{a},y))\}\bigr)$$
   where $y= {\rm fresh}(\Delta)$.
Now by inductive hypothesis
\begin{eqnarray*}
\lefteqn{(Qy:A[\bar{a}/\Gamma])\bigl(\theta\{(\langle\Delta,y:A[\bar{a}/\Gamma]\rangle,\langle \Gamma,x:A\rangle,(\bar{a},y))\}\bigr)} \\
&\equiv_{\alpha}&
(Qy:A[\bar{a}/\Gamma])\bigl(\theta[\bar{a},y/ \Gamma,x:A]\bigr) \\
&\equiv_{\alpha}& ((Qx:A)\theta)[\bar{a}/\Gamma].
\end{eqnarray*}
which establishes the result.
$\qed$

\begin{theorem} Let $(\Sigma,\Pi)$ be a normal signature, and let $T$ be a theory on standard form over $(\Sigma,\Pi)$.
If $(\phi \sequentunder{\Gamma} \psi) \in {\rm Thm}(\Sigma,\Pi,T)$, then
$(\phi \sequentunder{\Gamma} \psi) \in {\rm Thm}^*(\Sigma^{(\infty)},\Pi,T)$.
\end{theorem}
{\flushleft \bf Proof.} Induction on derivations. The base cases are the axioms: the logical axioms (ref), axioms for $\top$ and $\bot$, and the non-logical axioms of $T$. These are clear. The inductive cases for propositional connectives follows by the inductive hypotheses,
as well as for the (cut) rule.  Consider the quantifier rules
\begin{itemize}
\item[(univ)] $$\infer{\phi \sequentunder{\Gamma} (\forall x:A)\psi}{ \phi\{{\bf p}_\Gamma(x:A)\} \sequentunder{\Gamma,x:A}\psi} \qquad
\infer{\phi\{{\bf p}_\Gamma(x:A)\} \sequentunder{\Gamma,x:A}\psi}{\phi \sequentunder{\Gamma} (\forall x:A)\psi}$$
\item[(exis)] $$\infer{(\exists x:A)\psi \sequentunder{\Gamma} \phi }{\psi\sequentunder{\Gamma,x:A}\phi\{{\bf p}_\Gamma(x:A)\}} \qquad
\infer{\psi  \sequentunder{\Gamma,x:A}\phi\{{\bf p}_\Gamma(x:A)\}}{(\exists x:A)\psi \sequentunder{\Gamma} \phi}$$
\end{itemize}
Since by Lemma \ref{subvssub}
$$(\phi\{{\bf p}_\Gamma(x:A)\} \; {\rm form}\; (\Gamma,x:A)) \equiv_\alpha 
(\phi \; {\rm form}\; (\Gamma,x:A))$$
we can use the corresponding rules (univ*) and (exis*) and the inductive hypothesis
to prove the quantifier cases. Similarly, using this lemma we can apply (subs*) in the case
of the (sub) rule. $\qed$

Let $\sigma: {\mathbb N} \to V$ be a fresh variable sequence for $\Sigma$. For every formula $(\phi \; {\rm form}\; (\Gamma)) \in {\rm Form}^*(\Sigma,\Pi)$ we introduce
a standardized formula $(\phi^\sigma \; {\rm form}\; (\Gamma^\sigma)) \in {\rm Form}^*(\Sigma,\Pi)$ by induction on formulas:
\begin{itemize}
\item If $(\phi \; {\rm form}\; (\Gamma)) \in {\rm Form}^*(\Sigma,\Pi)$ is atomic, then
 $(\phi[\sigma(\Gamma)/\Gamma] \; {\rm form}\; (\Gamma^\sigma)) \in {\rm Form}^*(\Sigma,\Pi),$
 so let $$\phi^\sigma=\phi[\sigma(\Gamma)/\Gamma].$$
\item Suppose that  $(\phi \; {\rm form}\; (\Gamma)), (\psi \; {\rm form}\; (\Gamma)) \in {\rm Form}^*(\Sigma,\Pi)$. Then we have $(\phi^\sigma \; {\rm form}\; (\Gamma^\sigma))$, $(\psi^\sigma \; {\rm form}\; (\Gamma^\sigma)) \in {\rm Form}^*(\Sigma,\Pi)$.
Let $$(\phi \bigcirc \psi)^\sigma = \phi^\sigma \bigcirc \psi^\sigma.$$ Hence
$((\phi \bigcirc \psi)^\sigma\; {\rm form}\; (\Gamma^\sigma)) \in {\rm Form}^*(\Sigma,\Pi)$.
\item Suppose that $(\phi \; {\rm form}\; (\Gamma,x:A)) \in {\rm Form}^*(\Sigma,\Pi)$. By
inductive hypothesis, $$(\phi^\sigma\; {\rm form}\; (\Gamma^\sigma,
\sigma(|\Gamma|):A[\sigma(\Gamma)/\Gamma])) \in {\rm Form}^*(\Sigma,\Pi).$$
Thus
$$\bigl(Q\, \sigma(|\Gamma|):A[\sigma(\Gamma)/\Gamma]\bigr)\phi^\sigma \; {\rm form}\; (\sigma(\Gamma))) \in {\rm Form}^*(\Sigma,\Pi).$$
Let $$((Q x:A)\phi)^\sigma = \bigl(Q\, \sigma(|\Gamma|):A[\sigma(\Gamma)/\Gamma]\bigr)\phi^\sigma.$$
\end{itemize}

\begin{lemma} Suppose that $(\phi \; {\rm form}\; (\Gamma)) \in {\rm Form}^*(\Sigma,\Pi)$,
and that $\sigma: {\mathbb N} \to V$ is a fresh variable sequence for $\Sigma$. Then:
$$\phi[\sigma(\Gamma)/\Gamma]\; {\rm form}\; (\Gamma^\sigma) 
\equiv_\alpha \phi^{\sigma}\; {\rm form}\; (\Gamma^\sigma)
$$
\end{lemma}
{\flushleft \bf Proof.} By induction on formulas. For the atomic case, and the propositional connectives case this is straightforward. Consider the quantifier case: Let $\phi=(Qx:A)\theta$, where $(\theta\; {\rm form}\; (\Gamma,x:A)) \in {\rm Form}^*(\Sigma,\Pi)$. We have
using the inductive hypothesis,
\begin{eqnarray*}
((Q x:A)\theta)^\sigma &=& \bigl(Q\, \sigma(|\Gamma|):A[\sigma(\Gamma)/\Gamma]\bigr)
\theta^\sigma \\
&\equiv_\alpha & \bigl(Q\, \sigma(|\Gamma|):A[\sigma(\Gamma)/\Gamma]\bigr)
\theta[\sigma(\Gamma,x:A)/\langle \Gamma, x:A\rangle] \\
&= & \bigl(Q\, \sigma(|\Gamma|):A[\sigma(\Gamma)/\Gamma]\bigr)
\theta[(\sigma(\Gamma),\sigma(|\Gamma|))/\langle \Gamma, x:A\rangle] \\
&\equiv_\alpha & ((Q x:A)\theta)[\sigma(\Gamma)/\Gamma]
\end{eqnarray*}
This finishes the induction. $\qed$

\medskip
If $T$ is a theory with respect to $(\Sigma,\Pi)$, then let 
$$T^\sigma = \{ (\phi^\sigma \sequent \psi^\sigma\;  (\Gamma^\sigma)): 
(\phi \sequent \psi\;  (\Gamma)) \in T\},$$
which is a theory with respect to the same signature.

\medskip
For every signature $\Sigma$, there is a canonical fresh variable sequence given
by $\sigma_\Sigma = \sigma$, 
$$\sigma(n) = {\rm fr}_\Sigma(\{\sigma(0),\ldots,\sigma(n-1)\}).$$
Then for a context $\Gamma=x_1:A_1,\ldots,x_n:A_n$, the map $\sigma(\Gamma): \Gamma^\sigma \to \Gamma$
is given by 
$$\sigma(\Gamma) = (\sigma(0),\ldots,\sigma(n-1)),$$
and its inverse is just
$$\sigma^{-1}(\Gamma) = (x_1,\ldots,x_n).$$

\begin{lemma} \label{sigmasub00}
 Let $\Sigma$ be a normal signature.
If ${\cal U} \in {\cal J}(\Sigma^{(\infty)})$ is a judgement where the context involved is on standard form with respect to $\Sigma$, then
${\cal U} \in {\cal J}(\Sigma)$.
\end{lemma}
{\flushleft \bf Proof.} By induction on derivations. The case (R1) is direct. 

As for case (R2): Suppose $(\Gamma, x:A\; {\rm context}) \in {\cal J}(\Sigma^{(\infty)})$ and that
$\Gamma, x:A$ is on standard form. Then $x={\rm fresh}(\Gamma)$, and by inductive hypothesis $(\Gamma {\rm context}) \in {\cal J}(\Sigma)$ and
$(A\; {\rm type}\; (\Gamma)) \in {\cal J}(\Sigma)$. Hence by (R2), $(\Gamma, x:A\; {\rm context}) \in {\cal J}(\Sigma)$. 

Case (R3): by inductive hypothesis.

Case (R4): Suppose that $(S(\bar{a}_{\bar{i}})\; {\rm type}\; (\Delta)) \in {\cal J}(\Sigma^{(\infty)})$ is derived by (R4), where $(\Gamma,S, \bar{i}) \in \Sigma^{(\infty)}$, and
by shorter derivations $(\bar{a}: \Delta \to \Gamma) \in {\cal J}(\Sigma^{(\infty)})$. Suppose
moreover that $\Delta$ is on standard form. Now
since $\Sigma^{(\infty)}$ and $\Sigma$ have the same declarations,  
$(\Gamma,S, \bar{i}) \in \Sigma$. Thus $\Gamma$ is on standard form, and since $\Delta$ is
on standard form, the judgements in $\bar{a}: \Delta \to \Gamma$ have contexts
on standard form. It follows by inductive hypothesis,
$(\bar{a}: \Delta \to \Gamma) \in {\cal J}(\Sigma)$. Hence by (R4),  
$(S(\bar{a}_{\bar{i}})\; {\rm type}\; (\Delta)) \in {\cal J}(\Sigma)$.

Case (R5): similar to (R4). $\qed$

\begin{lemma} \label{sigmasub0}
 Let $\Sigma$ be a normal signature, and let $\sigma=\sigma_\Sigma$.
If $(\bar{t}: \Delta \to \Gamma) \in {\cal J}(\Sigma^{(\infty)})$, then
$$(\sigma^{-1}(\Gamma) \circ \bar{t} \circ \sigma(\Delta): \Delta^\sigma \to \Gamma^\sigma) \in {\cal J}(\Sigma).$$
If $|\Gamma|=m$ and  $|\Delta| =n$,
$$\sigma^{-1}(\Gamma) \circ \bar{t} \circ \sigma(\Delta) = (t_1[\bar{\sigma}(n)/\Delta],\ldots,
t_m[\bar{\sigma}(n)/\Delta]),$$
where $\bar{\sigma}(n) = \sigma(0),\ldots,\sigma(n-1)$.
\end{lemma}
{\flushleft \bf Proof.} We have by composition
$$(\sigma^{-1}(\Gamma) \circ \bar{t} \circ \sigma(\Delta): \Delta^\sigma \to \Gamma^\sigma) \in {\cal J}(\Sigma^{(\infty)}).$$
But all the judgements in this compound judgement have contexts on standard form with respect to $\Sigma$,
so by Lemma \ref{sigmasub00}, 
$$(\sigma^{-1}(\Gamma) \circ \bar{t} \circ \sigma(\Delta): \Delta^\sigma \to \Gamma^\sigma) \in {\cal J}(\Sigma).$$
$\qed$

\medskip
Write 
$$(\Delta,\Gamma,\bar{t})^\sigma=_{\rm def}
(\Delta^\sigma, \Gamma^\sigma, \sigma^{-1}(\Gamma) \circ \bar{t} \circ \sigma(\Delta)).$$
Consider the canonical projection
$${\bf p}_\Gamma(x:A) = (\langle \Gamma,x:A\rangle, \Gamma, {\rm OV}(\Gamma)).$$
Thus 
\begin{equation} \label{ssubeq1}
({\bf p}_\Gamma(x:A))^\sigma = 
(\langle \Gamma^\sigma, \sigma(|\Gamma|):A[\sigma(\Gamma)/\Gamma] \rangle, \Gamma^\sigma, \bar{\sigma}(\Gamma)) =
{\bf p}_{\Gamma^\sigma}(\sigma(|\Gamma|):A[\sigma(\Gamma)/\Gamma]).
\end{equation}

\begin{lemma} \label{sigmasub}
 Let $(\Sigma,\Pi)$ be a normal signature. Suppose
that $(\phi \; {\rm form}\; (\Gamma)) \in {\rm Form}(\Sigma^{(\infty)},\Pi)$ and that
$(\bar{t}: \Delta \to \Gamma) \in {\cal J}(\Sigma^{(\infty)})$. Then for  
$\sigma=\sigma_\Sigma$:
$$(\phi^\sigma)\{(\Delta,\Gamma,\bar{t})^\sigma\} =
(\phi^\sigma)\{(\Delta^\sigma, \Gamma^\sigma, \sigma^{-1}(\Gamma) \bar{t}\sigma(\Delta))\} = (\phi[\bar{t}/\Gamma])^\sigma.$$
\end{lemma}
{\flushleft \bf Proof.} By induction on the formula $\phi$. 

Case $\phi=P$ atomic:
\begin{eqnarray*}
\lefteqn{(P^\sigma)\{(\Delta,\Gamma,\bar{t})^\sigma\}}  \\
&=& (P[\sigma(\Gamma)/\Gamma])\{(\Delta,\Gamma,\bar{t})^\sigma\} \\
&=& (P[\sigma(\Gamma)/\Gamma])[\sigma^{-1}(\Gamma) \bar{t}\sigma(\Delta)/\sigma(\Gamma)] \\
&=& (P[\bar{t}/\Gamma])[\sigma(\Delta)/\Delta] \\
&=& (P[\bar{t}/\Gamma])^\sigma.
\end{eqnarray*}

Case $\phi=\psi \bigcirc \theta$: Using the inductive hypothesis we get
\begin{eqnarray*}
\lefteqn{((\psi \bigcirc \theta)^\sigma)\{(\Delta,\Gamma,\bar{t})^\sigma\}}  \\
&=& (\psi^\sigma \bigcirc \theta^\sigma)\{(\Delta,\Gamma,\bar{t})^\sigma\} \\
&=& (\psi^\sigma)\{(\Delta,\Gamma,\bar{t})^\sigma\}
 \bigcirc (\theta^\sigma)\{(\Delta,\Gamma,\bar{t})^\sigma\} \\
&=&  (\psi[\bar{t}/\Gamma])^\sigma \bigcirc  (\theta[\bar{t}/\Gamma])^\sigma \\
&=&  (\psi[\bar{t}/\Gamma] \bigcirc  \theta[\bar{t}/\Gamma])^\sigma \\
&=&   ((\psi\bigcirc\theta)[\bar{t}/\Gamma])^\sigma
\end{eqnarray*}

Case $\phi=(Q x:A)\psi$:  We have
\begin{eqnarray*}
\lefteqn{(((Q x:A)\psi)^\sigma)\{(\Delta,\Gamma,\bar{t})^\sigma\}}  \\
&=& (\bigl(Q\, \sigma(|\Gamma|):A[\sigma(\Gamma)/\Gamma]\bigr)\psi^\sigma) 
\{(\Delta,\Gamma,\bar{t})^\sigma\}\\
&=& \bigl(Q\, y:A[\sigma(\Gamma)/\Gamma][\sigma^{-1}(\Gamma) \bar{t}\sigma(\Delta)/\sigma(\Gamma)]\bigr) \\
&& \qquad \psi^\sigma 
\{(\langle \Delta^\sigma, y:A[\sigma(\Gamma)/\Gamma][\sigma^{-1}(\Gamma) \bar{t}\sigma(\Delta)/\sigma(\Gamma)] \rangle , \\
&& \qquad \qquad \qquad \langle \Gamma^\sigma, \sigma(|\Gamma|):A[\sigma(\Gamma)/\Gamma]\rangle, 
(\sigma^{-1}(\Gamma) \bar{t}\sigma(\Delta),y))\}
\end{eqnarray*}
where $y= {\rm fresh}(\Delta^\sigma)$. The right hand side is
\begin{equation} \label{Qxpsi}
\begin{aligned}
\lefteqn{\bigl(Q\, y:A[\bar{t}/\Gamma][\sigma(\Delta)/\Delta]\bigr)} \\
& \qquad \psi^\sigma 
\{(\langle \Delta^\sigma, y:A[\bar{t}/\Gamma][\sigma(\Delta)/\Delta] \rangle, \\
& \qquad \qquad \qquad \langle \Gamma^\sigma, \sigma(|\Gamma|):A[\sigma(\Gamma)/\Gamma]\rangle, 
(\sigma^{-1}(\Gamma) \bar{t}\sigma(\Delta),y))\}\\
\end{aligned}
\end{equation}
Now consider the context map, where $z \in {\rm Fresh}(\Delta)$,
$$f=_{\rm def}(\langle \Delta, z:A[\bar{t}/\Gamma]\rangle,\langle\Gamma,x:A\rangle, (\bar{t},z)) \in {\cal J}(\Sigma^{(\infty)}).$$
Then
\begin{eqnarray*}
f^{\sigma} &=& \bigl(\langle \Delta^\sigma, \sigma(|\Delta|):A[\bar{t}/\Gamma][\sigma(\Delta)/\Delta]\rangle,\langle\Gamma^\sigma,\sigma(|\Gamma|):A[\sigma(\Gamma)/\Gamma]\rangle, \\
&& \qquad \qquad \sigma^{-1}(\langle\Gamma,x:A\rangle) \circ (\bar{t},z) \circ \sigma(\langle \Delta, z:A[\bar{t}/\Gamma]\rangle)\bigr)
\end{eqnarray*}
By Lemma \ref{sigmasub0} with $|\Gamma|=m$ and  $|\Delta| =n$,
\begin{eqnarray*}
\lefteqn{\sigma^{-1}(\langle\Gamma,x:A\rangle) \circ \bar{t} \circ \sigma(\langle \Delta, z:A[\bar{t}/\Gamma]\rangle)}  \\
&=& (t_1[\bar{\sigma}(n+1)/\langle \Delta, z:A[\bar{t}/\Gamma]\rangle],\ldots,
t_m[\bar{\sigma}(n+1)/\langle \Delta, z:A[\bar{t}/\Gamma]\rangle, z[\bar{\sigma}(n+1)/\langle \Delta, z:A[\bar{t}/\Gamma]\rangle]) \\
&=& (t_1[\bar{\sigma}(n)/\Delta],\ldots,
t_m[\bar{\sigma}(n)/\Delta]\rangle, \sigma(n))\\
&=& (\sigma^{-1}(\Gamma) \bar{t}\sigma(\Delta), y). 
\end{eqnarray*}
Here we have used that  $\sigma(n)=\sigma(|\Delta|) = {\rm fresh}(\Delta^\sigma)=y$.
Thus 
$$f^{\sigma} = \bigl(\langle \Delta^\sigma, y:A[\bar{t}/\Gamma][\sigma(\Delta)/\Delta]\rangle,\langle\Gamma^\sigma,\sigma(|\Gamma|):A[\sigma(\Gamma)/\Gamma]\rangle,
(\sigma^{-1}(\Gamma) \bar{t}\sigma(\Delta), y)\bigr),$$
which occurs  in (\ref{Qxpsi}).

The inductive hypothesis gives
$$\psi^\sigma\{f^\sigma\} = (\psi[\bar{t},z/\Gamma,x:A])^\sigma$$
Combining the above equalities we get
\begin{eqnarray*}
\lefteqn{(((Q x:A)\psi)^\sigma)\{(\Delta,\Gamma,\bar{t})^\sigma\}}  \\
&=&
\bigl(Q\, y:A[\bar{t}/\Gamma][\sigma(\Delta)/\Delta]\bigr)(\psi^\sigma\{f^\sigma\}) \\
&=& (((Q z:A[\bar{t}/\Gamma])\psi[\bar{t},z/\langle\Gamma,x:A\rangle])^\sigma \\
&=& (((Q x:A)\psi)[\bar{t}/\Gamma])^\sigma
\end{eqnarray*}
$\qed$

\begin{corollary} \label{sigmasubprj}
 Let $(\Sigma,\Pi)$ be a normal signature. Suppose
that $(\phi \; {\rm form}\; (\Gamma)) \in {\rm Form}(\Sigma^{(\infty)},\Pi)$. 
Then for  $\sigma=\sigma_\Sigma$:
$$(\phi^\sigma)\{({\bf p}_\Gamma(x:A))^\sigma\} = \phi^\sigma.$$
\end{corollary}
{\flushleft \bf Proof.}  By the theorem
$$(\phi^\sigma)\{({\bf p}_\Gamma(x:A))^\sigma\} =  (\phi[{\rm OV}(\Gamma)/\Gamma])^\sigma = \phi^\sigma.$$
$\qed$

\begin{theorem}  Let $(\Sigma,\Pi)$ be a normal signature, and let $T$ be a theory over $(\Sigma,\Pi)$.
If $(\phi \sequentunder{\Gamma} \psi) \in {\rm Thm}^*(\Sigma^{(\infty)},\Pi,T)$, then
$(\phi^\sigma \sequentunder{\Gamma^\sigma} \psi^\sigma) \in {\rm Thm}(\Sigma,\Pi,T^\sigma)$,
where $\sigma=\sigma_\Sigma$.
\end{theorem}
{\flushleft \bf Proof.}  Let ${\rm Thm}={\rm Thm}(\Sigma,\Pi,T^\sigma)$ and 
${\rm Thm}^*={\rm Thm}^*(\Sigma^{(\infty)},\Pi,T)$. The proof goes by induction on derivations in ${\rm Thm}^*$. The cases for the axioms are clear 
and so are the cases for propositional connectives and  the cut rule.
Consider the quantificational rules. For $(\phi\; {\rm form} \; (\Gamma)), 
(\psi\; {\rm form} \; (\Gamma, x:A)) \in {\rm Form}^*(\Sigma,\Pi)$:

\begin{itemize}
\item[(univ*)] $$\infer{\phi \sequentunder{\Gamma} (\forall x:A)\psi}{ \phi \sequentunder{\Gamma,x:A}\psi} \qquad
\infer{\phi \sequentunder{\Gamma,x:A}\psi}{\phi \sequentunder{\Gamma} (\forall x:A)\psi}$$
\item[(exis*)] $$\infer{(\exists x:A)\psi \sequentunder{\Gamma} \phi }{\psi\sequentunder{\Gamma,x:A}\phi} \qquad
\infer{\psi  \sequentunder{\Gamma,x:A}\phi}{(\exists x:A)\psi \sequentunder{\Gamma} \phi}$$
\end{itemize}

Case (univ*, left): Suppose that $(\phi \sequentunder{\Gamma} (\forall x:A)\psi) \in {\rm Thm}^*$
has been derived by (univ,left). Then $(\phi \sequentunder{\Gamma,x:A}\psi) \in {\rm Thm}^*$ by a shorter derivation, and  by inductive hypothesis 
$(\phi^\sigma \sequentunder{(\Gamma,x:A)^\sigma}\psi^\sigma) \in {\rm Thm}$. But 
$$(\Gamma,x:A)^\sigma = \Gamma^\sigma, \sigma(|\Gamma|):A[\sigma(\Gamma)/\Gamma].$$
We need to show
$$(\phi^\sigma \sequentunder{\Gamma^\sigma} ((\forall x:A)\psi)^\sigma) \in {\rm Thm}.$$
Now $$((\forall x:A)\psi)^\sigma = \bigl(\forall \, \sigma(|\Gamma|):A[\sigma(\Gamma)/\Gamma]\bigr)\phi^\sigma,$$
so we need only show 
$$(\phi^\sigma\{{\bf p}_{\Gamma^\sigma}(\sigma(|\Gamma|):A[\sigma(\Gamma)/\Gamma])\} \sequentunder{\Gamma^\sigma, \sigma(|\Gamma|):A[\sigma(\Gamma)/\Gamma]} \psi^\sigma) \in {\rm Thm}.$$
If we can show
\begin{equation} \label{phip}
\phi^\sigma\{{\bf p}_{\Gamma^\sigma}(\sigma(|\Gamma|):A[\sigma(\Gamma)/\Gamma])\}  = \phi^\sigma
\end{equation}
we are done.  But this follows from (\ref{ssubeq1}) and Corollary \ref{sigmasubprj}.

Case (univ*, right): using (\ref{phip}) as above this follows by inductive hypothesis
and (univ, right).

Case (exis*, left): using (\ref{phip}) as above this follows by inductive hypothesis
and (exis, left).

Case (exis*, right): using (\ref{phip}) as above this follows by inductive hypothesis
and (exis, right).

Case (subs*): Suppose that $(\phi[\bar{t}/\Gamma] \sequentunder{\Delta} 
\psi[\bar{t}/\Gamma]) \in {\rm Thm}^*$ has been derived by
the substitution rule from $(\phi \sequentunder{\Gamma} 
\psi)$ and $(\bar{t}:\Delta \to \Gamma) \in {\cal J}(\Sigma^{(\infty)})$. By the inductive hypothesis
 $$(\phi^\sigma \sequentunder{\Gamma^\sigma} 
\psi^\sigma) \in {\rm Thm}.$$
Thus  
$$(\sigma^{-1}(\Gamma) \circ \bar{t} \circ \sigma(\Delta) :\Delta^\sigma \to \Gamma^\sigma) \in {\cal J}(\Sigma),$$
and hence by (subs)
$$((\phi^\sigma)\{(\Delta^\sigma, \Gamma^\sigma, \sigma^{-1}(\Gamma) \bar{t}\sigma(\Delta))\} \sequentunder{\Delta^\sigma} 
(\psi^\sigma)\{(\Delta^\sigma, \Gamma^\sigma, \sigma^{-1}(\Gamma) \bar{t} \sigma(\Delta))\}) \in {\rm Thm}.$$
Then by Lemma \ref{sigmasub} we have
$$((\phi[\bar{t}/\Gamma])^\sigma \sequentunder{\Delta^\sigma} 
(\psi[\bar{t}/\Gamma])^\sigma) \in {\rm Thm},$$
as required. $\qed$

\pagebreak
{\flushleft \bf Bibliography}
\begin{itemize}
\item[] P.\ Aczel (2004). {\em Predicate logic with dependent sorts.} Unpublished manuscript 2004.
\item[] P.\ Aczel (2008). {\em Predicate logic over a type setup.} 
Slides from talk at Categorical and Homotopical Studies of Proof Theory, Barcelone, February 2008 . URL: http://mat.uab.cat/\~{ }kock/crm/hocat/prooftheory/aczel.pdf
\item[] P.\ Aczel (2011). {\em Generalised type setups for dependently sorted logic.} 
Slides from talk at TACL Marseille 2011. URL: http://www.cs.man.ac.uk/\~{ }petera/Recent-Slides/Marseille-2011-slides\_pap.pdf
\item[] P.\ Aczel (2012). {\em Syntax and Semantics - another look,
especially for dependent type theories}. Slides from a talk at the Newton Institute,
May 10, 2012. URL: 
https://www.newton.ac.uk/files/seminar/20120510163017001-153074.pdf
\item[] P.\ Aczel and M.\ Rathjen (2010). {\em Constructive Set Theory.} Manuscript.
\item[]  J.F.\ Belo (2008a). Dependently Sorted Logic. In: M. Miculan, I. Scagnetto, and F. Honsell (Eds.): {\em TYPES 2007,} LNCS 4941, pp. 33--50, Springer 2008. 
\item[] J.F.\ Belo (2008b) {\em Foundations of Dependently Sorted Logic.} PhD Thesis, Manchester 2008.
\item[] D.S.\ Bridges and F. Richman (1987). {\em Varieties of Constructive Mathematics.} Cambridge University Press.
\item[]  A.\ Buisse and P.\ Dybjer (2008). The interpretation of intuitionistic type theory in locally cartesian closed categories --- an intuitionistic perspective. {\em Electronic Notes in Theoretical Computer Science} 218(2008), 21--32.
\item[] J.\ Cartmell (1978).  {\em Generalized algebraic theories and contextual categories.}  DPhil.\ Thesis Oxford 1978.
\item[] J.\ Cartmell (1986).  Generalized algebraic theories and contextual categories. {\em Annals of Pure and Applied Logic,}  32(1986), 209--243.
\item[] P.\ Clairambault and P.\ Dybjer (2014). The biequivalence of locally cartesian closed categories and Martin-L\"of type theories. {\em Mathematical Structures in Computer Science} 24(2014), 1--54.
\item[] P.\ Dybjer (1996). Internal Type Theory.  In:  {\em Types for Proofs and Programs} (eds.\ S.\ Berardi and M.\ Coppo)
Lecture Notes in Computer Science Volume 1158, 1996, pp 120--134.
\item[] N.\ Gambino and P.\ Aczel (2006). The generalized type-theoretic interpretation of constructive set theory. {\em Journal of Symbolic Logic} 71(2006), 67--103.
\item[] M.\ Hofmann (1997). Syntax and semantics of dependent types. In: {\em Semantics and Logics of Computation} (eds.\ A.\ Pitts and P.\ Dybjer), 
Cambridge University Press 1997.
\item[] J.M.E.\ Hyland, P.T.\ Johnstone and A.M.\ Pitts (1980). Tripos theory.
{\em Mathematical Proceedings of the Cambridge Philosophical Society} 88(1980),  
205--232.
\item[] B.\ Jacobs (1999). {\em Categorical Logic and Type Theory}. Elsevier 1999.
\item[] P.T.\ Johnstone (2002). {\em Sketches of an Elephant: A Topos Theory Sketchbook, vol 1}. Oxford University Press 2002.
\item[] F.W.\ Lawvere (1964), An elementary theory of the category of sets, {\em  Proceedings of the National Academy of Science of the U.S.A.} 52 (1964), 1506 --1511.
\item[] M.E.\ Maietti and G.\ Sambin (2005). Towards a minimalist foundation for constructive mathematics. In: {\em From Sets and Types to Topology and Analysis} (eds. L.\ Crosilla and P.\ Schuster) Oxford University Press 2005, 91--114.
\item[] M.E.\ Maietti (2009). A minimalist two-level foundation for constructive mathematics. {\em Annals of Pure and Applied Logic} 160(2009), 319--354.
\item[] M.\ Makkai (1995). First-order logic with dependent sorts, with applications to category theory. Preprint 1995, version November 6. 201 pp.  Available from Makkai's webpages.
\item[] M.\ Makkai (1998). Towards a categorical foundation of mathematics. In:
{\em Logic Colloquium '95} (eds.\ J.A.\ Makowsky and E.V.\ Ravve) Lecture Notes in Logic, vol.\ 11, Association for Symbolic Logic 1998, 153--190.
\item[] M.\ Makkai (2013). The theory of abstract sets based on first-order logic with dependent types. Preprint 2013. Available from Makkai's webpages.
\item[] E.\ Palmgren (2012a).
Constructivist and structuralist foundations: Bishop's and Lawvere's theories of sets. {\em Annals of Pure and Applied Logic,} 163(2012), 1384--1399.
\item[] E.\ Palmgren (2012b). Proof-relevance of families of setoids and identity in type theory. {\em Archive for Mathematical Logic,} 51(2012), 35--47.
\item[] E.\ Palmgren (2016). {\em Categories with families, FOLDS and logic enriched type theory.} arXiv:1605.01586 [math.LO] (Earlier preprint version of the present paper.)
\item[] R.A.G.\ Seely (1983). Hyperdoctrines, natural deduction and the Beck condition. {\em Zeitschrift f\"ur mathematisches Logik und Grundlagen der Mathematik} 29(1983), 
505--542.
\end{itemize}

\pagebreak
\appendix
\section{Appendix: Proofs}

\begin{customthm}{\ref{UniqType}}
 (Unique typing lemma)  Consider derivations relative to a fixed signature $\Sigma$.
 Let $\Delta = y_1:B_1,\ldots, y_n: B_n$ be a context.
\begin{itemize}
\item[(a)] If $a:A \; (\Delta)$ and $a:A' \; (\Delta)$, then $A \equiv A'$
\item[(b)] If $(a_1,\ldots,a_n), (b_1,\ldots,b_n): \Gamma \rightarrow \Delta$ are contexts maps, with $a_i \equiv b_i$
for each $y_i \in {\rm TV}(\Delta)$, then $(a_1,\ldots,a_n) \equiv (b_1,\ldots,b_n)$.
\end{itemize}
\end{customthm}
{\flushleft \bf Proof.} We prove (a) and (b) simultaneously, by induction on derivations
in ${\cal J}(\Sigma)$. 

(a): Case $a=x_i$: Then $a:A \; (\Delta)$ and $a:A' \; (\Delta)$ can only have been derived by (R2) so $A \equiv A_i \equiv A' $.

Case $a=f(c_1,\ldots,c_k)$: Suppose that the declaration of $f$ in the signature $\Sigma$ is 
$$(\langle x_1:A_1\ldots, x_n:A_n\rangle,f,x_{i_1},\ldots,x_{i_k},U)$$ 
and that the last lines of the respective derivations,  obtained by rule (R5), are
$$f(a_{i_1},\ldots,a_{i_k}) : U[a_1,\ldots,a_n/x_1,\ldots,x_n] \;(\Delta)$$
and 
$$f(b_{i_1},\ldots,b_{i_k}) : U[b_1,\ldots,b_n/x_1,\ldots,x_n] \;(\Delta)$$
where $(c_1,\ldots,c_k) \equiv (a_{i_1},\ldots,a_{i_k}) \equiv  (b_{i_1},\ldots,b_{i_k})$
and with shorter derivations
$$(a_1,\ldots,a_n): \Delta \to \langle x_1:A_1,\ldots,x_n:A_n\rangle,$$
$$(b_1,\ldots,b_n): \Delta \to \langle x_1:A_1,\ldots,x_n:A_n\rangle.$$
Now, since ${\rm TV}(\langle x_1:A_1\ldots, x_n:A_n\rangle) \subseteq \{x_{i_1},\ldots,x_{i_k}\}$, we get by the inductive hypothesis
for (b) that  $(a_1,\ldots,a_n) \equiv (b_1,\ldots,b_n)$. Thus 
$$U[a_1,\ldots,a_n/x_1,\ldots,x_n] \equiv  U[b_1,\ldots,b_n/x_1,\ldots,x_n] $$
as required.

(b):  We do secondary induction on the length of the context $\Delta$.  If  the length is zero there is nothing to prove. Suppose
$\Delta = y_1:B_1,\ldots, y_{n+1}: B_{n+1}$, and that 
\begin{equation} \label{abmaps}
(a_1,\ldots,a_{n+1}), (b_1,\ldots,b_{n+1}): \Gamma \rightarrow \Delta
\end{equation} with 
\begin{equation} \label{abequ}
a_i \equiv b_i
\end{equation}
whenever $y_i \in {\rm TV}(\Delta)$. We have
\begin{equation} \label{TVDequ}
{\rm TV}(\Delta) =({\rm TV}(y_1:B_1,\ldots, y_{n}: B_{n}) \setminus {\rm V}(B_{n+1})) \cup \{y_{n+1}\}.
\end{equation}
Since $y_{n+1}$ is a member of this set,  $a_{n+1} \equiv b_{n+1}$.
By assumption (\ref{abmaps}) follows
$$a_{n+1}: B_{n+1}[a_1,\ldots,a_n/y_1,\ldots,y_n] \quad (\Gamma)$$
and
$$b_{n+1}: B_{n+1}[b_1,\ldots,b_n/y_1,\ldots,y_n] \quad (\Gamma).$$
Now since $a_{n+1} \equiv b_{n+1}$, part (a) gives
$$B_{n+1}[a_1,\ldots,a_n/y_1,\ldots,y_n] \equiv B_{n+1}[b_1,\ldots,b_n/y_1,\ldots,y_n]$$
From this follows that $a_i \equiv b_i$ for each $y_i \in {\rm v}(B_{n+1})$. Hence by (\ref{abequ}) and (\ref{TVDequ}) it follows that 
$a_i \equiv b_i$ for each $y_i \in {\rm TV}(y_1:B_1,\ldots, y_{n}: B_{n})$. Thus applying the
induction hypothesis for this shortening of $\Delta$, we get
$$(a_1,\ldots,a_n) \equiv (b_1,\ldots,b_n).$$
Hence $(a_1,\ldots,a_{n+1}) \equiv (b_1,\ldots,b_{n+1})$ as required. $\qed$

\begin{customthm}{\ref{substitution}}(Substitution lemma) Consider derivations in ${\cal J}(\Sigma)$ for a fixed signature $\Sigma$. Assume that $\vdash_{\bar{\ell}} \bar{s}:\Theta \to \Gamma$.
\begin{itemize}
\item[(a)] If $\vdash_k B \; {\rm type}\; (\Gamma)$, then 
$\vdash_{k+ {\rm max}(\bar{\ell})} B[\bar{s}/\Gamma] \; {\rm type}\; (\Theta)$.
\item[(b)] If $\vdash_k b: B\; (\Gamma)$, then 
$\vdash_{k+ {\rm max}(\bar{\ell})} b[\bar{s}/\Gamma] :B[\bar{s}/\Gamma] \;  (\Theta)$.
\end{itemize}
\end{customthm}
{\flushleft \bf Proof.}  The argument goes by induction on derivations in ${\cal J}(\Sigma)$. 
 Suppose that $\Gamma = x_1:A_1,\ldots,x_n:A_n$ and $\bar{s} =s_1,\ldots,s_n$.
 Then the assumption $\vdash_{\bar{\ell}} \bar{s}:\Theta \to \Gamma$ implies that
\begin{equation} \label{subst}
\begin{array}{l} 
 \vdash_{\ell_3} s_1: A_1\quad (\Theta)  \\ 
 \vdash_{\ell_4} s_2: A_2[s_1/x_1] \quad (\Theta)  \\ 
 \quad \vdots \\
\vdash_{\ell_{n+2}}  s_n: A_n[s_1,\ldots,s_{n-1}/x_1,\ldots,x_{n-1}] \quad (\Theta)  
 \end{array}
 \end{equation}
 
 Only rules R3 -- R5 can be used to obtain the hypotheses of (a) and (b), so we consider those.

{\flushleft \em Rule R3: } The conclusion of the rule is 
$$ x_i:A_i \; (x_1:A_1,\ldots,x_n:A_n).$$
By (\ref{subst}) we have
$$s_i: A_i[s_1,\ldots,s_{i-1}/x_1,\ldots,x_{i-1}] \;(\Theta).$$
Noting again that $A_i[s_1\ldots,s_{i-1}/x_1,\ldots,x_{i-1}]  = A_i [s_1,\ldots,s_n/x_1,\ldots,x_n]$ and that
$$s_i=x_i[s_1,\ldots,s_n/x_1,\ldots,x_n]$$
we obtain the desired conclusion
$$  x_i[s_1,\ldots,s_n/x_1,\ldots,x_n] :  A_i[s_1,\ldots,s_n/x_1,\ldots,x_n]\quad (\Theta).$$
Hence, in this case
$$\vdash_{\ell_{i+2}} b[\bar{s}/\Gamma] :B[\bar{s}/\Gamma] \;  (\Theta).$$
Now $\ell_{i+2} \le k+ {\rm max}(\bar{\ell})$, so also the upper bound of proof height is established for this case.
      
{\flushleft \em Rule R4: } Consider an application of the rule:

$$\infer[({\rm R4})]{S(b_{i_1},\ldots,b_{i_q}) \; {\rm type}\;(\Gamma) }{ 
(b_1,\ldots,b_m): \Gamma \to\, y_1:B_1\ldots, y_m:B_m},
$$
where $(y_1:B_1\ldots, y_m:B_m, i_1,\ldots,i_q, S) \in \Sigma$, and such that
$$\vdash_k S(b_{i_1},\ldots,b_{i_q}) \; {\rm type}\;(\Gamma),$$
and
\begin{equation} \label{R4sub}
\vdash_{\bar{r}} (b_1,\ldots,b_m): \Gamma \to\, y_1:B_1\ldots, y_m:B_m.
\end{equation}
Thus $k\ge 1+ {\rm max}(\bar{r})$. For each $i=1,\ldots,m$, we have
$$\vdash_{r_{i+2}} b_i :  B_i[b_1,\ldots,b_{i-1}/y_1,\ldots,y_{i-1}] \; (\Gamma).$$
By the inductive hypothesis using (\ref{R4sub})
\begin{equation} \label{R4b}
\vdash_{r_{i+2} + {\rm max}(\bar{\ell})}  b_i[s_1,\ldots,s_n/x_1,\ldots,x_n]: B_i[b_1,\ldots,b_{i-1}/y_1,\ldots,y_{i-1}][s_1,\ldots,s_n/x_1,\ldots,x_n]
 \; (\Theta)
\end{equation}
By syntactic substitution the expression in (\ref{R4b}) is
\begin{equation} \label{R4bb}
\begin{split}
& b_i[s_1,\ldots,s_n/x_1,\ldots,x_n]: \\
&B_i[b_1[s_1,\ldots,s_n/x_1,\ldots,x_n],\ldots,b_{i-1}[s_1,\ldots,s_n/x_1,\ldots,x_n]/y_1,\ldots,y_{i-1}] \; (\Theta)
\end{split}
\end{equation}
Thus we may apply (R4) to (\ref{R4b})  for $i=1,\ldots,m$ and obtain
$$\vdash_{t} S(b_{i_1}[s_1,\ldots,s_n/x_1,\ldots,x_n],\ldots,b_{i_q}[s_1,\ldots,s_n/x_1,\ldots,x_n]) \; {\rm type} \; (\Theta)$$
where 
\begin{equation} \label{R4t}
t= 1+ {\rm max}(r_1,r_2,r_3+{\rm max}(\bar{\ell}),\ldots,
r_{n+2}+{\rm max}(\bar{\ell})).
\end{equation}
But this is nothing else than
$$\vdash_{t} S(b_{i_1},\ldots,b_{i_q})[s_1,\ldots,s_n/x_1,\ldots,x_n] \; {\rm type}\; (\Theta)$$
and the estimate $t$ is bound by
$$k+{\rm max}(\bar{\ell}) \ge 1+ {\rm max}(\bar{r}) +{\rm max}(\bar{\ell}) \ge t$$
as desired.

{\flushleft \em Rule R5: } An application of the rule looks like this:

$$\infer[({\rm R5})]{f(b_{i_1},\ldots,b_{i_q}): U[b_1,\ldots,b_m/y_1,\ldots,y_m]\;(\Gamma)}{ 
 (b_1,\ldots,b_m): \Gamma \to\,  y_1:B_1\ldots, y_m:B_m &
 U[b_1,\ldots,b_m/y_1,\ldots,y_m]\; {\rm type}\; (\Gamma)}$$
 where $(y_1:B_1\ldots, y_m:B_m,f,i_1,\ldots,i_q,U) \in \Sigma$. 
By inductive hypothesis and the property of substitution we have, exactly as above, (\ref{R4bb}) for all $i=1,\ldots,m$. Moreover the inductive hypothesis and substitution property gives
$$U[b_1[s_1,\ldots,s_n/x_1,\ldots,x_n],\ldots,b_m[s_1,\ldots,s_n/x_1,\ldots,x_n]/y_1,\ldots,y_m]\; {\rm type}\; (\Theta).$$
Thus applying (R5) to these judgements, one obtains
\begin{equation} 
\begin{split}
&\vdash_t f(b_{i_1}[s_1,\ldots,s_n/x_1,\ldots,x_n],\ldots,b_{i_q}[s_1\ldots,s_n/x_1,\ldots,x_n]): \\
&
\qquad U[b_{i_1}[s_1,\ldots,s_n/x_1,\ldots,x_n],\ldots,b_{i_q}[s_1,\ldots,s_n/x_1,\ldots,x_n]/y_1,\ldots,y_m]\; (\Theta)
\end{split}
\end{equation}
which is the same as
\begin{equation} 
\begin{split}
&\vdash_t  f(b_{i_1},\ldots,b_{i_k})[s_1,\ldots,s_n/x_1,\ldots,x_n]: \\
&
\qquad U[b_1,\ldots,b_m/y_1,\ldots,y_m][s_1,\ldots,s_n/x_1,\ldots,x_n] \; (\Theta)
\end{split}
\end{equation}
which gives the required conclusion, as $k+{\rm max}(\bar{\ell}) \ge t$.
$\qed$

\begin{customthm}{\ref{weakening}}(Weakening lemma) 
 Consider derivations in ${\cal J}(\Sigma)$ for a fixed signature $\Sigma$ with unrestricted
 variable system. Suppose that 
$\vdash B \; {\rm type}\;(\Gamma)$. Let $y$ be a variable not in ${\rm v}(\Gamma, \Theta)$. 
\begin{itemize}
\item[(a)] If $\vdash \Gamma, \Theta \; {\rm context}$, then $\vdash\Gamma, y:B, \Theta \; {\rm context}$ 
\item[(b)] If $\vdash A \; {\rm type}\; (\Gamma,\Theta)$, then $\vdash A \; {\rm type}\; (\Gamma,y:B, \Theta)$
\item[(c)] If $\vdash a:A\; (\Gamma,\Theta)$, then $\vdash a:A \; (\Gamma,y:B, \Theta)$
\end{itemize}
\end{customthm}
{\flushleft \bf Proof.}  Induction on derivations. Note that $x \in \varphi(FV(\Delta))$ is equivalent to $x \in V \setminus FV(\Delta)$.

{\flushleft \em Rule R1.} The conclusion is $\langle \rangle \; {\rm context}$, so $\Gamma$ and $\Delta$ must be empty contexts. We
have by assumption $B \; {\rm type}\; ()$, so we can conclude
$$y:B \; {\rm context}$$
as desired.

{\flushleft \em Rule R2.} The rule application  has the form
$$\infer[({\rm R2})]{\Gamma,\Theta' ,x:A\; {\rm context}}{\Gamma, \Theta' \; {\rm context}\;&
   A\; {\rm type} \;(\Gamma,\Theta')}$$
 where $\Theta= \Theta',x:A$, and $x \notin {\rm v}(\Gamma,\Theta')$.
Now $y \notin {\rm v}(\Gamma,\Theta') \supseteq {\rm v}(\Gamma,\Theta)$, so we have by inductive hypothesis
 $$\Gamma,y:B, \Theta' \; {\rm context} \; \mbox{ and } A\; {\rm type}  \; (\Gamma,y:B,\Theta') .$$
 Clearly $y$ is not in ${\rm v}(\Gamma,\Theta)$ so R2 may be applied to obtain
 $$\Gamma,y:B,\Theta' ,x:A\; {\rm context}$$
as required.

{\flushleft \em Rule R3.} Let rule application have the form
$$\infer[({\rm R3})]{x_i:A_i\; (x_1:A_1\ldots, x_n:A_n)  }{x_1:A_1\ldots, x_n:A_n \; {\rm context} }$$
where $\Gamma= x_1:A_1\ldots, x_j:A_j$ and $\Theta=  x_{j+1}:A_{j+1}\ldots, x_n:A_n$.
By inductive hypothesis
$$ x_1:A_1\ldots, x_j:A_j,y:B, x_{j+1}:A_{j+1}\ldots, x_n:A_n  \; {\rm context}.$$
By applying R3 we obtain the desired
$$x_i:A_i\; (x_1:A_1\ldots, x_j:A_j,y:B, x_{j+1}:A_{j+1}\ldots, x_n:A_n).$$

{\flushleft \em Rule R4.}  Suppose we have the rule instance

$$\infer[({\rm R4})]{S(a_{i_1},\ldots,a_{i_k}) \; {\rm type}\;(\Gamma,\Theta)}{ 
(a_1,\ldots,a_n): \Gamma,\Theta  \to \Delta
 }$$
 where $(\Delta,i_1,\ldots,i_k,S) \in \Sigma$.
 Invocation of the inductive hypotheses gives us the premisses of the rule but where $\Gamma,\Theta$ have been replaced by $\Gamma,y:B,\Theta$. By R4 we may thus conclude
 $$S(a_{i_1},\ldots,a_{i_k}) \; {\rm type}\; (\Gamma,y:B,\Theta) $$
 as required.
 
{\flushleft \em Rule R5.}  Suppose we have the application
$$\infer[({\rm R5})]{f(a_{i_1},\ldots,a_{i_k}) : U[a_1,\ldots,a_n/x_1,\ldots,x_n]  \;(\Gamma,\Theta)}{ 
(a_1,\ldots,a_n): \Gamma,\Theta  \to  \Delta
& U[a_1,\ldots,a_n/x_1,\ldots,x_n]  \; {\rm type} \; (\Gamma,\Theta)}$$
 where $(\Delta,i_1,\ldots,i_k,f, U) \in \Sigma$.
 Again, invocation of the inductive hypotheses gives us the premisses of the rule but where $\Gamma,\Theta$ have been replaced by $\Gamma,y:B,\Theta$. By R5 we may thus conclude
 $$   f(a_{i_1},\ldots,a_{i_k}) : U[a_1,\ldots,a_n/\Delta] \; (\Gamma,y:B,\Theta)$$
as desired. 

This concludes the proof.
$\qed$

\begin{customthm}{\ref{strengthening}}(Strengthening lemma) 
Consider derivations in ${\cal J}(\Sigma)$ for a fixed signature $\Sigma$ unrestricted
 variable system.
Let $y$ be a variable not in ${\rm v}(\Theta, A, a)$. 
\begin{itemize}
\item[(a)] If $\vdash \Gamma, y:B, \Theta \; {\rm context}$, then $\vdash \Gamma, \Theta \; {\rm context}$. 
\item[(b)] If $\vdash A \; {\rm type}\; (\Gamma,y:B, \Theta)$, then  $\vdash A \; {\rm type}\; (\Gamma,\Theta)$.
\item[(c)] If $\vdash a:A \; (\Gamma,y:B, \Theta)$, then $\vdash a:A\; (\Gamma,\Theta)$.
\end{itemize}
\end{customthm}
{\flushleft \bf Proof.}   Induction on derivations.

{\flushleft \em Rule R1.} The conclusion is $\langle \rangle \; {\rm context}$ so this case does not apply.

{\flushleft \em Rule R2.}  In case $\Theta$ is empty the rule application has the form
 $$\infer[({\rm R2})]{\Gamma ,y:B\; {\rm context}}{\Gamma \; {\rm context}\;&
   B\; {\rm type} \;(\Gamma)}$$
 where  $y \in V \setminus {\rm v}(\Gamma)$.
 Thus we have
 $$\Gamma \; {\rm context}$$
by a shorter derivation.

In case $\Theta= y:B,\Theta' $ is not empty, the rule application has the form
$$\infer[({\rm R2})]{\Gamma,y:B,\Theta' ,x:A\; {\rm context}}{\Gamma,y:B, \Theta' \; {\rm context}\;&
   A\; {\rm type} \;(\Gamma,y:B, \Theta')}$$
where $x \in V \setminus {\rm v}(\Gamma,B,\Theta')$.
By inductive hypothesis we have 
$$ \Gamma, \Theta' \; {\rm context} \qquad
   A\; {\rm type} \;(\Gamma, \Theta') $$
 so applying (R2) to this we get the desired judgement
$$\Gamma,\Theta' ,x:A\; {\rm context}$$

{\flushleft \em Rule R3.} Suppose rule application has the form
$$\infer[({\rm R3})]{x_i:A_i\; (x_1:A_1\ldots, x_n:A_n)  }{x_1:A_1\ldots, x_n:A_n \; {\rm context} }$$
where $\Gamma= x_1:A_1\ldots, x_{j-1}:A_{j-1}$, $x_j:A_j = y:B$ and $\Theta=  x_{j+1}:A_{j+1}\ldots, x_n:A_n$.  
Now $y\, (=x_j)$ is not in ${\rm v}(\Theta, A_i,x_i)$. 
By inductive hypothesis
$$\Gamma,\Theta  \; {\rm context}$$
and
$$ A_i \;  {\rm type} \; (\Gamma,\Theta).$$
In particular by the variable condition  $i \ne j$, so we can apply (R3) to obtain the desired
$$x_i:A_i\; (\Gamma,\Theta).$$

{\flushleft \em Rule R4.}  Suppose we have the rule instance

$$\infer[({\rm R4})]{S(a_1,\ldots,a_n) \; {\rm type}\;(\Gamma,y:B,\Theta)}{  
(a_1,\ldots,a_n): \Gamma,y:B,\Theta  \to  \langle x_1:A_1\ldots, x_n:A_n\rangle
}$$
where  $(\langle x_1:A_1\ldots, x_n:A_n\rangle,1,\ldots,n,S) \in \Sigma$.
 Reading out the assumption
$(a_1,\ldots,a_n): \Gamma,y:B,\Theta  \to  \langle x_1:A_1\ldots, x_n:A_n\rangle$ it
says that
\begin{equation} \label{GBTcontext}
\Gamma,y:B,\Theta  \; {\rm context}
\end{equation}
and for $i=1,\ldots,n$
\begin{equation} \label{aAGBT}
a_i: A_i[a_1,\ldots,a_{i-1}/x_1,\ldots,x_{i-1}] \; (\Gamma,y:B,\Theta).
\end{equation}
By assumption $y$ is not in ${\rm v}(\Theta, S(a_1,\ldots,a_n))$.  Here we use the fact that the signature is regular.
Thus $y$ is not in ${\rm v}(\Theta, a_i,  A_i[a_1,\ldots,a_{i-1}/x_1,\ldots,x_{i-1}])$, so by inductive hypothesis:
$$\Gamma, \Theta  \; {\rm context}$$
and for $i=1,\ldots,n$
$$a_i: A_i[a_1,\ldots,a_{i-1}/x_1,\ldots,x_{i-1}] \; (\Gamma,\Theta).$$
We can then apply R4 to obtain
$$S(a_1,\ldots,a_n) \; {\rm type}\;(\Gamma,\Theta)$$
as desired.
  
{\flushleft \em Rule R5.}  Suppose we have the application
$$\infer[({\rm R5})]{f(a_1,\ldots,a_n) : U[a_1,\ldots,a_n/x_1,\ldots,x_n]  \;(\Gamma,y:B,\Theta)}{ 
(a_1,\ldots,a_n): \Gamma,y:B,\Theta  \to  \langle x_1:A_1\ldots, x_n:A_n\rangle 
}$$
where $(\langle x_1:A_1\ldots, x_n:A_n\rangle,f,1,\ldots,n,U) \in \Sigma$.
Spelling out the assumption about the context map, we get as above (\ref{GBTcontext}) and (\ref{aAGBT}). By assumption $y$ is not in 
$${\rm v}(\Theta, f(a_1,\ldots,a_n), U[a_1,\ldots,a_n/x_1,\ldots,x_n])= {\rm v}(\Theta, a_1,\ldots,a_n, U[a_1,\ldots,a_n/x_1,\ldots,x_n]).$$
Thus we can apply the inductive hypothesis to obtain
 $$\Gamma, \Theta   \; {\rm context}$$
 and 
  $$(a_1,\ldots,a_n): \Gamma,\Theta  \to  \langle x_1:A_1\ldots, x_n:A_n\rangle .$$
Thus applying R5, we obtain 
$$f(a_1,\ldots,a_n) : U[a_1,\ldots,a_n/x_1,\ldots,x_n]  \;(\Gamma,\Theta)$$
This concludes the proof.
$\qed$

\begin{customthm}{ \ref{interchange}}(Interchange lemma)
Consider derivations in ${\cal J}(\Sigma)$ for a fixed signature $\Sigma$ with unrestricted
 variable system. Suppose that $x \notin {\rm v}(C)$. 
\begin{itemize}
\item[(a)] If $\vdash \Gamma,x;B, y:C, \Theta \; {\rm context}$, then $\vdash \Gamma,y:C, x:B, \Theta \; {\rm context} $. 
\item[(b)] If $\vdash A \; {\rm type}\; (\Gamma,x:B, y:C, \Theta) $, then  $\vdash A \; {\rm type}\; (\Gamma,y:C, x:B,\Theta)$.
\item[(c)] If $\vdash a:A \; (\Gamma, x:B,y:C, \Theta) $, then $\vdash a:A\; (\Gamma,y:C, x:B,\Theta)$.
\end{itemize}
\end{customthm}
{\flushleft \bf Proof.}   Induction on derivations.

{\flushleft \em Rule R1.} The conclusion is $\langle \rangle \; {\rm context}$  and hence this case does not apply.

{\flushleft \em Rule R2.}  In case $\Theta$ is empty the rule application has the form
 $$\infer[({\rm R2})]{\Gamma ,x:B, y:C\; {\rm context}}{\Gamma,x:B \; {\rm context}\;&
   C\; {\rm type} \;(\Gamma, x:B)}$$
   where  $y \in V \setminus {\rm v}(\Gamma,x:B)$.
   Using the Strengthening lemma on 
   $C\; {\rm type} \;(\Gamma, x:B),$
we obtain  $C\; {\rm type} \;(\Gamma)$. As $y \notin {\rm v}(B)$,  the Weakening lemma on $\Gamma,x:B \; {\rm context}$ gives
$$\Gamma,y:C, x:B \; {\rm context}$$
as required. Assume now that $\Theta = \Theta', z:A$. Then the rule application has the form
$$\infer[({\rm R2})]{\Gamma ,x:B, y:C, \Theta', z:A\; {\rm context}}{\Gamma ,x:B, y:C, \Theta' \; {\rm context}\;&
   A\; {\rm type} \;(\Gamma, x:B, y:C, \Theta') }$$
where $z \in V \setminus {\rm v}(\Gamma,x:B,, y:C, \Theta')$.
Applying the inductive hypothesis to the premisses, and then applying R2 gives the required result.

{\flushleft \em Rule R3.} Suppose the rule application has the form
$$\infer[({\rm R3})]{x_i:A_i\; (x_1:A_1\ldots, x_n:A_n)  }{x_1:A_1\ldots, x_n:A_n \; {\rm context} }$$
where $\Gamma= x_1:A_1\ldots, x_{j-2}:A_{j-2}$,  $x_{j-1}:A_{j-1} = x:B$, $x_{j}:A_{j} = y:C$,  
and $\Theta=  x_{j+1}:A_{j+1}\ldots, x_n:A_n$. By inductive hypothesis
$$\Gamma,y:C, x:B, \Theta \; {\rm context} \qquad A_i \; {\rm type}\; (\Gamma,y:C, x:B, \Theta).$$ 
An application of R3 gives
$$x_i:A_i \;(\Gamma,y:C, x:B, \Theta).$$

{\flushleft \em Rule R4.} Assume that the rule application is
$$\infer[({\rm R4})]{ S(a_{i_1},\ldots,a_{i_k}) \; {\rm type} \;(\Gamma,x:B, y:C, \Theta)}{ 
(a_1,\ldots,a_n): 
(\Gamma, x:B, y:C, \Theta) \to \Delta 
 }$$
 where $(\Delta,S,i_1,\ldots,i_k) \in \Sigma$
 and $\Delta = \langle x_1:A_1\ldots, x_n:A_n\rangle$.
Spelling out the assumption on the context map it says 
\begin{equation} \label{BCTcontext}
\Gamma,x:B, y:C, \Theta \; {\rm context}
\end{equation}
and for all $i=1,\ldots,n$,
\begin{equation} \label{aABCT}
a_i: A_i[a_1,\ldots,a_{i-1}/x_1,\ldots,x_{i-1}] \; (\Gamma,x:B, y:C, \Theta).
\end{equation}
By inductive hypothesis we can switch $x:B$ and $y:C$ in these judgements and apply R4 to 
achieve
$$S(a_{i_1},\ldots,a_{i_k}) \; {\rm type} \;(\Gamma, y:C, x:B, \Theta)$$

{\flushleft \em Rule R5.} 
 Suppose we have the rule application
 $$\infer[({\rm R5})]{f(a_{i_1},\ldots,a_{i_k}) : U[a_1,\ldots,a_n/\Delta]  \;(\Gamma, x:B, y:C,\Theta)}{(a_1,\ldots,a_n): \Gamma,x:B, y:C, \Theta  \to  \Delta & U[a_1,\ldots,a_n/\Delta]  \; {\rm type} \; (\Gamma, x:B, y:C,\Theta)
}$$
where 
$(\Delta,f,i_1,\ldots,i_k,U) \in \Sigma$ and $\Delta = \langle x_1:A_1\ldots, x_n:A_n\rangle$.
 Reading out the assumption on the context map gives as above (\ref{BCTcontext}) and (\ref{aABCT}).
 By inductive assumption we can again switch $x:B$ and $y:C$ in the judgements and thereby
get 
$$(a_1,\ldots,a_n): \Gamma, y:C, x:B, \Theta  \to  \langle x_1:A_1\ldots, x_n:A_n\rangle$$
and
$$U[a_1,\ldots,a_n/x_1,\ldots,x_n]  \; {\rm type} \; (\Gamma, y:C, x:B, \Theta).$$
Hence by R5 we obtain:
$$f(a_{i_1},\ldots,a_{i_k}) : U[a_1,\ldots,a_n/x_1,\ldots,x_n]  \;(\Gamma, y:C, x:B,\Theta),$$
as was required. $\qed$

\begin{customlemma}{\ref{expbyfunct}} Suppose that $\Sigma$ is a signature such that ${\cal J}^*(\Sigma) = {\cal J}(\Sigma)$. Assume that $D=(\Delta,f,\bar{i},U)$ is a predeclaration, where $f$ is not declared in $\Sigma$, and where $(U \; {\rm type}\; (\Delta)) \in {\cal J}^*(\Sigma)$. Then for
$\Sigma' = \Sigma \cup \{D\}$
$${\cal J}^*(\Sigma') = {\cal J}(\Sigma').$$
\end{customlemma}
{\flushleft \bf Proof. }  The  ${\cal J}^*(\Sigma') \supseteq {\cal J}(\Sigma')$ is direct as observed above. We prove ${\cal J}^*(\Sigma') \subseteq {\cal J}(\Sigma')$. Denote by
$R$ the application of R5* with the function symbol $f$.
The proof goes by double induction on the property $P(m,n):$
$$\forall {\cal U}(\Sigma' \vdash^{*,R}_m {\cal U} \; \& \; 
\Sigma' \vdash^{*}_n {\cal U} \Longrightarrow \Sigma' \vdash  {\cal U}).$$
Note that for $m=0$, $\Sigma' \vdash^{*,R}_m {\cal U}$ implies that $R$ is not used
so $\Sigma \vdash {\cal U}$. Hence by the assumption 
${\cal J}^*(\Sigma) = {\cal J}(\Sigma)$, we have ${\cal U} \in {\cal J}(\Sigma) \subseteq {\cal J}(\Sigma')$ as required.

To prove the inductive step assume that $\forall n\; P(m,n)$. We prove 
$\forall n\; P(m+1,n)$. For $n=0$, $\Sigma' \vdash^{*}_n {\cal U}$ implies 
that ${\cal U} = (\langle \rangle \; {\rm context})$, so trivially $\Sigma' \vdash  {\cal U}$.
Now assume that $P(m+1,n)$ holds. We show that $P(m+1,n+1)$: assume
that 
\begin{equation} \label{P1}
 \Sigma' \vdash^{*,R}_{m+1} {\cal U} 
\end{equation}
\begin{equation} \label{P2}
\Sigma' \vdash^{*}_{n+1} {\cal U}.
\end{equation}
We establish $\Sigma' \vdash {\cal U}$ by case distinction on the last rule 
applied.

Rule (R1): Then ${\cal U}= (\langle \rangle \; {\rm context})$, so trivially $\Sigma' \vdash {\cal U}$.

Rule (R2): Suppose that this last rule application  is
$$\infer{\Gamma,x:A\; {\rm context}}{ \Gamma\; {\rm context} &
    A\; {\rm type}\; (\Gamma)}.$$
Then $\Sigma' \vdash^{*,R}_{m+1} \Gamma\; {\rm context}$  and 
$\Sigma' \vdash^{*,R}_{m+1} A\; {\rm type}\; (\Gamma)$.
Also 
$\Sigma' \vdash^{*}_{n} \Gamma\; {\rm context}$ and
$\Sigma' \vdash^{*}_{n} A\; {\rm type}\; (\Gamma)$. By $P(m+1,n)$ we get
$\Sigma' \vdash \Gamma\; {\rm context}$ and 
$\Sigma' \vdash A\; {\rm type}\; (\Gamma)$. Thus using (R2)
$$\Sigma' \vdash \Gamma,x:A\; {\rm context}$$
as required.

Rule (R3): Assume the last rule application gives
$$\infer{x_i:A_i  \;(x_1:A_1,\ldots, x_n:A_n)}{  x_1:A_1,\ldots, x_n:A_n \; {\rm context}}.$$
Then $\Sigma' \vdash^{*,R}_{m+1} x_1:A_1,\ldots, x_n:A_n \; {\rm context}$  and
$\Sigma' \vdash^{*}_{n} x_1:A_1,\ldots, x_n:A_n \; {\rm context}.$
By $P(m+1,n)$ we get
$\Sigma' \vdash x_1:A_1,\ldots, x_n:A_n \; {\rm context}$. Hence applying (R3)
yields $$\Sigma' \vdash x_i:A_i  \;(x_1:A_1,\ldots, x_n:A_n).$$ 

Rule (R4): Suppose the last rule application is
$$\infer[ \begin{array}{l} (\Delta,S, \bar{i}) \mbox{ in  } \Sigma\\
  \end{array}   ]{S(\bar{a}_{\,\bar{i}}) \; {\rm type} \;(\Gamma) }{\bar{a}: \Gamma \to \Delta 
}.$$
Then $\Sigma' \vdash^{*,R}_{m+1} \bar{a}: \Gamma \to \Delta 
$  and
$\Sigma' \vdash^{*}_{n} \bar{a}: \Gamma \to \Delta.$
By $P(m+1,n)$ this gives
$\Sigma' \vdash \bar{a}: \Gamma \to \Delta$. Hence by (R4):
$\Sigma' \vdash S(\bar{a}_{\,\bar{i}}) \; {\rm type} \;(\Gamma)$.

Rule (R5*): Next assume the last rule application is
$$\infer[ \begin{array}{l} (\Delta',f',\bar{j},U') \mbox{ in  } \Sigma' \\
 \end{array}]{ f'(\bar{a}_{\,\bar{j}}) : U'[\bar{a}/\Delta'] \;(\Gamma) }{ 
 \bar{a}: \Gamma \to\Delta' \\
}$$
If $f' \ne f$ we have $\Sigma' \vdash^{*,R}_{m+1} \bar{a}: \Gamma \to \Delta' 
$   and $\Sigma' \vdash^{*}_{n} \bar{a}: \Gamma \to \Delta' .$
By $P(m+1,n)$ this gives
$\Sigma' \vdash \bar{a}: \Gamma \to \Delta' $.  Since $f' \ne f$, we must have
$(\Delta',f',\bar{j},U') \in \Sigma$. Thus  $\Sigma \vdash U'\; {\rm type}\; (\Delta')$,
since ${\cal J}(\Sigma) = {\cal J}^*(\Sigma)$.
We can now use the Substitution Lemma to obtain 
$\Sigma' \vdash U'[\bar{a}/\Delta']\; {\rm type}\; (\Gamma)$.
Now using (R5) we get as desired
$\Sigma' \vdash f'(\bar{a}_{\, \bar{j}}):U'[\bar{a}/\Delta']\; {\rm type}\; (\Gamma)$.

If $f'=f$, then $\Delta'=\Delta$ and $\bar{j}=\bar{i}$, and  
we have $\Sigma' \vdash^{*,R}_{m} \bar{a}: \Gamma \to \Delta' 
$   and $\Sigma' \vdash^{*}_{n} \bar{a}: \Gamma \to \Delta'$. Since 
$(U\; {\rm type} \;(\Delta)) \in {\cal J}^*(\Sigma)$, no applications of $R$
are need, so
$$\Sigma \vdash^{*,R}_0 U\; {\rm type} \;(\Delta)$$
and further, for some $k$:
$$\Sigma \vdash^{*}_k U\; {\rm type} \;(\Delta).$$
By the Substitution Lemmas we have
$$\Sigma' \vdash^{*,R}_{0+m} U[\bar{a}/\Delta]\; {\rm type} \;(\Gamma)$$
and 
$$\Sigma' \vdash^{*}_{k+m} U[\bar{a}/\Delta]\; {\rm type} \;(\Gamma).$$
Applying $P(m,k+m)$, we obtain
$$\Sigma' \vdash U[\bar{a}/\Delta]\; {\rm type} \;(\Gamma).$$
Now using (R5) we get
$$\Sigma' \vdash f'(\bar{a}_{\, \bar{j}}):U'[\bar{a}/\Delta']\; {\rm type}\; (\Gamma)$$
as desired.  $\qed$

\begin{customlemma}{\ref{expbytype}} Suppose that $\Sigma$ is a signature such that ${\cal J}^*(\Sigma) = {\cal J}(\Sigma)$. Assume that $D=(\Delta,\bar{i},S)$ is a predeclaration, where $S$ is not declared in $\Sigma$, and where $(\Delta \; {\rm context}) \in {\cal J}^*(\Sigma)$. Then for
$\Sigma' = \Sigma \cup \{D\}$
$${\cal J}^*(\Sigma') = {\cal J}(\Sigma').$$
\end{customlemma}
{\flushleft \bf Proof.} The  ${\cal J}^*(\Sigma') \supseteq {\cal J}(\Sigma')$ is direct as observed above. We prove ${\cal J}^*(\Sigma') \subseteq {\cal J}(\Sigma')$. 
The proof goes by double induction on the property $P(n):$
$$\forall {\cal U}(\Sigma' \vdash^{*}_n {\cal U} \Longrightarrow \Sigma' \vdash  {\cal U}).$$
For $n=0$, $\Sigma' \vdash^{*}_n {\cal U}$ implies 
that ${\cal U} = (\langle \rangle \; {\rm context})$, so trivially $\Sigma' \vdash  {\cal U}$.
Now assume that $P(n)$ holds. We show that $P(n+1)$: assume
that 
\begin{equation} \label{P2}
\Sigma' \vdash^{*}_{n+1} {\cal U}.
\end{equation}
We establish $\Sigma' \vdash {\cal U}$ by case distinction on the last rule 
applied.

Rule (R1): Then ${\cal U}= (\langle \rangle \; {\rm context})$, so trivially $\Sigma' \vdash {\cal U}$.

Rule (R2): Suppose that this last rule application  is
$$\infer{\Gamma,x:A\; {\rm context}}{ \Gamma\; {\rm context} &
    A\; {\rm type}\; (\Gamma)}.$$
Then 
Also 
$\Sigma' \vdash^{*}_{n} \Gamma\; {\rm context}$ and
$\Sigma' \vdash^{*}_{n} A\; {\rm type}\; (\Gamma)$. By $P(n)$ we get
$\Sigma' \vdash \Gamma\; {\rm context}$ and 
$\Sigma' \vdash A\; {\rm type}\; (\Gamma)$. Thus using (R2)
$$\Sigma' \vdash \Gamma,x:A\; {\rm context}$$
as required.

Rule (R3): Assume the last rule application gives
$$\infer{x_i:A_i  \;(x_1:A_1,\ldots, x_n:A_n)}{  x_1:A_1,\ldots, x_n:A_n \; {\rm context}}.$$
Then 
$\Sigma' \vdash^{*}_{n} x_1:A_1,\ldots, x_n:A_n \; {\rm context}.$
By $P(n)$ we get
$\Sigma' \vdash x_1:A_1,\ldots, x_n:A_n \; {\rm context}$. Hence applying (R3)
yields $$\Sigma' \vdash x_i:A_i  \;(x_1:A_1,\ldots, x_n:A_n).$$ 

Rule (R4): Suppose the last rule application is
$$\infer[ \begin{array}{l} (\Delta,S, \bar{i}) \mbox{ in  } \Sigma\\
  \end{array}   ]{S(\bar{a}_{\,\bar{i}}) \; {\rm type} \;(\Gamma) }{\bar{a}: \Gamma \to \Delta 
}.$$
Then
$\Sigma' \vdash^{*}_{n} \bar{a}: \Gamma \to \Delta.$
By $P(n)$ this gives
$\Sigma' \vdash \bar{a}: \Gamma \to \Delta$. Hence by (R4):
$\Sigma' \vdash S(\bar{a}_{\,\bar{i}}) \; {\rm type} \;(\Gamma)$.

Rule (R5*): Next assume the last rule application is
$$\infer[ \begin{array}{l} (\Delta',f,\bar{j},U) \mbox{ in  } \Sigma' \\
 \end{array}]{ f(\bar{a}_{\,\bar{j}}) : U[\bar{a}/\Delta'] \;(\Gamma) }{ 
 \bar{a}: \Gamma \to\Delta' \\
}$$
We have $\Sigma' \vdash^{*}_{n} \bar{a}: \Gamma \to \Delta' .$
By $P(n)$ this gives
$\Sigma' \vdash \bar{a}: \Gamma \to \Delta' $.  We must have
$(\Delta',f,\bar{j},U) \in \Sigma$. Thus  $\Sigma \vdash U\; {\rm type}\; (\Delta')$,
since ${\cal J}(\Sigma) = {\cal J}^*(\Sigma)$.
We can now use the Substitution Lemma to obtain 
$\Sigma' \vdash U[\bar{a}/\Delta']\; {\rm type}\; (\Gamma)$.
Now using (R5) we get 
$$\Sigma' \vdash f(\bar{a}_{\, \bar{j}}):U[\bar{a}/\Delta']\; {\rm type}\; (\Gamma).$$
as desired.  $\qed$

\dontshow{
\begin{customthm}{\ref{F2signature}}
For any FOLDS vocabulary $K$, the pre-signature $\Sigma_K$ is a signature. 
\end{customthm}
{\flushleft \bf Proof.} To prove the theorem it is enough to prove that for each $k=1,\ldots, N$, the pre-context $\Gamma_{B_k}$ is a context relative
to the signature $\{{\cal D}_{B_1},\ldots,{\cal D}_{B_{k-1}}\}$. For $k=1$, $n(B_1)=0$, so $\Gamma_{B_1}$ is the empty context, and hence ${\cal D}_{B_1}$ is a type declaration
relative to the empty signature $[]$. Suppose the statement holds for all $k=1,\ldots,m$. We show $\Gamma_{B_{m+1}}$ is a
context relative to the signature $\{{\cal D}_{B_1},\ldots,{\cal D}_{B_{m}}\}$. By construction
\begin{equation}\label{Bcontext}
\Gamma_{B_{m+1}} = x_1: C_1({\rm OV}( \Gamma_{C_1})^{x_1}),\ldots, x_n: C_n({\rm OV}( \Gamma_{C_n})^{x_n})
\end{equation}
where $x_i = x_i^{B_{m+1}}$,  $C_i = {\rm cod}(x_i)$ and $n=n(B_{m+1})$.
We prove by induction that each initial segment of (\ref{Bcontext}) is a context
\begin{equation}\label{inBcontext}
\Gamma_{B_{m+1}} |_k \;\; {\rm context}  
\end{equation}
relative to $\{{\cal D}_{B_1},\ldots,{\cal D}_{B_{m}}\}$
where $k=0,\ldots,n$. For $k=0$, this is clear. Suppose (\ref{inBcontext}) holds for all $k=0,\ldots, r$.
To prove $\Gamma_{B_{m+1}} |_{r+1} \;\; {\rm context}$, it suffices to prove
$$C_{r+1}({\rm OV}( \Gamma_{C_{r+1}})^{x_{r+1}}) \; {\rm type} \qquad (\Gamma_{B_{m+1}} |_r) $$
Since $C_{r+1} <^* B_{m+1}$, we have that $D_{C_{r+1}}$ is declared in  $\{{\cal D}_{B_1},\ldots,{\cal D}_{B_{m}}\}$,
and hence it is enough to prove that the following is a context map
\begin{equation} \label{cmap2be}
{\rm OV}( \Gamma_{C_{r+1}})^{x_{r+1}}:  (\Gamma_{B_{m+1}} |_r)  \rightarrow  \Gamma_{C_{r+1}}.
\end{equation}
Now suppose that
\begin{equation} \label{Ccontext}
\Gamma_{C_{r+1}} = y_1: D_1({\rm OV}( \Gamma_{D_1})^{y_1}),\ldots, y_s: D_s({\rm OV}( \Gamma_{D_s})^{y_s})
\end{equation}
where  $y_i = x_i^{C_{r+1}}$,  $D_i = {\rm cod}(y_i)$ and $s=n(C_{r+1})$.  Note that $s \le r$. As $${\rm v}(D_i({\rm OV}( \Gamma_{D_i})^{y_i})) \subseteq \{y_1,\ldots,y_{i-1}\},$$
 we can
write
\begin{equation} \label{zsdef}
y_{f(i,1)},\ldots, y_{f(i,t(i))} =_{\rm def} {\rm OV}(\Gamma_{D_i})^{y_i}
\end{equation}
where $f(i,k) \in \{y_1,\ldots,y_{i-1}\}$ and $t(i) \le i-1$.
Spelling out (\ref{cmap2be}) as basic judgements it becomes the conjunction of the following judgements
$$
\begin{array}{lll} 
 &(C)&\Gamma_{B_{m+1}} |_r  \; {\rm context} \\
 &(E_1)&  y_1x_{r+1}:D_1(y_{f(1,1)},\ldots, y_{f(1,t(1))}) \; (\Gamma_{B_{m+1}} |_r) \\ 
&(E_2)&  y_2x_{r+1}: D_2(y_{f(2,1)},\ldots, y_{f(2,t(2))})[y_1x_{r+1}/y_1] \;(\Gamma_{B_{m+1}} |_r) \\
&&\quad \vdots \\
&(E_s)& y_sx_{r+1}:D_s(y_{f(s,1)},\ldots, y_{f(s,t(s))})[y_1x_{r+1},\ldots,y_{s-1}x_{r+1}/y_1,\ldots,y_{s-1}] \;(\Gamma_{B_{m+1}} |_r) 
 \end{array}
 $$
 First $(C)$ follows from the inductive hypothesis.  By performing the syntactic substitutions,  the problem is reduced to proving
 $$
\begin{array}{lll} 
 &(E_1)&  y_1x_{r+1}:D_1(y_{f(1,1)},\ldots, y_{f(1,t(1))}) \; (\Gamma_{B_{m+1}} |_r) \\ 
&(E_2)&  y_2x_{r+1}: D_2(y_{f(2,1)}x_{r+1},\ldots, y_{f(2,t(2))}x_{r+1}) \;(\Gamma_{B_{m+1}} |_r) \\
&&\quad \vdots \\
&(E_s)& y_sx_{r+1}:D_s(y_{f(s,1)}x_{r+1},\ldots, y_{f(s,t(s))}x_{r+1}) \;(\Gamma_{B_{m+1}} |_r) 
 \end{array}
 $$

By the form of (\ref{Ccontext}) it follows that $t(1)=0$, and 
$$D_1 \; {\rm type} \quad ().$$ 
Hence by weakening
$$D_1 \; {\rm type} \quad(\Gamma_{B_{m+1}} |_r).$$ 
As $y_1x_{r+1}: B_{m+1} \rightarrow D_1$ is a non-identity map, there must be some $x_j$ with $x_j=y_1x_{r+1}$. Let $j$ be the smallest such. We have 
$D_1 = C_j$. Moreover $D_1 <^* C_{r+1}$, so $j \le r$. Thus
$$y_1x_{r+1}: D_1 \quad(\Gamma_{B_{m+1}} |_r).$$ 
This verifies $(E_1)$.

For every $i=1,\ldots,s$, we have a non-identity arrow $y_ix_{r+1}$ from $B_{m+1}$ to $D_i$ via $C_{r+1}$. Since  $D_i <^* C_{r+1}$ there is 
$g: \{1,\ldots,s\}\rightarrow \{1,\ldots,r\}$ such that
$$x_{g(i)} = y_ix_{r+1} \qquad C_{g(i)} = D_i.$$
Now since $i<j$ implies $C_{g(i)} <^* C_{g(j)}$ it follows that $g$ is strictly increasing.
The remaining problem is now to check that
$$
\begin{array}{lll} 
 &(E_1)&  x_{g(1)}:C_{g(1)} \; (\Gamma_{B_{m+1}} |_r) \\ 
&(E_2)&  x_{g(2)}: C_{g(2)}(x_{g(f(2,1))},\ldots,x_{g(f(2,t(2)))}) \;(\Gamma_{B_{m+1}} |_r) \\
&&\quad \vdots \\
&(E_s)& x_{g(s)}: C_{g(s)}(x_{g(f(s,1))},\ldots,x_{g(f(s,t(s)))}) \;(\Gamma_{B_{m+1}} |_r) 
 \end{array}
 $$
Case $(E_1)$ follows by the above. Let $i \in \{2,\ldots,s\}$. We prove $(E_i)$:
$$x_{g(i)}: C_{g(i)}(x_{g(f(i,1))},\ldots,x_{g(f(i,t(i)))}) \;(\Gamma_{B_{m+1}} |_r).$$
We have 
\begin{eqnarray*}
C_{g(i)}(x_{g(f(i,1))},\ldots,x_{g(f(i,t(i)))}) &=& D_i(y_{f(i,1)}x_{r+1},\ldots,y_{f(i,t(i))}x_{r+1})  \\
 &=& D_i({\rm OV}(\Gamma_{D_i})^{y_i})^{x_{r+1}} \\
  &=& D_i({\rm OV}(\Gamma_{D_i})^{x_{g(i)}}) \\
  &=& C_{g(i)}({\rm OV}(\Gamma_{C_{g(i)}})^{x_{g(i)}})
 \end{eqnarray*}
By (\ref{Ccontext}) we have
$$ x_{g(i)}: C_{g(i)}({\rm OV}(\Gamma_{C_{g(i)}})^{x_{g(i)}}) \qquad (\Gamma_{B_{m+1}} |_r)$$
so we are done. $\qed$
}

\begin{customthm}{\ref{catofcwfs}}
The cwfs and cwf morphisms form a category.
\end{customthm}
{\flushleft \bf Proof.} 
We first check that composition as defined above is
well-defined. We have
$$(F' \circ F)(\top_{\cal C}) =F'(\top_{{\cal C}'})= \top_{{\cal C}''}.$$
Suppose $f: \Delta \rightarrow \Gamma$,
and $A \in \Ty(\Gamma)$. Then
\begin{eqnarray*} 
(\sigma' \circ \sigma)_{\Delta}(A\{f\}) &=& \sigma'_{F(\Delta)}(\sigma_\Delta(A\{f\})) \\
&=&
\sigma'_{F(\Delta)}(\sigma_\Gamma(A)\{F(f)\}) \\
&=&
\sigma'_{F(\Gamma)}(\sigma_\Gamma(A))\{F'(F(f))\} \\
&=&
(\sigma' \circ \sigma)_\Gamma(A)\{(F' \circ F)(f))\} \\
\end{eqnarray*}
\begin{eqnarray*}
(F' \circ F)(\Gamma.A) &= & F'(F(\Gamma).\sigma_{\Gamma}(A)) \\
                                     &= & F'(F(\Gamma)).\sigma'_{F(\Gamma)}(\sigma_{\Gamma}(A)) \\
                                     &= & (F'\circ F)(\Gamma).(\sigma' \circ \sigma)_{\Gamma}(A) \\
\end{eqnarray*}
\begin{eqnarray*}
(F' \circ F)({\rm p}_{\Gamma}(A)) 
      & = & F'({\rm p}_{F(\Gamma)}(\sigma_{\Gamma}(A))) \\
      & = & {\rm p}_{F'(F(\Gamma))}(\sigma'_{F(\Gamma)}(\sigma_{\Gamma}(A))) \\
      & = & {\rm p}_{(F' \circ F)(\Gamma)}((\sigma' \circ \sigma)_{\Gamma}(A)) \\
\end{eqnarray*}
\begin{eqnarray*}
 ((\theta' \circ \theta)_{\Gamma,A}(a))\{(F'\circ F)(f)\} & = &
    (\theta'_{F(\Gamma),\sigma_{\Gamma}(A)} (\theta_{\Gamma,A}(a)))\{(F'(F(f))\} \\
    & = & \theta'_{F(\Delta),\sigma_\Delta(A\{f\})}(\theta_{\Gamma,A}(a)\{F(f)\}) \\
    & = & \theta'_{F(\Delta),\sigma_\Delta(A\{f\})}(\theta_{\Delta,A\{f\}}(a\{f\})) \\
    & = & (\theta' \circ \theta)_{F(\Delta),\sigma_\Delta(A\{f\})}(a\{f\}) \\
    & = & (\theta' \circ \theta)_{\Delta,A\{f\}}(a\{f\}) \\  
\end{eqnarray*}
\begin{eqnarray*} 
(\theta' \circ \theta)_{\Gamma.A, A\{{\rm p}(A)\}}({\rm v}_A) &=& 
\theta'_{F(\Gamma.A),\sigma_{\Gamma.A}(A\{{\rm p}(A)\})}(\theta_{\Gamma.A, A\{{\rm p}(A)\}}({\rm v}_A)) \\
&=& 
\theta'_{F(\Gamma.A),\sigma_{\Gamma.A}(A\{{\rm p}(A)\})}({\rm v}_{\sigma_\Gamma(A)}) \\
&=& {\rm v}_{\sigma_{F(\Gamma)}(\sigma_\Gamma(A))} \\
&=& {\rm v}_{(\sigma' \circ \sigma)_\Gamma(A)}.
\end{eqnarray*}
Suppose that $f:\Delta \to \Gamma$ and $a \in \Tm(\Delta, A\{f\})$. Then
\begin{eqnarray*} 
(F' \circ F)(\langle f,a\rangle_A) &=& F'(\bigl\langle F(f),\theta_{\Delta, A\{f\}}(a) \bigr\rangle_{\sigma_\Gamma(A)}) \\
&=& \bigl\langle F'(F(f)),\theta_{F(\Delta), \sigma_\Gamma(A)\{F(f)\}}(\theta_{\Delta, A\{f\}}(a))\bigr\rangle_{\sigma'_{F(\Gamma)}(\sigma_\Gamma(A))} \\
&=& \bigl\langle F'(F(f)),\theta'_{F(\Delta), \sigma_\Delta(A\{f\})}(\theta_{\Delta, A\{f\}}(a))\bigr\rangle_{\sigma'_{F(\Gamma)}(\sigma_\Gamma(A))} \\
&=& \bigl\langle (F' \circ F)(f),(\theta' \circ \theta)_{\Delta, A\{f\}}(a)\bigr\rangle_{(\sigma' \circ \sigma)_\Gamma(A)}
\end{eqnarray*} 
This proves that the composition of two cwf morphisms is again a cwf morphism.
The identity laws are
\begin{equation} \label{id1}
  (I_{{\cal C}'}, \iota^{{\cal C}'}, \varepsilon^{{\cal C}'}) \circ (F, \sigma, \theta) = (F, \sigma, \theta)
\end{equation}
and
\begin{equation} \label{id2}
    (F, \sigma, \theta) \circ (I_{{\cal C}}, \iota^{{\cal C}}, \varepsilon^{{\cal C}}) = (F, \sigma, \theta)
\end{equation}
These hold obviously for the functor components of the equalities.
As for (\ref{id1}) we
$$(\iota^{{\cal C}'} \circ \sigma)_{\Gamma}(A) = \iota^{{\cal C}'}_{F(\Gamma)}(\sigma_{\Gamma}(A)) = 1_{\Ty(F(\Gamma))}(\sigma_{\Gamma}(A)) = \sigma_{\Gamma}(A),$$
and
$$(\varepsilon^{{\cal C}'} \circ \theta)_{\Gamma,A} = \varepsilon^{{\cal C}'}_{F(\Gamma),\sigma_{\Gamma}(A)} \circ \theta_{\Gamma,A}= 1_{\Tm F(\Gamma),\sigma_{\Gamma}(A)} \circ \theta_{\Gamma,A} =  \theta_{\Gamma,A}$$
As for (\ref{id2}) we have similarly
$$(\sigma \circ \iota^{{\cal C}})_{\Gamma}(A) = 
\sigma_{I(\Gamma)}(\iota^{{\cal C}}_\Gamma(A)) =
\sigma_{\Gamma}(1_{\Ty(\Gamma)}(A)) = \sigma_{\Gamma}(A),$$
$$(\theta \circ \varepsilon^{{\cal C}})_{\Gamma,A} = 
 \theta_{I(\Gamma),\iota_{\Gamma}(A)} \circ \varepsilon_{\Gamma,A}
 = 
 \theta_{\Gamma,A} \circ 1_{\Tm(\Gamma,A)} = \theta_{\Gamma, A}.$$

To prove associativity assume that 
$(F'', \sigma'', \theta'') : ({\cal C}'',\Ty'',\Tm'') \to ({\cal C}''',\Ty''',\Tm''')$
is a third cwf morphism. We shall thus show
$$((F'', \sigma'', \theta'') \circ (F', \sigma', \theta')) \circ (F, \sigma, \theta)=
(F'', \sigma'', \theta'') \circ ((F', \sigma', \theta') \circ (F, \sigma, \theta)).$$
The equality of the first component is just associativity of functor composition.  The second component equality is
\begin{eqnarray*}
((\sigma'' \circ \sigma') \circ \sigma)_\Gamma(A) &= & 
(\sigma'' \circ \sigma')_{F(\Gamma)}(\sigma_\Gamma(A)) \\
&= & 
\sigma''_{F'(F(\Gamma))} (\sigma'_{F(\Gamma)}(\sigma_\Gamma(A))) \\
&=& \sigma''_{(F' \circ F)(\Gamma)} ((\sigma' \circ \sigma)_\Gamma(A)) \\
&=& (\sigma'' \circ (\sigma' \circ \sigma))_\Gamma(A) \\
\end{eqnarray*}
The third equality is given by
\begin{eqnarray*}
((\theta'' \circ \theta') \circ \theta)_{\Gamma,A} &= &  
(\theta'' \circ \theta')_{F(\Gamma),\sigma_{\Gamma}(A)} \circ \theta_{\Gamma,A}\\
&= &  (\theta''_{F'(F(\Gamma)), \sigma'_{F(\Gamma)}(\sigma_\Gamma(A))} \circ \theta'_{F(\Gamma),\sigma_{\Gamma}(A)}) \circ \theta_{\Gamma,A}\\ 
&= & \theta''_{F'(F(\Gamma)), \sigma'_{F(\Gamma)}(\sigma_\Gamma(A))} \circ (\theta'_{F(\Gamma),\sigma_{\Gamma}(A)} \circ \theta_{\Gamma,A}) \\
&= & \theta''_{(F' \circ F)(\Gamma), (\sigma' \circ \sigma)_\Gamma(A)} \circ (\theta' \circ \theta)_{\Gamma,A} \\
&= & (\theta'' \circ (\theta' \circ \theta))_{\Gamma,A} 
\end{eqnarray*}

$\qed$       

\begin{customlemma}{\ref{varstandard}} 
Suppose that $\Sigma$ is a signature, and that
$\sigma: {\mathbb N} \to V$ is a fresh sequence of variables. If $(\Gamma\; {\rm context}) \in {\cal J}(\Sigma)$, then variable ordering $\sigma$ for $\Gamma$ forms
a context map
$$(\sigma(\Gamma): \Gamma^\sigma \to \Gamma) \in {\cal J}(\Sigma),$$
whose inverse 
$(\sigma^{-1}(\Gamma):\Gamma \to \Gamma^\sigma) \in {\cal J}(\Sigma)$, 
is given by $\sigma^{-1}(\Gamma) = {\rm OV}(\Gamma)$.
\end{customlemma}
{\flushleft \bf Proof.} First we prove $\sigma(\Gamma)$ is a context map. This is proved by induction on derivations using the Substitution lemma above: case (R1) is trivial. Suppose $(\Gamma,x:A \; {\rm context}) \in {\cal J}(\Sigma)$. Then
by shorter derivations $(\Gamma \; {\rm context}) \in {\cal J}(\Sigma)$  and
$(A \;{\rm type}\;(\Gamma)) \in {\cal J}(\Sigma)$. Hence by inductive hypothesis 
$(\sigma(\Gamma):\Gamma^\sigma \to \Gamma) \in {\cal J}(\Sigma)$. Thus by the Substitution lemma
$(A[\sigma(\Gamma)/\Gamma] \;{\rm type}\; (\Gamma^\sigma)) \in {\cal J}(\Sigma)$.
Also by inductive hypothesis $(\Gamma^\sigma \; {\rm context}) \in {\cal J}(\Sigma)$.
So by (R2) and  $\sigma(|\Gamma|) \in {\rm Fresh}(\Gamma^\sigma)$,
\begin{equation} \label{vs1}
(\langle \Gamma^\sigma, \sigma(|\Gamma|):A[\sigma(\Gamma)/\Gamma]\rangle \; {\rm context}) \in {\cal J}(\Sigma).
\end{equation}
Thus $(\langle \Gamma,x:A \rangle^\sigma \; {\rm context}) \in {\cal J}(\Sigma)$.
By (\ref{vs1}), and (R3) follows
\begin{equation} \label{vs2}
\sigma(|\Gamma|):A[\sigma(\Gamma)/\Gamma] \; \bigl(\langle \Gamma^\sigma, \sigma(|\Gamma|):A[\sigma(\Gamma)/\Gamma]\rangle\bigr).
\end{equation}
Then by $(\sigma(\Gamma):\Gamma^\sigma \to \Gamma) \in {\cal J}(\Sigma)$  and (\ref{vs2}) it
follows that
$$(\sigma(\langle \Gamma,x:A \rangle):\langle \Gamma,x:A \rangle^\sigma \to \langle \Gamma,x:A \rangle) \in {\cal J}(\Sigma)$$
as required. 

Next we show that for $(\Gamma\; {\rm context}) \in {\cal J}(\Sigma)$,
\begin{eqnarray*}
 \sigma(\Gamma) \circ \sigma^{-1}(\Gamma) & = & 1_{\Gamma^\sigma}, \\
 \sigma^{-1}(\Gamma) \circ \sigma(\Gamma) & = & 1_{\Gamma}. \\
\end{eqnarray*}
Suppose that $\Gamma=x_1:A_1,\ldots,x_n:A_n$, so 
$$\sigma^{-1}(\Gamma) = {\rm OV}(\Gamma)= (x_1,\ldots,x_n)$$
Then 
$$\sigma(\Gamma) = (y_1,\ldots, y_n) = {\rm OV}(\Gamma^\sigma)$$
where $$y_k=\sigma(k)$$ for $k=1,\ldots,n$.
Thus 
\begin{eqnarray*}
\sigma(\Gamma) \circ \sigma^{-1}(\Gamma) &= & (y_1[x_1,\ldots,x_n/y_1,\ldots,y_n],\ldots,y_n[x_1,\ldots,x_n/y_1,\ldots,y_n]) \\
&= &(x_1,\ldots,x_n) = 1_{\Gamma}.
\end{eqnarray*}
and 
\begin{eqnarray*}
\sigma^{-1}(\Gamma) \circ \sigma(\Gamma)   &= & (x_1[y_1,\ldots,y_n/x_1,\ldots,x_n],\ldots,x_n[y_1,\ldots,y_n/x_1,\ldots,x_n]) \\
&= &(y_1,\ldots,y_n) = 1_{\Gamma^\sigma}.
\end{eqnarray*}

Finally we prove that $(\sigma^{-1}(\Gamma):\Gamma \to \Gamma^\sigma) \in {\cal J}(\Sigma)$. With the above notation
$$\Gamma^\sigma = y_1:A_1,y_2:A_2[y_1/x_1],\ldots,y_n:A_n[y_1,\ldots,y_{n-1}/x_1,\ldots,x_{n-1}].$$
Since both  $(\Gamma\; {\rm context}), (\Gamma^\sigma\; {\rm context})  \in {\cal J}(\Sigma)$, we need only to show in ${\cal J}(\Sigma)$
\begin{equation} 
\begin{array}{lll} 
 &  x_1: A_1 \; (\Gamma) \\ 
&  x_2: A_2[y_1/x_1][x_1/y_1] \;(\Gamma) \\
&\quad \vdots \\
& x_n: A_n[y_1,\ldots,y_{n-1}/x_1,\ldots,x_{n-1}][x_1,\ldots,x_{n-1}/y_1,\ldots,y_{n-1}] \;(\Gamma)
 \end{array}
 \end{equation}
But each of the substitutions cancel so the judgements follow by the assumption rules.
$\qed$

\dontshow{
\section{Appendix: Type constructions on cwfs}

We shall consider the type constructions $\Pi$, $\Sigma$, $+$ and $\nattype_k$  below.

\medskip
A CwF {\em supports $\nattype_k$-types} if for every $\Gamma \in {\cal C}$ there is $\nattype_{k,\Gamma} \in \Ty(\Gamma)$ and $k$ elements
\begin{equation} \label{Nkintro}
0_{k,\Gamma},\ldots,(k-1)_{k,\Gamma} \in \Tm(\Gamma, \nattype_{k,\Gamma})
\end{equation}
such that for any $f: \Delta \rightarrow \Gamma$,
\begin{itemize}
\item[($\nattype_k$-subst)] $\nattype_{k,\Gamma}\{f\} = \nattype_{k,\Delta}$
\item[($i_k$-subst)] $i_{k,\Gamma}\{f\} = i_{k,\Delta}$.
\end{itemize}
Moreover for any $C \in \Ty(\Gamma.\nattype_{k,\Gamma})$ and any $M_i \in \Tm(\Gamma, C\{\langle 1_\Gamma,i_{k,\Gamma}\rangle\})$ ($i=0,\ldots,k-1$) and any 
$P \in \Tm(\Gamma,\nattype_{k,\Gamma})$, there is
\begin{equation} \label{Nkelim}
{\rm R}_{k,\Gamma,C}(P,M_0,\ldots,M_{k-1}) \in \Tm(\Gamma, C\{\langle 1_\Gamma, P \rangle\}
\end{equation}
which is such that
\begin{itemize}
\item[($\nattype_k$-conv)] ${\rm R}_{k,\Gamma,C}(i_{k,\Gamma},M_0,\ldots,M_{k-1}) = M_i$
\end{itemize}
and for $f: \Delta \rightarrow \Gamma$,
\begin{itemize}
\item[($\nattype_k$-subst)] ${\rm R}_{k,\Gamma,C}(P,M_0,\ldots,M_{k-1})\{f\}= 
{\rm R}_{k,\Delta,C\{{\rm q}(f)\}}(P\{f\},M_0\{f\},\ldots,M_{k-1}\{f\})$.
\end{itemize}

\begin{lemma} \label{N0N1L} Let ${\cal C}$ be a cwf that supports $\nattype_0$ and $\nattype_1$. Suppose that 
$\Gamma \in {\cal C}$ and $A \in \Ty(\Gamma)$.
\begin{itemize}
\item[(a)]  $\Tm(\Gamma.\nattype_0, A\{{\rm p}(\nattype_0)\})$ is inhabited,
\item[(b)] $\Tm(\Gamma.A, \nattype_1\{{\rm p}(A)\})$ is inhabited.
\end{itemize}
\end{lemma}

{\flushleft \bf Proof.}  As for (a): Let $C=A\{{\rm p}(\nattype_{\Gamma, 0})
{\rm p}(\nattype_{\Gamma.\nattype_0, 0})\}$. So 
$C \in \Ty(\Gamma.\nattype_0.\nattype_0)$. Now, 
$$v_{\nattype_0} \in \Tm(\Gamma.\nattype_{\Gamma, 0}, \nattype_{\Gamma, 0}\{{\rm p}(\nattype_{\Gamma, 0})\})= \Tm(\Gamma.\nattype_{\Gamma, 0}, \nattype_{\Gamma.\nattype_0, 0}).$$
So by $\nattype_0$-elimination (\ref{Nkelim})

$${\rm R}_{0,\Gamma.\nattype_0,C}(v_{\nattype_0}) \in \Tm(\Gamma.\nattype_0, C\{\langle 1_{\Gamma.\nattype_0}, v_{\nattype_0} \rangle\}).$$ 
But 
$$C\{\langle 1_{\Gamma.\nattype_0}, v_{\nattype_0} \rangle\}=
A\{{\rm p}(\nattype_{\Gamma, 0})
{\rm p}(\nattype_{\Gamma.\nattype_0, 0})\}\{\langle 1_{\Gamma.\nattype_0}, v_{\nattype_0} \rangle\} =A\{{\rm p}(\nattype_{\Gamma, 0})\}$$
Hence $\Tm(\Gamma.\nattype_0, A\{{\rm p}(\nattype_0)\})$ is inhabited.

As for (b):  By $\nattype_1$-introduction (\ref{Nkintro}),
$$0_{1,\Gamma.A} \in \Tm(\Gamma.A, \nattype_{\Gamma.A, 1}).$$
Since $\nattype_{1,\Gamma}\{{\rm p}(A)\} =\nattype_{\Gamma.A, 1}$ we have that
 $$\Tm(\Gamma.A, \nattype_{1,\Gamma}\{{\rm p}(A)\})$$
 is inhabited as desired. $\qed$.

\medskip
A cwt {\em supports $\Sigma$-types} if for $A \in \Ty(\Gamma)$ and $B \in \Ty(\Gamma.A)$ there is 
a type  $\Sigma(A,B) \in \Ty(\Gamma)$, and for $M \in\Tm(\Gamma,A)$ and
$N \in \Tm(\Gamma, B\{\langle1_\Gamma, M\rangle_A\})$ there is an element ${\rm Pair}_{A,B}(M,N) \in \Tm(\Gamma, \Sigma(A,B))$. These constructions should satisfy for any $f: \Delta \to \Gamma$ 
\begin{itemize}
\item[($\Sigma$-subst)] $\Sigma(A,B)\{f\} = \Sigma(A\{f\},B\{f.A\})$,
\item[(Pair-subst)] ${\rm Pair}_{A,B}(M,N)\{f\} = {\rm Pair}_{A\{f\},B\{f.A\}}(M\{f\}, N\{f\})$
\end{itemize}
Next we need to specify the
elimination operation ${\rm E}$. First we construct a substitution
$${\rm pair}_{A,B}: \Gamma.A.B \rightarrow \Gamma.\Sigma(A,B)$$
by 
$${\rm pair}_{A,B} =_{\rm def} \langle {\rm p}(A.B), 
{\rm Pair}_{A\{{\rm p}(A.B)\} ,B\{{\rm p}(A.B).A\}}({\rm v}_{A.B},{\rm v}_B) \rangle.$$
For any $C \in \Ty(\Gamma.\Sigma(A,B))$, for any $P \in \Tm(\Gamma, \Sigma(A,B))$ and for any $K \in \Tm(\Gamma.A.B, C\{{\rm pair}_{A,B}\})$ there is an element
$${\rm E}_{A,B,C}(P,K) \in \Tm(\Gamma, C\{\langle 1_\Gamma,P\rangle_{\Sigma(A,B)}\})$$
such that
\begin{itemize}
\item[($\Sigma$-conv)] ${\rm E}_{A,B,C}({\rm Pair}_{A,B}(M,N),K)= 
K\{\langle\langle 1_\Gamma, M\rangle_A,N\rangle_B\}$.
\end{itemize}
Moreover for any $f: \Delta \rightarrow \Gamma$, we require
\begin{itemize}
\item[(${\rm E}$-subst)] ${\rm E}_{A,B,C}(P,K)\{f\} = 
{\rm E}_{A\{f\},B\{f.A\},C\{f.\Sigma(A,B)\}}(P\{f\}, K\{f.A.B\}).$
\end{itemize}

\medskip
For $A, B \in \Ty(\Gamma)$ we introduce the notation
$$(A \times B) =_{\rm def} \Sigma(A, B\{{\rm p}(A)\}) \in \Ty(\Gamma)$$

\begin{example} \label{firstproj}
{\em First projection. Let $C=A\{{\rm p}(\Sigma(A,B))\}$ in the above.
Then for $P \in \Tm(\Gamma,\Sigma(A,B))$, we have
$${\rm E}_{A,B,C}(P,K) \in
 \Tm(\Gamma, C\{\langle 1_\Gamma,P\rangle_{\Sigma(A,B)}\}) =
 \Tm(\Gamma, A),$$
where
$$K \in \Tm(\Gamma.A.B, C\{{\rm pair}_{A,B}\}) = 
\Tm(\Gamma.A.B, A\{{\rm p}(A.B)\})
.$$
Thus we may let $K={\rm v}_{A.B}$. Then
\begin{eqnarray*}
{\rm E}_{A,B,C}({\rm Pair}_{A,B}(M,N),{\rm v}_{A.B}) &= &
{\rm v}_{A.B}\{\langle\langle 1_\Gamma, M\rangle_A,N\rangle_B\}\\
&= &
{\rm v}_{A}\{{\rm p}(B)\}\{\langle\langle 1_\Gamma, M\rangle_A,N\rangle_B\}\\
&= &
{\rm v}_{A}\{\langle 1_\Gamma, M\rangle_A\} = M,\\
\end{eqnarray*}
as required. Let $$\pi_{1,A,B}(P) =_{\rm def} 
{\rm E}_{A,B,A\{{\rm p}(\Sigma(A,B))\}}(P,{\rm v}_{A.B}) \in \Tm(\Gamma, A).$$
We have 
$$\pi_{1,A,B}({\rm Pair}_{A,B}(M,N)) =M.$$
}
\end{example}
\begin{example}\label{secondproj}{\em Second projection. Suppose $A \in \Ty(\Gamma)$ and 
$B \in \Ty(\Gamma.A)$. Using the first projection we construct
$$\langle {\rm p}(\Sigma(A,B)), \pi_1({\rm v}_{\Sigma(A,B)}) \rangle : \Gamma.\Sigma(A,B) \to
\Gamma.A$$
Let $$C= B\{\langle {\rm p}(\Sigma(A,B)), \pi_1({\rm v}_{\Sigma(A,B)}) \rangle\} \in 
\Ty(\Gamma.\Sigma(A,B)).$$
Now for $P \in \Tm(\Gamma,\Sigma(A,B))$
\begin{eqnarray*}
{\rm E}_{A,B,C}(P,K) &\in &
\Tm(\Gamma, C\{\langle 1_\Gamma,P\rangle_{\Sigma(A,B)}\})\\
&=&
\Tm(\Gamma,B\{\langle {\rm p}(\Sigma(A,B)), \pi_1({\rm v}_{\Sigma(A,B)}) \rangle\}\{\langle 1_\Gamma,P\rangle_{\Sigma(A,B)}\}) \\
&=&
\Tm(\Gamma,B\{\langle 1_\Gamma, {\rm E}({\rm v}_{\Sigma(A,B)},{\rm v}_{A.B})\{\langle 1_\Gamma,P\rangle_{\Sigma(A,B)}\} \rangle\}) \\
&=&
\Tm(\Gamma,B\{\langle 1_\Gamma, {\rm E}({\rm v}_{\Sigma(A,B)}\{\langle 1_\Gamma,P\rangle_{\Sigma(A,B)}\},{\rm v}_{(A.B)\{\langle 1_\Gamma,P\rangle_{\Sigma(A,B)}\}})\rangle\}) \\
&=&
\Tm(\Gamma,B\{\langle 1_\Gamma, {\rm E}(P,{\rm v}_{(A.B)\{\langle 1_\Gamma,P\rangle_{\Sigma(A,B)}\}})\rangle\}) \\
&=&
\Tm(\Gamma,B\{\langle 1_\Gamma, \pi_{1,A,B}(P)\rangle\}) \\
\end{eqnarray*}
where 
\begin{eqnarray*}
K &\in& \Tm(\Gamma.A.B, B\{\langle {\rm p}(\Sigma(A,B)), \pi_1({\rm v}_{\Sigma(A,B)}) \rangle\}\{{\rm pair}_{A,B}\}) \\
&=& \Tm(\Gamma.A.B, B\{\langle {\rm p}(\Sigma(A,B)), \pi_1({\rm v}_{\Sigma(A,B)}) \rangle\}\{\langle {\rm p}(A.B), 
{\rm Pair}({\rm v}_{A.B},{\rm v}_B) \rangle\}) \\
&=& \Tm(\Gamma.A.B, B\{\langle {\rm p}(A.B), \pi_1({\rm v}_{\Sigma(A,B)})\{\langle {\rm p}(A.B), {\rm Pair}({\rm v}_{A.B},{\rm v}_B) \rangle\} \rangle\}) \\
&=& \Tm(\Gamma.A.B, B\{\langle {\rm p}(A.B), {\rm E}({\rm v}_{\Sigma(A,B)},{\rm v}_{A.B})\{\langle {\rm p}(A.B), {\rm Pair}({\rm v}_{A.B},{\rm v}_B) \rangle\} \rangle\}) \\
&=&  \Tm(\Gamma.A.B, B\{\langle {\rm p}(A.B), 
{\rm E}({\rm v}_{\Sigma(A,B)}\{\langle {\rm p}(A.B), {\rm Pair}({\rm v}_{A.B},{\rm v}_B) \rangle\} \rangle\},  \\
&& \qquad {\rm v}_{A.B}\{\langle {\rm p}(A.B), {\rm Pair}({\rm v}_{A.B},{\rm v}_B) \rangle\} \rangle.A.B\})) \\
&=&  \Tm(\Gamma.A.B, B\{\langle {\rm p}(A.B), 
{\rm E}({\rm Pair}({\rm v}_{A.B},{\rm v}_B), {\rm v}_{(A.B)\{\langle {\rm p}(A.B), {\rm Pair}({\rm v}_{A.B},{\rm v}_B) \rangle\}}))\\
&=&  \Tm(\Gamma.A.B, B\{\langle {\rm p}(A.B),{\rm v}_{A.B}\rangle\})\\
&=&  \Tm(\Gamma.A.B, B\{{\rm p}(B)\})\\
\end{eqnarray*}
We may let $K={\rm v}_B$. Define $\pi_{2,A,B}(P)=_{\rm def} {\rm E}(P,{\rm v}_B) \in
\Tm(\Gamma,B\{\langle 1_\Gamma, \pi_{1,A,B}(P)\rangle\})$. Then
$$
\pi_{2,A,B}({\rm Pair}_{A,B}(M,N))=
{\rm E}_{A,B,C}({\rm Pair}_{A,B}(M,N),{\rm v}_B) = 
{\rm v}_B\{\langle\langle 1_\Gamma, M\rangle_A,N\rangle_B\}
=  N, $$
and
\begin{eqnarray*}
\pi_{2,A,B}({\rm Pair}_{A,B}(M,N)) &\in&
\Tm(\Gamma,B\{\langle 1_\Gamma, {\rm E}({\rm Pair}_{A,B}(M,N),{\rm v}_{(A.B)\{\langle 1_\Gamma,{\rm Pair}(M,N)\rangle_{\Sigma(A,B)}\}})\rangle\}) \\
&=& \Tm(\Gamma,B\{\langle 1_\Gamma, M\rangle\})
\end{eqnarray*}
as required.
}
\end{example}

\begin{example}{\em For $A \in \Ty(\Gamma)$, $B \in \Ty(\Gamma.A)$, we define a mapping
$${\rm unpair}_{A,B} : \Gamma.\Sigma(A,B) \to \Gamma.A.B$$
First note that ${\rm p}_{\Sigma(A,B)}:  \Gamma.\Sigma(A,B) \to \Gamma$,
and 
\begin{eqnarray*}
{\rm v}_{\Sigma(A,B)} &\in& \Tm(\Gamma.\Sigma(A,B), \Sigma(A,B)\{{\rm p}_{\Sigma(A,B)}\}) \\
&=& \Tm(\Gamma.\Sigma(A,B), \Sigma(A\{{\rm p}_{\Sigma(A,B)}\},B\{{\rm p}_{\Sigma(A,B)}.A\})).
\end{eqnarray*}
Thus 
$$\pi_1({\rm v}_{\Sigma(A,B)}) \in 
\Tm(\Gamma.\Sigma(A,B), A\{{\rm p}_{\Sigma(A,B)}\})$$
and hence
$$\langle {\rm p}_{\Sigma(A,B)},  \pi_1({\rm v}_{\Sigma(A,B)}) \rangle :
\Gamma.\Sigma(A,B) \to \Gamma.A.$$
Further 
\begin{eqnarray*}\pi_{2}({\rm v}_{\Sigma(A,B)}) &\in&
\Tm(\Gamma.\Sigma(A,B),B\{{\rm p}_{\Sigma(A,B)}.A\}\{\langle 1_\Gamma, \pi_{1}({\rm v}_{\Sigma(A,B)})\rangle\}) \\
&=& \Tm(\Gamma.\Sigma(A,B),B\{\langle {\rm p}_{\Sigma(A,B)},  \pi_1({\rm v}_{\Sigma(A,B)}) \rangle\})
\end{eqnarray*}
so we define
$${\rm unpair}_{A,B}=_{\rm def} \langle\langle {\rm p}_{\Sigma(A,B)},  \pi_1({\rm v}_{\Sigma(A,B)}) \rangle, \pi_{2}({\rm v}_{\Sigma(A,B)})\rangle :
\Gamma.\Sigma(A,B) \to \Gamma.A.B$$
}
\end{example}

\medskip
A cwf {\em supports $+$-types} if for $A,B \in \Ty(\Gamma)$ there is
a type $A+B \in \Ty(\Gamma)$, and for each $M \in \Tm(\Gamma, A)$
there is ${\rm inl}_{A,B}(M) \in \Tm(\Gamma,A+B)$, and for
each $N \in \Tm(\Gamma, B)$
there is ${\rm inr}_{A,B}(N) \in \Tm(\Gamma,A+B)$.
For all $f: \Delta \rightarrow \Gamma$, these construction should satisfy
\begin{itemize}
\item[(+-subst)] $(A+B)\{f\} = A\{f\}+B\{f\}$
\item[(inl-subst)] ${\rm inl}_{A,B}(M)\{f\} = {\rm inl}_{A\{f\},B\{f\}}(M\{f\})$
\item[(inr-subst)] ${\rm inr}_{A,B}(N)\{f\} = {\rm inr}_{A\{f\},B\{f\}}(N\{f\})$
\end{itemize}

Furthermore for each $C \in \Ty(\Gamma.A+B)$, each $P \in \Tm(\Gamma,A+B)$ and each pair 
$$K_1 \in \Tm(\Gamma.A,C\{\langle {\rm p}A,{\rm inl}_{A\{{\rm p}A\},B\{{\rm p}A\}}({\rm v}_A)\rangle_{A+B}\})$$ and  
$$K_2 \in  \Tm(\Gamma.B,C\{\langle {\rm p}B,{\rm inr}_{A\{{\rm p}B\},B\{{\rm p}B\}}({\rm v}_B)\rangle_{A+B}\})$$
there is $${\rm D}_{A,B,C}(P,K_1,K_2) \in \Tm(\Gamma,C\{\langle 1_\Gamma,P\rangle\})$$
such that
\begin{itemize}
\item[(+-conv1)] ${\rm D}_{A,B,C}({\rm inl}_{A,B}(M),K_1,K_2) = 
K_1\{\langle 1_\Gamma,M \rangle_A\}$
\item[(+-conv2)] ${\rm D}_{A,B,C}({\rm inr}_{A,B}(N),K_1,K_2) = K_2\{\langle 1_\Gamma,N \rangle_B\}$.
\end{itemize}
Moreover for any $f: \Delta \rightarrow \Gamma$, the construction $D$ should satisfy
\begin{itemize}
\item[(D-subst)]  ${\rm D}_{A,B,C}(P,K_1,K_2)\{f\} =
{\rm D}_{A\{f\},B\{f\},C\{f.(A+B)\}}(P\{f\},K_1\{f.A\},K_2\{f.B\})$
\end{itemize}

\begin{lemma}\label{TimesL} Let ${\cal C}$ be a cwf that supports $\Sigma$. Suppose that 
$\Gamma \in {\cal C}$ and $A,B,C \in \Ty(\Gamma)$.
\begin{itemize}
\item[(a)] If $\Tm(\Gamma.A, B\{{\rm p}\})$ and $\Tm(\Gamma.A, C\{{\rm p}\})$ are inhabited, then 
so is $\Tm(\Gamma.A, (B \times C)\{{\rm p}\})$
\item[(b)] $\Tm(\Gamma.B \times C, B\{{\rm p}\})$ and  $\Tm(\Gamma.B \times C, C\{{\rm p}\})$
are inhabited.
\item[(c)] For any $f: \Delta \to \Gamma$,
$$(B \times C)\{f\} = (B\{f\} \times C\{f\}).$$
\end{itemize}
\end{lemma}
{\flushleft \bf Proof.} As for (a): suppose that 
$b \in \Tm(\Gamma.A, B\{{\rm p}(A)\})$ and 
$c \in \Tm(\Gamma.A, C\{{\rm p}(A)\})$.
Find an element in 
\begin{eqnarray*}
\Tm(\Gamma.A, (B \times C)\{{\rm p}A\}) &=& 
\Tm(\Gamma.A,   \Sigma(B, C\{{\rm p}B\})\{{\rm p}A\}) \\
&=& 
\Tm(\Gamma.A,   \Sigma(B\{{\rm p}A\}, C\{{\rm p}B\}\{{\rm p}(A).B\})).
\end{eqnarray*}

We have 
\begin{eqnarray*}
&&\Tm(\Gamma.A, C\{{\rm p}B\}\{{\rm p}(A).B\}\{\langle1_{\Gamma.A}, b\rangle_{B\{{\rm p}A\}}\})
\\
&&= \Tm(\Gamma.A, C\{{\rm p}B\}\{\langle{\rm p}(A){\rm p}(B\{{\rm p}A\}), {\rm v}_{B\{{\rm p}A\}}\rangle\}\{\langle1_{\Gamma.A}, b\rangle_{B\{{\rm p}A\}}\}) \\
&&= \Tm(\Gamma.A, C\{{\rm p}(B)\}\{\langle{\rm p}(A){\rm p}(B\{{\rm p}A\})\langle1_{\Gamma.A}, b\rangle, {\rm v}_{B\{{\rm p}A\}} \{\langle1_{\Gamma.A}, b\rangle\} \rangle\}) \\
&&= \Tm(\Gamma.A, C\{{\rm p}B\}\{\langle{\rm p}A, b \rangle\}) \\
&&= \Tm(\Gamma.A, C\{{\rm p}A\}).
\end{eqnarray*}
Thus $$c \in \Tm(\Gamma.A, C\{{\rm p}B\}\{{\rm p}(A).B\}\{\langle1_{\Gamma.A}, b\rangle_{B\{{\rm p}A\}}\}),$$
so
$${\rm Pair}(b,c) \in \Sigma(B\{{\rm p}A\}, C\{{\rm p}B\}\{{\rm p}(A).B\})$$
as required.

We prove (c) before (b): For $f: \Delta \to \Gamma$, we have
\begin{eqnarray*}
(B \times C)\{f\} &=& \Sigma(B, C\{{\rm p}B\})\{f\} \\
&=& \Sigma(B\{f\}, C\{{\rm p}B\}\{f.B\}) \\
&=& 
\Sigma(B\{f\}, C\{f\}\{{\rm p}(B\{f\})\}) \\
&=&  B\{f\} \times C\{f\} \\
\end{eqnarray*}

As for (b):  We have using (c):
\begin{eqnarray*}
{\rm v}_{B\times C} &\in& \Tm(\Gamma.B \times C, (B\times C)\{{\rm p}(B\times C)\}) \\
&=& \Tm(\Gamma.B \times C, B\{{\rm p}(B\times C)\} \times C\{{\rm p}(B\times C)\}) \\
&=& \Tm(\Gamma.B \times C, \Sigma(B\{{\rm p}(B\times C)\}, 
     C\{{\rm p}(B\times C)\}\{{\rm p}\{B\{{\rm p}(B\times C)\}\}\}))
\end{eqnarray*}
By (\ref{firstproj}), 
$$\pi_1({\rm v}_{B\times C}) \in 
\Tm(\Gamma.B \times C, B\{{\rm p}(B\times C)\}).$$
Thus $\Tm(\Gamma.B \times C, B\{{\rm p}(B\times C)\})$ is inhabited as required.
By (\ref{secondproj}), 
\begin{eqnarray*}
\pi_2({\rm v}_{B\times C}) &\in &
\Tm(\Gamma.B \times C, 
C\{{\rm p}(B\times C)\}\{{\rm p}\{B\{{\rm p}(B\times C)\}\}\}
\{\langle 1_{\Gamma.B \times C}, {\rm E}({\rm v}_{B\times C},{\rm v})\rangle\}) \\
&=& \Tm(\Gamma.B \times C, 
C\{{\rm p}(B\times C)\}).
\end{eqnarray*}
Thus also $\Tm(\Gamma.B \times C, 
C\{{\rm p}(B\times C)\})$ is inhabited as was to be shown.
$\qed$.

\begin{lemma} \label{PlusL} Let ${\cal C}$ be a cwf that supports $+$. Suppose that 
$\Gamma \in {\cal C}$ and $A,B,C \in \Ty(\Gamma)$.
\begin{itemize}
\item[(a)] If $\Tm(\Gamma.B, A\{{\rm p}\})$ and $\Tm(\Gamma.C, A\{{\rm p}\})$ are inhabited, then 
so is $\Tm(\Gamma.(B + C), A\{{\rm p}\})$.
\item[(b)] $\Tm(\Gamma.B, (B + C)\{{\rm p}\})$ and  $\Tm(\Gamma.C, (B + C)\{{\rm p}\})$
are inhabited.
\end{itemize}
\end{lemma}
{\flushleft \bf Proof.}
As for (a): Suppose $b \in \Tm(\Gamma.B, A\{{\rm p}B\})$ and 
$c \in \Tm(\Gamma.C, A\{{\rm p}C\})$. We have ${\rm p}(B+C): \Gamma.(B+C) \to \Gamma$, so 
$${\rm p}(B+C).B :\Gamma.(B+C).B\{{\rm p}(B+C)\}  \to \Gamma.B$$
and
$${\rm p}(B+C).C :\Gamma.(B+C).C\{{\rm p}(B+C)\}  \to \Gamma.C$$
and by substitution
$$b\{{\rm p}(B+C).B\} \in \Tm(\Gamma.(B+C).B\{{\rm p}(B+C)\}, A\{{\rm p}(B)({\rm p}(B+C).B)\})$$
and
$$c\{{\rm p}(B+C).C\}  \in \Tm(\Gamma.(B+C).C\{{\rm p}(B+C)\}, A\{{\rm p}(C)({\rm p}(B+C).C)\}) $$

Let $$U = A\{{\rm p}(B+C)\}\{{\rm p}((B+C)\{{\rm p}(B+C)\})\} \in \Ty(\Gamma.(B+C).(B+C)\{{\rm p}(B+C)\}).$$
Thus  
$${\rm v}_{(B+C)} \in \Tm(\Gamma.(B+C), (B+C)\{{\rm p}(B+C)\})
= \Tm(\Gamma.(B+C), B\{{\rm p}(B+C)\} + C\{{\rm p}(B+C)\})$$
For
\begin{eqnarray*}
K_1 &\in& \Tm(\Gamma.(B+C).B\{{\rm p}(B+C)\},U\{\langle {\rm p}(B\{{\rm p}(B+C)\}),{\rm inl}({\rm v}_{B\{{\rm p}(B+C)\}})\rangle\})\\
&=& \Tm(\Gamma.(B+C).B\{{\rm p}(B+C)\},A\{{\rm p}(B+C)\}\{{\rm p}(B\{{\rm p}(B+C)\})\}) \\
&=& \Tm(\Gamma.(B+C).B\{{\rm p}(B+C)\},A\{{\rm p}(B)({\rm p}(B+C).B)\}) 
\end{eqnarray*}
\begin{eqnarray*}
K_2 &\in& \Tm(\Gamma.(B+C).C\{{\rm p}(B+C)\},U\{\langle {\rm p}(C\{{\rm p}(B+C)\}),{\rm inl}({\rm v}_{C\{{\rm p}(B+C)\}})\rangle\})\\
&=& \Tm(\Gamma.(B+C).C\{{\rm p}(B+C)\},A\{{\rm p}(B+C)\}\{{\rm p}(C\{{\rm p}(B+C)\})\}) \\
&=& \Tm(\Gamma.(B+C).C\{{\rm p}(B+C)\},A\{{\rm p}(C)({\rm p}(B+C).C)\}) 
\end{eqnarray*}
there is $${\rm D}({\rm v}_{(B+C)},K_1,K_2) \in \Tm(\Gamma.(B+C),U\{\langle 1_\Gamma,{\rm v}_{(B+C)}\rangle\})=
\Tm(\Gamma.(B+C),A\{{\rm p}(B+C)\}).$$
Hence letting $K_1=b\{{\rm p}(B+C).B\}$ and $K_2=c\{{\rm p}(B+C).C\}$, we have
$${\rm D}({\rm v}_{(B+C)},b\{{\rm p}(B+C).B\},c\{{\rm p}(B+C).C\}) \in \Tm(\Gamma,(B+C),A\{{\rm p}(B+C)\})$$
as required.

As for (b): We have ${\rm v}_B \in {\rm Tm}(\Gamma.B, B\{{\rm p}B\})$
and  ${\rm v}_C \in {\rm Tm}(\Gamma.C, C\{{\rm p}C\})$. Thus
$${\rm inl}_{B\{{\rm p}B\}, C\{{\rm p}C\}}({\rm v}_B) \in 
{\rm Tm}(\Gamma.B, B\{{\rm p}B\})+C\{{\rm p}B\})$$
and
$${\rm inr}_{B\{{\rm p}C\}, C\{{\rm p}C\}}({\rm v}_C) \in 
{\rm Tm}(\Gamma.C, B\{{\rm p}B\})+C\{{\rm p}C\}).$$
Now $B\{{\rm p}B\})+C\{{\rm p}B\} =(B+C)\{{\rm p}B\}$ and
$B\{{\rm p}C\})+C\{{\rm p}C\} =(B+C)\{{\rm p}C\}$, so
${\rm Tm}(\Gamma.B,(B+C)\{{\rm p}B\})$  and
${\rm Tm}(\Gamma.C,(B+C)\{{\rm p}C\})$  are inhabited as required.

$\qed$

\medskip
A cwf {\em supports $\Pi$-types} if for $A \in \Ty(\Gamma)$ and $B \in \Ty(\Gamma.A)$ there is 
a type  $\Pi(A,B) \in \Ty(\Gamma)$, and moreover for every $P \in \Tm(\Gamma.A,B)$ there is an
element $\lambda_{A,B}(P) \in  \Tm(\Gamma, \Pi(A, B))$, and furthermore  for any 
$M \in  \Tm(\Gamma, \Pi(A,B))$ and any $N \in  \Tm(\Gamma,A)$ this is an element
${\rm App}_{A,B}(M,N) \in  \Tm(\Gamma, B\{\langle1_\Gamma, N\rangle_A\})$, such that the following equations hold for any
$f: \Theta \to \Gamma$:
\begin{itemize}
\item[($\beta$-conv)] ${\rm App}_{A,B}(\lambda_{A,\tau}(P),N) = P\{\langle1_\Gamma, N\rangle_A\}$,  
\item[($\Pi$-subst)] $\Pi(A,B)\{f\} = \Pi(A\{f\},B\{f.A\})$,
\item[($\lambda$-subst)] $\lambda_{A,B}(P)\{f\}= \lambda_{A\{f\},B\{f.A\}}(P\{f.A\})$,
\item[(${\rm App}$-subst)] ${\rm App}_{A,B}(M,N)\{f\} = {\rm App}_{A\{f\},B\{f.A\}}(M\{f\},N\{f\})$.
\end{itemize}

\medskip
For $A, B \in \Ty(\Gamma)$ we introduce the notation
$$(A \rightarrow B) =_{\rm def} \Pi(A, B\{{\rm p}(A)\}) \in \Ty(\Gamma)$$

\begin{lemma} \label{CCL} Let ${\cal C}$ be a cwf that supports $\Sigma$ and $\Pi$. Suppose that 
$\Gamma \in {\cal C}$ and $A,B,C \in \Ty(\Gamma)$.
\begin{itemize}
\item[(a)] $\Tm(\Gamma.A \times B, C\{{\rm p}\})$ is inhabited, if, and only if, 
$\Tm(\Gamma.A, (B \rightarrow C)\{{\rm p}\})$ is inhabited.
\item[(b)] For any $f: \Delta \to \Gamma$,
$$(B \rightarrow C)\{f\} = (B\{f\} \rightarrow C\{f\}).$$
\end{itemize}
\end{lemma}
{\flushleft \bf Proof.}  First we prove (b): For $f: \Delta \to \Gamma$, we have
\begin{eqnarray*}
(B \rightarrow C)\{f\} &=& \Pi(B, C\{{\rm p}B\})\{f\} \\
&=& \Pi(B\{f\}, C\{{\rm p}B\}\{f.B\}) \\
&=& 
\Pi(B\{f\}, C\{f\}\{{\rm p}(B\{f\})\}) \\
&=&  B\{f\} \rightarrow C\{f\} \\
\end{eqnarray*}

As for (a): Suppose first that $h \in \Tm(\Gamma.A \times B, C\{{\rm p}(A \times B)\})$.
To find an element $\lambda(P)$ of 
\begin{eqnarray*}
\Tm(\Gamma.A, (B \rightarrow C)\{{\rm p}A\}) &=&
\Tm(\Gamma.A, B\{{\rm p}A\} \rightarrow C\{{\rm p}A\}) \\
&=& \Tm(\Gamma.A, \Pi(B\{{\rm p}A\}, C\{{\rm p}A\}\{{\rm p}(B\{{\rm p}A\})\}))
\end{eqnarray*}
it is sufficient to find
$$P \in \Tm(\Gamma.A.B\{{\rm p}A\}, C\{{\rm p}A\}\{{\rm p}(B\{{\rm p}A\})\}).$$
We have 
$${\rm pair}_{A,B\{{\rm p}A\}}: \Gamma.A.B\{{\rm p}A\} \rightarrow \Gamma.\Sigma(A,B\{{\rm p}A\})=\Gamma.A \times B,$$
so
$$h\{{\rm pair}_{A,B\{{\rm p}A\}}\} \in
 \Tm(\Gamma.A.B\{{\rm p}A\}, C\{{\rm p}(A \times B)\}\{{\rm pair}_{A,B\{{\rm p}A\}}\}).$$
 Now
 $${\rm pair}_{A,B\{{\rm p}A\}} =_{\rm def} \langle {\rm p}(A.(B\{{\rm p}A\})), 
{\rm Pair}({\rm v}_{A.(B\{{\rm p}A\})},{\rm v}_{B\{{\rm p}A\}}) \rangle,$$
which entails that
$$C\{{\rm p}(A \times B)\}\{{\rm pair}_{A,B\{{\rm p}A\}}\} = 
C\{{\rm p}(A.(B\{{\rm p}A\}))\} = C\{{\rm p}A\}\{{\rm p}(B\{{\rm p}A\})\}).$$
Thus we can take $P= h\{{\rm pair}_{A,B\{{\rm p}A\}}\}$ to prove
that $\Tm(\Gamma.A, (B \rightarrow C)\{{\rm p}A\})$ is inhabited.

Conversely, suppose that 
$$k \in \Tm(\Gamma.A, (B \rightarrow C)\{{\rm p}A\})
= \Tm(\Gamma.A, \Pi(B\{{\rm p}A\}, C\{{\rm p}A\}\{{\rm p}(B\{{\rm p}A\})\})).$$
Need to find an element in
$$\Tm(\Gamma.A \times B, C\{{\rm p}(A \times B)\}).$$
As above we have
$$\pi_1({\rm v}_{A\times B}) \in 
\Tm(\Gamma.A \times B, A\{{\rm p}(A\times B)\}).$$
Hence 
$$\langle {\rm p}(A\times B), \pi_1({\rm v}_{A\times B}) \rangle:
\Gamma.A \times B \to \Gamma.A
$$
Writing $f$ for this map, we have 
\begin{eqnarray*}
k\{f\} &\in& \Tm(\Gamma.A \times B, \Pi(B\{{\rm p}A\}, C\{{\rm p}A\}\{{\rm p}(B\{{\rm p}A\})\})\{f\})) \\
&=& \Tm(\Gamma.A \times B, \Pi(B\{{\rm p}A\}\{f\}, C\{{\rm p}A\}\{{\rm p}(B\{{\rm p}A\})\}\{f.B\{{\rm p}A\}\})) \\
&=& \Tm(\Gamma.A \times B, \Pi(B\{{\rm p}(A\times B)\}, C\{{\rm p}A\}\{f\}\{{\rm p}\{B\{{\rm p}(A)f\}\}\})) \\
&=& \Tm(\Gamma.A \times B, \Pi(B\{{\rm p}(A\times B)\}, C\{{\rm p}(A\times B)\}\{{\rm p}\{B\{{\rm p}(A\times B)\}\}\})) \\
\end{eqnarray*}
Furthermore as above we have
$$\pi_2({\rm v}_{A\times B})\in \Tm(\Gamma.A \times B, 
B\{{\rm p}(A\times B)\}).$$
Thus
\begin{eqnarray*}
{\rm App}(k\{f\},\pi_2({\rm v}_{A\times B})) &\in& 
 \Tm(\Gamma.A \times B, C\{{\rm p}(A\times B)\}\{{\rm p}\{B\{{\rm p}(A\times B)\}\}\}
 \{\langle1_{\Gamma.A \times B}, \pi_2({\rm v}_{A\times B})\rangle\}) \\
 &=&  \Tm(\Gamma.A \times B, C\{{\rm p}(A\times B)\}).
  \end{eqnarray*}
 This proves that $\Tm(\Gamma.A \times B, C\{{\rm p}(A\times B)\})$ is inhabited as
 required. $\qed$

\begin{lemma} \label{forallL} Suppose that ${\cal C}$ is a cwf that supports $\Pi$. Let
$\Gamma \in {\cal C}$, $S, Q \in \Ty(\Gamma)$ and $R \in \Ty(\Gamma.S)$.
Then
$$\Tm(\Gamma.Q, \Pi(S,R)\{{\rm p}Q\})$$
is inhabited, if and only if,
$$\Tm(\Gamma.S.Q\{{\rm p}S\}, R\{{\rm p}(Q\{{\rm p}S\})\})$$
inhabited.
\end{lemma}
{\flushleft \bf Proof.}  Assume that  
$$f \in \Tm(\Gamma.Q, \Pi(S,R)\{{\rm p}Q\})= 
\Tm(\Gamma.Q, \Pi(S\{{\rm p}Q\},R\{{\rm p}(Q).S\})).$$
We have the map ${\rm p}(S).Q: \Gamma.S.Q\{{\rm p}S\} \to \Gamma.Q$ and substituting this into $f$ yields
\begin{eqnarray*}
f\{{\rm p}(S).Q\} &\in& \Tm(\Gamma.S.Q\{{\rm p}S\}, 
\Pi(S\{{\rm p}Q\},R\{{\rm p}(Q).S\})\{{\rm p}(S).Q\}) \\
&=& \Tm(\Gamma.S.Q\{{\rm p}S\}, 
\Pi(S\{{\rm p}Q\}\{{\rm p}(S).Q\},R\{{\rm p}(Q).S\}\{{\rm p}(S).Q.S\{{\rm p}Q\}\})) \\
&=& \Tm(\Gamma.S.Q\{{\rm p}S\}, 
\Pi(S\{{\rm p}S\}\{{\rm p}(Q\{{\rm p}S\})\},R\{{\rm p}(Q).S\}\{{\rm p}(S).Q.S\{{\rm p}Q\}\})) \\
\end{eqnarray*}
We have 
$${\rm v}_{S.Q\{{\rm p}S\}} \in \Tm(\Gamma.S.Q\{{\rm p}S\},S\{{\rm p}(S.Q\{{\rm p}S\})\})
= \Tm(\Gamma.S.Q\{{\rm p}S\},S\{{\rm p}S\}\{{\rm p}(Q\{{\rm p}S\})\}).$$
Thus 
\begin{eqnarray*}
&&{\rm App}(f\{{\rm p}(S).Q\}, {\rm v}_{S.Q\{{\rm p}S\}}) \in
\\
&=&
\Tm(\Gamma.S.Q\{{\rm p}S\},
R\{{\rm p}(Q).S\}\{{\rm p}(S).Q.S\{{\rm p}Q\}\}\{\langle 1_{\Gamma.S.Q\{{\rm p}S\}},{\rm v}_{S.Q\{{\rm p}S\}}\rangle\}) \\
&=&
\Tm(\Gamma.S.Q\{{\rm p}S\},
R\{\langle{\rm p}(Q) \circ {\rm p}(S).Q ,{\rm v}_{S.Q\{{\rm p}(S)\}}\rangle\}) \\
&=&
\Tm(\Gamma.S.Q\{{\rm p}S\},
R\{\langle{\rm p}(S) \circ {\rm p}\{Q\{{\rm p}S\}\} ,{\rm v}_{S}\{{\rm p}(Q\{{\rm p}S\})\rangle\}) \\
&=&
\Tm(\Gamma.S.Q\{{\rm p}S\},
R\{{\rm p}(Q\{{\rm p}S\}))
\end{eqnarray*}
as required.

As for the converse, suppose $g \in \Tm(\Gamma.S.Q\{{\rm p}S\},
R\{{\rm p}(Q\{{\rm p}S)\})$. Construct a map
$\Gamma.Q.S\{{\rm p}Q\} \to \Gamma.S.Q\{{\rm p}S\}$ as follows.
We have $${\rm p}(Q.S\{{\rm p}Q\}):\Gamma.Q.S\{{\rm p}Q\} \to \Gamma$$
and
$${\rm v}_{S\{{\rm p}Q\}} \in \Tm(\Gamma.Q.S\{{\rm p}Q\}, S\{{\rm p}(Q.S\{{\rm p}Q\})\})$$
so 
$$\langle {\rm p}(Q.S\{{\rm p}Q\}), {\rm v}_{S\{{\rm p}Q\}} \rangle : \Gamma.Q.S\{{\rm p}Q\} \to \Gamma.S.$$
Now
\begin{eqnarray*}
{\rm v}_{Q.S\{{\rm p}Q\}} &\in & \Tm(\Gamma.Q.S\{{\rm p}Q\}, Q\{{\rm p}S\}\{\langle {\rm p}(Q.S\{{\rm p}Q\}), {\rm v}_{S\{{\rm p}Q\}} \rangle\}) \\
&=& \Tm(\Gamma.Q.S\{{\rm p}Q\}, Q\{{\rm p}(Q.S\{{\rm p}Q\})\}).
\end{eqnarray*}
Hence
$$\langle \langle {\rm p}(Q.S\{{\rm p}Q\}), {\rm v}_{S\{{\rm p}Q\}} \rangle, {\rm v}_{Q.S\{{\rm p}Q\}} \rangle:
\Gamma.Q.S\{{\rm p}Q\} \to \Gamma.S.Q\{{\rm p}S\}.
$$
Substituting this in to $g$,
\begin{eqnarray*}
&&g\{\langle \langle {\rm p}(Q.S\{{\rm p}Q\}), {\rm v}_{S\{{\rm p}Q\}} \rangle, {\rm v}_{Q.S\{{\rm p}Q\}} \rangle\} \\
&& \in \Tm(\Gamma.Q.S\{{\rm p}Q\} ,
R\{{\rm p}(Q\{{\rm p}S)\}\{\langle \langle {\rm p}(Q.S\{{\rm p}Q\}), {\rm v}_{S\{{\rm p}Q\}} \rangle, {\rm v}_{Q.S\{{\rm p}Q\}} \rangle\}) \\
&& = \Tm(\Gamma.Q.S\{{\rm p}Q\} ,
R\{\langle {\rm p}(Q.S\{{\rm p}Q\}), {\rm v}_{S\{{\rm p}Q\}} \rangle\}) \\
&&= \Tm(\Gamma.Q.S\{{\rm p}Q\} ,
R\{{\rm p}(Q).S\})
\end{eqnarray*}
Thus by lambda abstraction
\begin{eqnarray*}
g\{\langle \langle {\rm p}(Q.S\{{\rm p}Q\}), {\rm v}_{S\{{\rm p}Q\}} \rangle, {\rm v}_{Q.S\{{\rm p}Q\}} \rangle\} &\in& \Tm(\Gamma.Q, \Pi(S\{{\rm p}Q\}, R\{{\rm p}(Q).S\})) \\
&=& \Tm(\Gamma.Q, \Pi(S, R)\{{\rm p}Q\})
\end{eqnarray*}
as required.
$\qed$

\begin{lemma} \label{existsL} Suppose that ${\cal C}$ is a cwf that supports $\Sigma$. Let
$\Gamma \in {\cal C}$, $S, Q \in \Ty(\Gamma)$ and $R \in \Ty(\Gamma.S)$.
Then
$$\Tm(\Gamma.\Sigma(S,R), Q\{{\rm p}(\Sigma(S,R))\})$$
is inhabited, if and only if,
$$\Tm(\Gamma.S.R, Q\{{\rm p}(S.R)\})$$
inhabited.
\end{lemma}
{\flushleft \bf Proof.} Suppose $f \in \Tm(\Gamma.\Sigma(S,R), Q\{{\rm p}(\Sigma(S,R))\})$.
Then
\begin{eqnarray*}
f\{{\rm pair}_{S,R}\} &\in  &\Tm(\Gamma.S.R, Q\{{\rm p}(\Sigma(S,R))\}\{{\rm pair}_{S,R}\})
\\
&= &\Tm(\Gamma.S.R, Q\{{\rm p}(S.R)\})
\end{eqnarray*}
as required.

Conversely, suppose $g  \in \Tm(\Gamma.S.R, Q\{{\rm p}(S.R)\})$. Then
\begin{eqnarray*}
g\{{\rm unpair}_{S,R}\}  &\in &\Tm(\Gamma.\Sigma(S,R), Q\{{\rm p}(S.R)\}\{{\rm unpair}_{S,R}\}) \\
&= & \Tm(\Gamma.\Sigma(S,R), Q\{{\rm p}(S.R)\}\{
\langle\langle {\rm p}_{\Sigma(S,R)},  \pi_1({\rm v}_{\Sigma(S,R)}) \rangle, \pi_{2}({\rm v}_{\Sigma(S,R)})\rangle\}) \\
&= & \Tm(\Gamma.\Sigma(S,R), Q\{{\rm p}_{\Sigma(S,R)}\}) \\
\end{eqnarray*}
as required.
$\qed$

\begin{customthm}{\ref{MLTThdoc}}  If ${\cal C}$ is a cwf which admits the type constructions $\Sigma, \Pi, +, {\rm N}_0$ and ${\rm N}_1$ then $({\cal C},{\rm Pr}_{\cal C})$ is a cwf with first-order doctrine.
\end{customthm}
{\flushleft \bf Proof.} Let the logical operations be defined according to the propositions-as-types principle:
\begin{itemize}
\item $\top_\Gamma = \nattype_{1,\Gamma}$,
\item $\bot_\Gamma = \nattype_{0,\Gamma}$,
\item $P \land Q = P \times Q$,
\item $P \lor Q = P + Q$,
\item $P \rightarrow Q = P \rightarrow Q$,
\item $\forall_S(R) = \Pi(S,R)$,
\item $\exists_S(R) = \Sigma(S,R)$.
\end{itemize}

Now the properties of the first-order doctrine follows by Lemmas
\ref{N0N1L},
\ref{TimesL},
\ref{PlusL},
 \ref{CCL},
\ref{forallL},
and
\ref{existsL}, and the substitution properties of the type constructions 
$\nattype_k$, $+$, $\Pi$ and $\Sigma$. $\qed$
}
\end{document}